\documentclass[utf8,latin1]{article}
\usepackage[T1]{fontenc}
\usepackage{lmodern}
\usepackage[utf8,latin1]{inputenc}
\usepackage{latexsym}                 
\usepackage{color}                 
\usepackage{amssymb}
\usepackage{amsmath}
\usepackage{amsthm}
\usepackage{graphicx}
\usepackage{hyperref}
\usepackage{verbatim}
\usepackage{mathptmx}      
\usepackage{bm}
\usepackage{matlab-prettifier}
\usepackage[section]{placeins}
\newcommand{\R}{\ensuremath{\mathbb{R}}}

\def\fa{\hbox{ for all }}

\def\Xn{X_n=\{x_1,\ldots,x_n\}}

\def\bod{\mathbf{d}}	
\def\boN{\mathbf{N}}	
\def\boz{\mathbf{z}}
\def\calh{{\cal H}}
\def\calO{\mathcal{O}}
\def\biglf{\par\bigskip\noindent} 
\def\eref#1{(\ref{#1}%
)}
\def\RSref#1{\ref{#1}%
}
\def\RSlabel#1{\label{#1}%
}
\def\RScite#1{\cite{#1}%
}
\newcommand{\bql}[1]{%
\begin{equation}\label{#1}%
}
\def\eq{\end{equation}}

\begin{document}

\begin{center}
  {\bf \Large Greedy Adaptive Local Recovery \\~\\ of Functions 
   in Sobolev Spaces}
~\\~\\
  {\bf Robert Schaback}\biglf 
Institut f\"ur Numerische und Angewandte Mathematik,\\
Georg-August-Universit\"at G\"ottingen, Lotzestra\ss{}e 16--18,\\
schaback@math.uni-goettingen.de
  \biglf
Version of \today
\end{center}
{\bf Abstract}\\
There are many ways to upsample functions from multivariate scattered
data locally,
using only a few neighbouring data points of the evaluation point.
The position and number
of the actually used data
points is not trivial, and many cases like Moving Least Squares
require point selections that guarantee local recovery of polynomials
up to a specified order. This paper suggests a kernel-based
greedy local algorithm
for point selection that has no such constraints. It
realizes the optimal $L_\infty$
convergence rates in Sobolev spaces using the minimal number of points
necessary for that purpose. On the downside, it does not care
for smoothness, relying on fast $L_\infty$ convergence to a smooth function.
The algorithm ignores near-duplicate points
automatically and works for quite irregularly distributed point sets
by proper selection of points. Its computational
complexity is constant for each evaluation point,
being dependent only on the Sobolev space parameters. 
Various numerical examples are provided.
As a byproduct, it turns out that the well-known instability of
global
kernel-based interpolation in the standard basis of kernel translates
arises already locally, independent of global kernel matrices
and small separation distances.
\biglf
\biglf
\section{Introduction}\RSlabel{SecIntro}
Throughout, we shall work on a bounded domain $\Omega\subset\R^d$
and consider recoveries of functions $f$ on $\Omega$
using values $f(x_i)$
on scattered locations $x_i\in X_N=\{x_1,\ldots,x_N\}
\subset\Omega\subset\R^d$. 
We focus on large point sets $X_N$, but avoid working with the whole set.
Instead, we calculate the recovery at each
point $z\in\Omega$ separately, using only the points
from a set $X(z)\subset X_N$, and with a fixed computational complexity
for each $z$. This strategy is well-known from Moving Least Squares
\cite{lancaster-salkauskas:1981-1,
  farwig:1986-1, bos-salkauskas:1989-1, levin:1998-1,wendland:2000-1,
  armentano:2001-1, fasshauer:2003, mirzaei-et-al:2012-1},
Shepard-type \cite{shepard:1968-1,franke-nielson:1991-1,dellaccio-tomaso:2016-1}
and Partition-of-Unity
\cite{babuska-melenk:1997-1,wendland:2002-1,%
heryudono-et-al:2016-1,larsson-et-al:2017-1,cavoretto:2021-1}
techniques, but here we ignore any tricks to ensure smoothness properties.
We focus on optimal convergence rates and minimal point sets
instead, and confine ourselves
to recoveries in spaces $\calh$ like Sobolev spaces, using
translates of the kernel $K$ that is reproducing in $\calh$.
\biglf
It is well-known that using the full set $X_N$
gives the least possible error, but how much is lost when using only the points
in $X(z)$ for recovery of $f(z)$?
Clearly, one loses smoothness and has a larger error while gaining
much better computational complexity. This tradeoff is the basic question here.
The proposed method shows optimal convergence rates
in $L_\infty$ at $\calO(1)$ computational complexity  at each evaluation point,
using the minimal possible number of points for such rates. 
\biglf
Our point selection strategy will not be based on geometry properties
like guaranteeing that sets are in general position with respect to
polynomials. This is a criterion that is unstable under variation of data
points. Instead, we select points by minimizing a continuous function
connected to the error, aiming at error-minimal selections
like in \RScite{davydov-schaback:2018-2}.
\biglf
The underlying theory is kernel-based
interpolation, as covered by the books \cite{buhmann:2003-1,wendland:2005-1,
  fasshauer-mccourt:2015-1}. Section \RSref{SecRec} provides the basic
notations, including the
{\em Power Function}, i.e.
the norm of the error functional in the kernel-based space
considered. As a prerequisite for our algorithm,
Section \RSref{SecNewBas} deals with the very useful {\em Newton basis}.
The algorithm follows in Section
\RSref{SecPoiSel}, as a stepwise adaptive greedy minimization of
the Power Function in terms of the Newton basis.
\biglf
Then there are various numerical examples. Sections \RSref{SecSingExa}
and \RSref{SecGlobExa} focus on the error locally and globally,
while Section \RSref{SecFunRep}
recovers functions and shows the grades of discontinuity of the recoveries.
Then 
Section \RSref{SecConv} demonstrates that the optimal convergence rates
are actually attained. The unexpected instabilities for large smoothness
parameters $m$ are explained in Section
\RSref{SecStab}, followed by conclusions and open problems
in Section \RSref{SecCOP}. 
\section{Recovery}\RSlabel{SecRec}
For a very large set $X_N=\{x_1,\ldots,x_N\}$ and a point $z\notin X_N$
we pick a set of $n$ points near $z$,
renumbered as $\Xn$.
Recovery of $f(z)$ from values at points of $\Xn$ can be written as 
$$
s_{f,X_n}(z)=\displaystyle{\sum_{j=1}^n u_j(z)f(x_j)   },\;f \in\calh, \,
  z\in\Omega 
$$
where the values $u_j(z)$ are just real numbers, because we do not vary $z$.
Lagrange conditions
$u_j(x_k)=\delta_{jk},\;1\leq j,k\leq n$ are not required.
The error functional is
$$
\epsilon_{X_n,z}\;:\; f\mapsto f(z)- s_{f,X_n}(z),
$$
written in terms of point evaluation functionals as
$$
\epsilon_{X_n,z}=\delta_z- \displaystyle{\sum_{j=1}^n u_j(z)\delta_{x_j}}.
$$
Its norm in $\calh$ is known as the {\em Power Function} $P_{X_n}(z)$.
By the standard dual representation
$$
K(x,y)=(\delta_x,\delta_y)_{\calh^*} \fa x,y\in \Omega,
$$
the squared norm is
$$
\|\epsilon_{X_n,z}\|^2_{\calh^*}=P_{X_n}^2(z)
=K(z,z)- 2\displaystyle{\sum_{j=1}^n u_j(z)K(z,x_j)}+
\displaystyle{\sum_{j,k=1}^n u_j(z)u_k(z)K(x_k,x_j)}.
$$
This makes the Power Function computable.
Due to the alternative definition
$$
P_{X_n}(z)=\displaystyle{\sup\{
  f(z) \;:\; \|f\|_\calh \leq 1, \;f(X_n)=\{0\}\}   }, 
$$
the Power Function decreases at all $z$ when the point set $X_n$ is enlarged. 
Details are in basic texts on kernels, e.g.
\cite{schaback:1997-3,buhmann:2003-1,wendland:2005-1,
  fasshauer-mccourt:2015-1}.
\biglf
Using the Power Function,
the error of the recovery of $f(z)$ by $s_{f,X_n}(z)$
at a point $z$ has the optimal bound 
\bql{eqerrbnd}
|f(z)-s_{f,X_n}(z)|\leq P_{X_n}(z)\|f\|_\calh \fa f\in\calh
\eq
in $\calh$, but we want to get away with fewer points forming a subset
$X(z)\subset X_n$. Therefore, the goal is to find a subset $X(z)$ of
$X_n\subset X_N$ such that
the difference between $P_{X(z)}(z)$ and its lower bound $P_{X_n}(z)$
is small. We tacitly assume that $N$ is too large to let
the even lower bound $P_{X_N}(z)$
be computable.
\biglf
Note that the error bound \eref{eqerrbnd} is {\em local}
or {\em pointwise}, and we shall use it after selection as
\bql{eqerrboundz}
|f(z)-s_{f,X(z)}(z)|\leq P_{X(z)}(z)\|f\|_\calh \fa f\in\calh.
\eq
We shall always know the value of $P_{X(z)}(z)$, and therefore we have an
error bound that extends to an $L_\infty$ error bound when varying $z$, because
$P^2_{X(z)}(z)\leq P^2_{\emptyset}(z)=K(z,z)$. This allows full control
of the $L_\infty$ error up to the unknown factor $\|f\|_\calh$. It is an old
unsolved problem to provide upper bounds on it. 
\biglf
We propose a greedy method to select the set $X(z)\subset X_n$ for fixed
$z$ by stepwise minmization of the Power Function as a
function of $x_1,\,x_2,\ldots,x_n$. This is similar to,
but different from the $P$-greedy point selection in
\RScite{DeMarchi-et-al:2005-1}. There, $x_{n+1}$
is picked as the argmax of $P_{X_n}^2(z)$
as a function of $z$. The paper 
\RScite{santin-haasdonk:2017-1} showed that the $P$-greedy
point selection method is
asymptotically optimal with respect to convergence rates.
\biglf
Here, we proceed differently. The point $z$ is fixed, and we select
$x_1,x_2,\ldots $ sequentially to produce smallest possible values
of $P_{\{x_1,\ldots,x_j\}}^2(z)$ for increasing $j$ until the error is small
enough. We can stop this at some bound on $j$ or by a lower threshold
on $P_{\{x_1,\ldots,x_j\}}^2(z)$. At any time, we have the error bound
\eref{eqerrboundz}
where $P_{X(z)}(z)$ is explicitly known.
\biglf
Of course, the final recovery $s_{f,X(z)}(z)$ will be a discontinuous function
of $z$. But if $P_{X(z)}(z)$ is small and controllable,
$s_{f,X(z)}(z)$is close to $f$ in
the $L_\infty$ norm, and its deviation from $f$ in $L_\infty$ may be as small
as $\|f-s_{f,X_N}\|_\infty$ if we can manage to let
$P_{X(z)}(z)$ for $X(z)\subset X_N$ behave like $\|f-s_{f,X_N}\|_\infty$.
\biglf
Here is an illustration for Sobolev spaces $W_2^m(\R^d)$ with $m>d/2$.
If sets $X_N$ have {\em fill distances}
$$
h(X_N,\Omega)=\displaystyle{ \sup_{z\in\Omega}\inf_{x_j\in X_N}\|z-x_j\|_2}
$$
and {\em separation distances}
$$
\sigma(X_N)=\displaystyle{\dfrac{1}{2}\inf_{x_j\neq x_k \in X_N}\|x_k-x_j\|_2}
$$
with {\em asymptotic regularity}
$$
0<c_0 \sigma(X_N) \leq h(X_N,\Omega)\leq C_0 \sigma(X_N),
$$
then by e.g. \RScite{wendland:2005-1}
$$
\|f-s_{f,X_N}\|_\infty
\leq \|P_{X_n}\|_\infty\|f\|_{W_2^m(\R^d)}
\leq C_1h^{m-d/2}(X_N,\Omega)\|f\|_{W_2^m(\R^d)}.
$$
Therefore our goal must be to ensure 
$$
P_{X(z)}(z)\leq C_2 h^{m-d/2}
$$
for $X(z)\subset X_N$ and all $z\in\Omega$. This is possible
for small $X(z)\neq X_N$,
and the examples in later sections will show how. We shall not need asymptotic
regularity for that, and we can get away with the minimal number of selected
points for getting $L_\infty$ convergence like $\calO(h^{m-d/2})$.
\biglf
The Power Function value $P_{X_n}(z)$ can be seen as a distance from $z$
to the set $X_n$, and therefore the above algorithm is a greedy method to find
nearest neighbours in the kernel metric. More explicitly, $P_{X_n}(z)$
is the Euclidean distance of the point evaluation functional $\delta_z$
to the space spanned by the functionals $\delta_{x_j},\,1\leq j\leq n$
in the dual $\calh^*$ of the Hilbert space $\calh$. The problem then is to
find a few functionals that already give a small distance to the space spanned
by them. This, in turn, is
a case of {\em $n$-term approximation} \cite{temlyakov:1998-1, temlyakov:1999-1}
that has good solutions by greedy
methods and a connection to sparsity techniques
\cite{cohen-et-al:2009-1}.
The Newton basis approach below will be a special implementation
adapted to kernel-based spaces.
\section{Newton Basis}\RSlabel{SecNewBas}
To study the variation of $P_{\{x_1,\ldots,x_j\}}^2(z)$
as a function of $x_j$, we use the Newton basis representation
dating back to \RScite{mueller-schaback:2009-1}, written here
via recursive kernels
\RScite{mouattamid-schaback:2009-1}.
The kernel recursion for points $x_1,x_2,\ldots$ is
$$
\begin{array}{rcl}
K_1(x,y)&:=& K(x,y)\\
K_{j+1}(x,y)&:=& K_j(x,y)-\dfrac{K_j(x,x_j)K_j(x_j,y)}{K_j(x_j,x_j)},\;j\geq
1,\;x,\,y\in\Omega. 
\end{array}
$$
It is easy to show by induction that
$$
K_{j+1}(x_k,y)=0=K_{j+1}(x,x_k)
$$
holds for all $x,\,y$ and all $1\leq k\leq j$.
The Newton basis function $N_j$  then is
$$
N_j(x)=\dfrac{K_j(x,x_j)}{\sqrt{K_j(x_j,x_j)}},\;j\geq 1,
$$
satisfying
$$
N_j(x_k)=0,\;1\leq k<j \hbox{ and } N_j(x_j)=\sqrt{K_j(x_j,x_j)},\,j\geq 1.
$$
Now the kernel recursion takes the form 
$$
K_{j+1}(x,y)=K_j(x,y)-N_j(x)N_j(y)=K(x,y)-
\displaystyle{\sum_{m=1}^j N_m(x)N_m(y), \,j\geq 1 } 
$$
cancelling the denominators, with the final residual
$$
K_{j+1}(x,y)=K(x,y)-\displaystyle{\sum_{m=1}^j N_m(x)N_m(y)  } ,\;x,y\in \Omega, \;j\geq 0.
$$
Inserting points for $x$ and $y$, we see that
$$
K(x_i,x_j)=\displaystyle{\sum_{m=1}^{\min(i,j)} N_m(x_i)N_m(x_j)  }\;1\leq i,j\leq n
$$
is a
Cholesky factorization of the positive definite symmetric kernel matrix.
It is in the background, but we prefer to work in terms of functions, not
matrices.  
The Newton basis recursion is
\bql{eqNNNKKK}
\begin{array}{rcl}
  N_j(x_j)^2&=&K_j(x_j,x_j),\\
  N_j(x)N_j(x_j)&=& K(x,x_j)-\displaystyle{\sum_{m=1}^{j-1} N_m(x)N_m(x_j),
  \;j\geq 1}.
\end{array} 
\eq

\section{Greedy Point Selection Algorithm}\RSlabel{SecPoiSel}
We want an  optimal point selection for recovery of $f(z)$ from values
$f(x_j),\;1\leq j\leq n$ by minimizing the Power Function. 
The first chosen point $x_1$ should minimize
$$
P^2_{x_1}(z) =K(z,z)-N^2_{1}(z) = K(z,z)-\dfrac{K(z,x_1)^2}{K(x_1,x_1)},
$$
as a function of $x_1$, via
$$
x_1=\arg\max_{x_j\in X_n}
\dfrac{K(z,x_j)^2}{K(x_j,x_j)}.
$$
By positive semidefiniteness of kernel matrices,
$$
K(z,x_j)^2\leq K(z,z)K(x_j,x_j)
$$
holds and leads to the choice of $x_j$ if $z=x_j$.  Then the process can be
stopped.
\biglf
If $x_1,\ldots,x_{j-1}$ are determined and $x_j$ is still not selected, the
Newton basis functions $N_k$ are determined for $k<j$, and we can assume
vectors and matrices
$$
\begin{array}{rcl}
  \boz&:=&(K_j(z,x_k),\,1\leq k\leq n)\in\R^n,\\
\bod&:=&(K_j(x_k,x_k),\,1\leq k\leq n))\in\R^n,\\
\boN&:=& (N_i(x_k),\,1\leq i<j,\,1\leq k\leq n).
\end{array} 
$$
The Power
Function at $z$ for varying  $x_j$ is
$$
P^2_{x_1,\ldots,x_j}(z)=K(z,z)-\sum_{k=1}^{j-1}N_k(z)^2-N_{j}(z)^2
=K(z,z)-\sum_{k=1}^{j-1}N_k(z)^2-\dfrac{K_j(z,x_j)^2}{K_j(x_j,x_j)}
$$
and the next point should be
\bql{eqdecision}
x_j=\arg\max_{x_k\in X_n, k\geq j}
\dfrac{K_j(z,x_k)^2}{K_j(x_k,x_k)}
=\arg\max_{j\leq k\leq n}
\dfrac{\boz_k^2}{\bod_k}.
\eq
The maximum value is $N_j^2(z)$ and we can construct $N_j$ by
the recursion \eref{eqNNNKKK} to get a new column in $\boN$.
For upgrading our vectors, we use 
$$
\begin{array}{rcl}
  K_{j+1}(z,x_k)
  &=&K_{j}(z,x_k)-N_{j}(z)N_{j}(x_k)\\
  K_{j+1}(x_k,x_k)
  &=&K_{j}(x_k,x_k)-  N_{j}^2(x_k)\\
  P_{j+1}^2(z)&=& P_{j}^2(z) -N_j^2(z)
\end{array}
$$
in terms of the Newton basis. Here, we still need to fix the sign of $N_j(z)$ to
match the sign of $K_j(z,x_j)=d_j$.
\biglf
The Lagrange coefficients $L_j(z)$ must recover the $N_k(z)$ exactly. They solve
the triangular system
\bql{eqNewtonrep}
\displaystyle{N_m(z)=\sum_{k=1}^jL_k(z)N_m(x_k),\;1\leq m\leq j}
\eq
that is cheaply solvable. Then one can monitor the 
Lebesgue constants
$$
\displaystyle{\sum_{k=1}^j |L_j(z)|   }
$$
that control the $L_\infty$ {\em evaluation stability}
\RScite{noorizadeghan-schaback:2024-1}
of the solution via
$$
|s_{f,X(z)}(z)-s_{g,X(z)}(z)|\leq \sum_{k=1}^j |L_j(z)||f(x_j)-g(x_j)|
\leq \sum_{k=1}^j |L_j(z)|\|f-g\|_\infty. 
$$
If $k$ steps are executed, storage grows like $\calO(kn)$, while
computational complexity is $\calO(k^2n)$. The next section will show that
one can keep $k$ much smaller than $n$ in most cases. In the 2D examples
below we get away with $n=5k$, letting the complexity for a single
recovery be $\calO(5k^3)$. 
Matrix conditions arise here only for the triangular matrices with entries
$N_m(x_k)$ used in \eref{eqNewtonrep}.
This will be relevant for the stability arguments of Section
\RSref{SecStab}. 
\section{Single Point Examples}\RSlabel{SecSingExa}
As described in Section \RSref{SecRec}, we consider very large sets
$X_N \subset \R^d$.
For each fixed point $z\in\R^d$ we have a subset $X_n$ of $X_N$
depending on $z$, taking possibly
nearest neighbours of $z$, 
to start the local greedy point
selection on $X_n$ to end up with an even smaller set
$X(z)\subseteq X_n$. Because everything is local,
we can set $z$ to be the origin and replace $X_n$ by $X_n-\{z\}$.
\biglf
Here, for illustration,
we work with $n=100$ random points
in $[-1,+1]^2$ depicted in the top left of Figure \RSref{FigLocalm3},
and used  in 
later cases as well. These are unduly many,
but the goal of this section
is to compare selections of $k<<n$ 
points with selection of all $n$ points. An even larger set $X_N$
may be in the background, on a larger domain. 
The greedy algorithm is run up to $n$
points to show how the error behaves for
small $k$ when compared to $k=n$. We want to find $k$ such that selecting $k$
points is not much different from selecting $n$ points.  
\biglf
Figure \RSref{FigLocalm3} shows results for the kernel
generating Sobolev space $W_2^3(\R^2)$.
The decay of the squared Power Function is plotted in the centre of the top
row. Large $n$ do not pay off, error-wise, because the curve flattens
dramatically. 
To compare with polynomial exactness orders $q=1,2,3,\ldots$  needing
$Q={q-1+d\choose d}=1,3,6,10,\ldots$ points in $\R^2$,
these cases are marked in red circles. The match between
the optimal rate $m-d/2$ in Sobolev space $W_2^m(\R^d)$ with polynomial
exactness of order $q$ is additionally marked here and later
with a red cross. Now $m=3$ leads to order $q=2$
and $Q=6$ points. This is our goal for $m=3$.
\biglf
In the plots of this section,
the third plot in the first
row shows the Lebesgue constants, the cases for
polynomial orders marked in red circles
again. The blue circles are the Lebesgue constants for using polynomials of the
appropriate orders at the selected points.
It can be expected that selecting the $Q$ best points is
sufficient for getting the optimal convergence rate. Taking more points
does not decrease the local error  substantially, see the top centre plot.
\biglf
The lower plots show the point selections for different numbers of points
belonging to different orders $q$. The central evaluation
point $z=(0,0)$ is marked with
a red cross. Note that the
selection does not take all nearest points.
\biglf
The low-regularity case $m=1.5$ is in Figure \RSref{FigLocalm1p5}, while
$m=6$ is in Figure \RSref{FigLocalm6}. The first case should not use more than
3 points. The last case is run at a scale of
$c=0.1$, taking $K(\|x-y\|_2/c)$, and using 100 regular data locations.
Here, taking 21 points is enough.  Polynomial Lebesgue constants
run out of hand, because the selected
points are not in general position wrt. polynomials. Due to regularity of the
point distribution, the greedy technique has several choices at various steps,
and therefore there are plenty of equivalent point selections.  
For scale $c=1$, the kernel matrix on 100 points exceeds the condition limit of $1.e14$.
But if one stops when the squared Power Function is below 1.e-8, one gets away
with 9 points, see Figure \RSref{FigLocalm6scal1}, without condition
problems.
\biglf
If the point is in the corner $(-1,+1)$, Figure \RSref{FigLocalm3corner} shows
the results for $m=3$ at scale 1.0, like Figure \RSref{FigLocalm3}. The
achievable squared Power Function now is about 0.004, while 
the central case had about 0.001. Using more than about 10 points does not help,
while roughly 4 were sufficient in the central situation. 
  \def\RSh{3.5cm}
  \def\RSw{3.5cm}
  \begin{figure}[hbtp]
    \begin{center}
 \includegraphics[width=\RSw,height=\RSh]{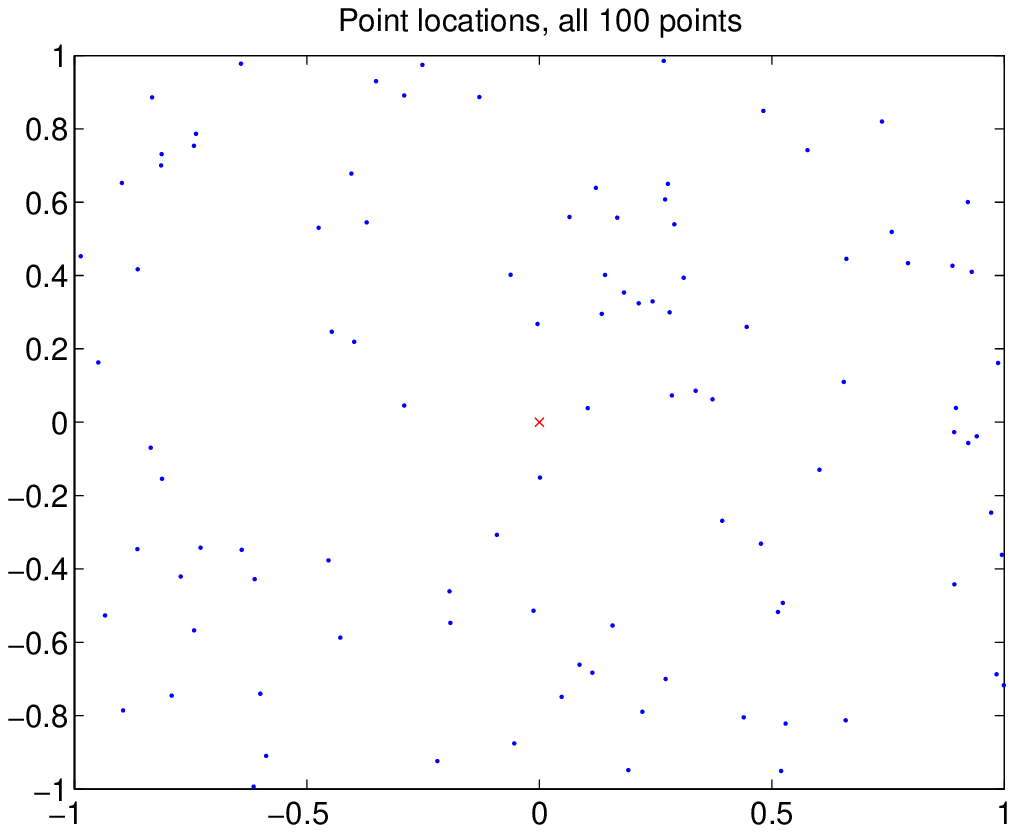}
 \includegraphics[width=\RSw,height=\RSh]{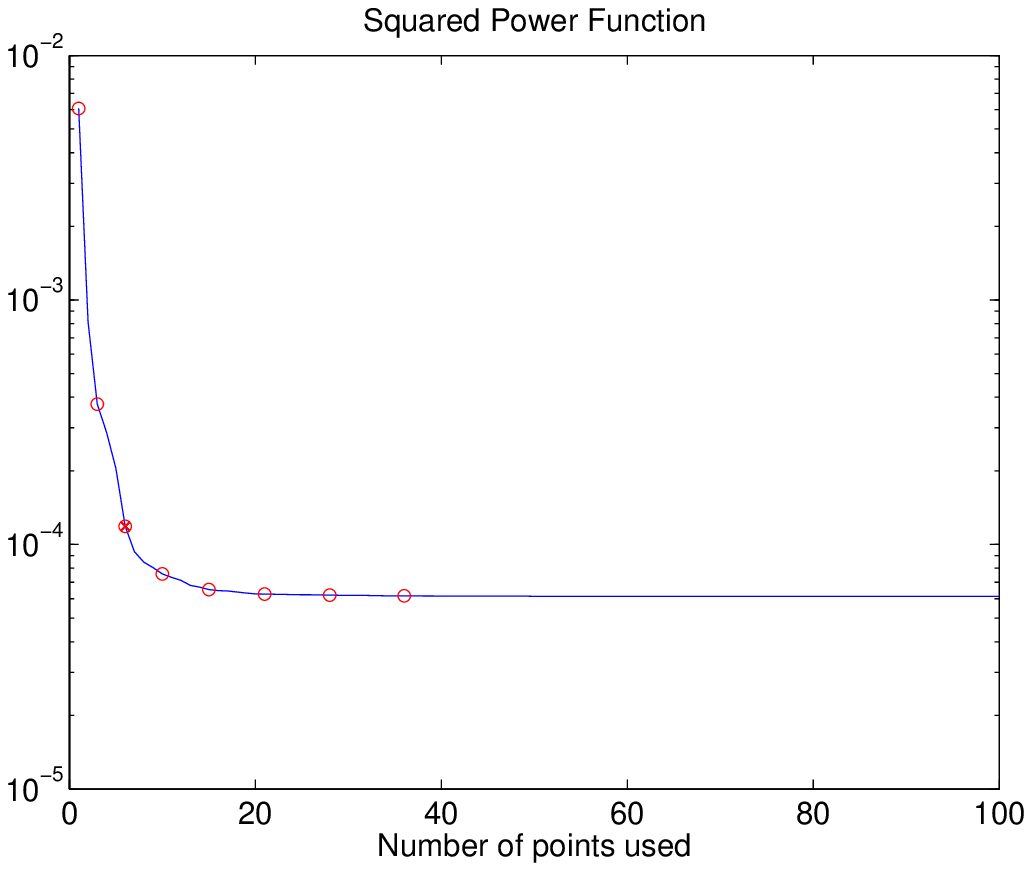}
 \includegraphics[width=\RSw,height=\RSh]{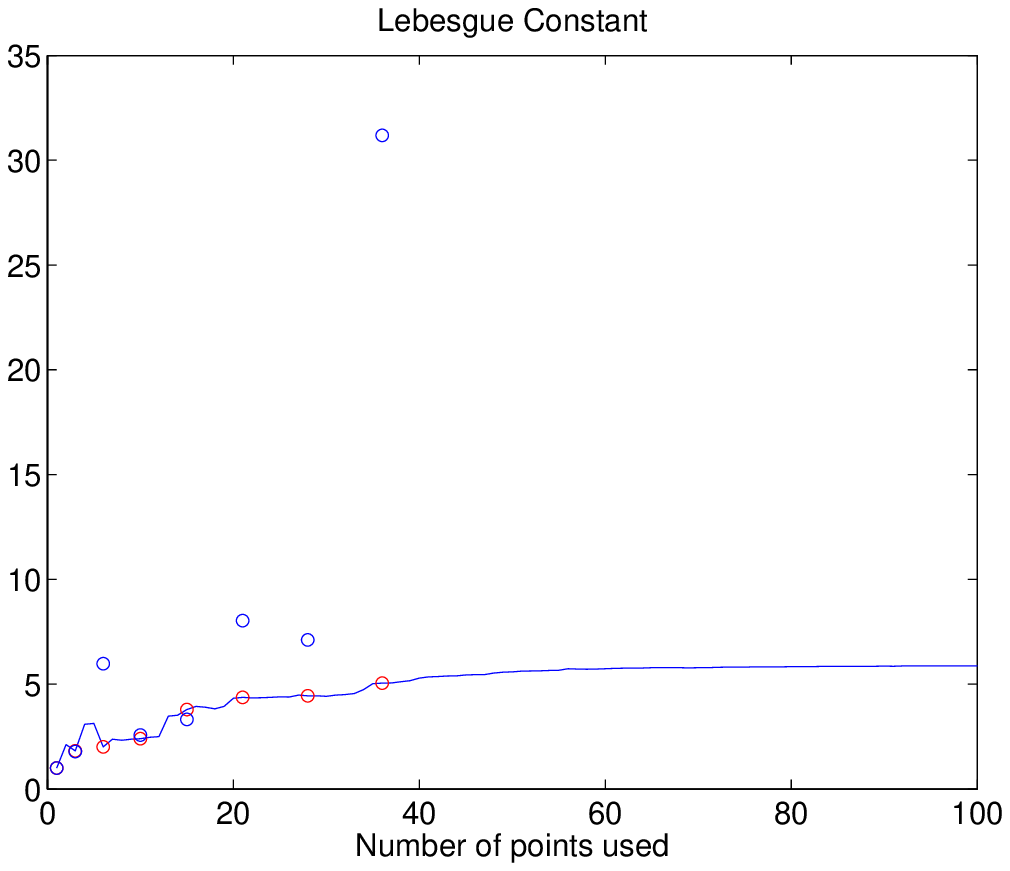}\\
\includegraphics[width=\RSw,height=\RSh]{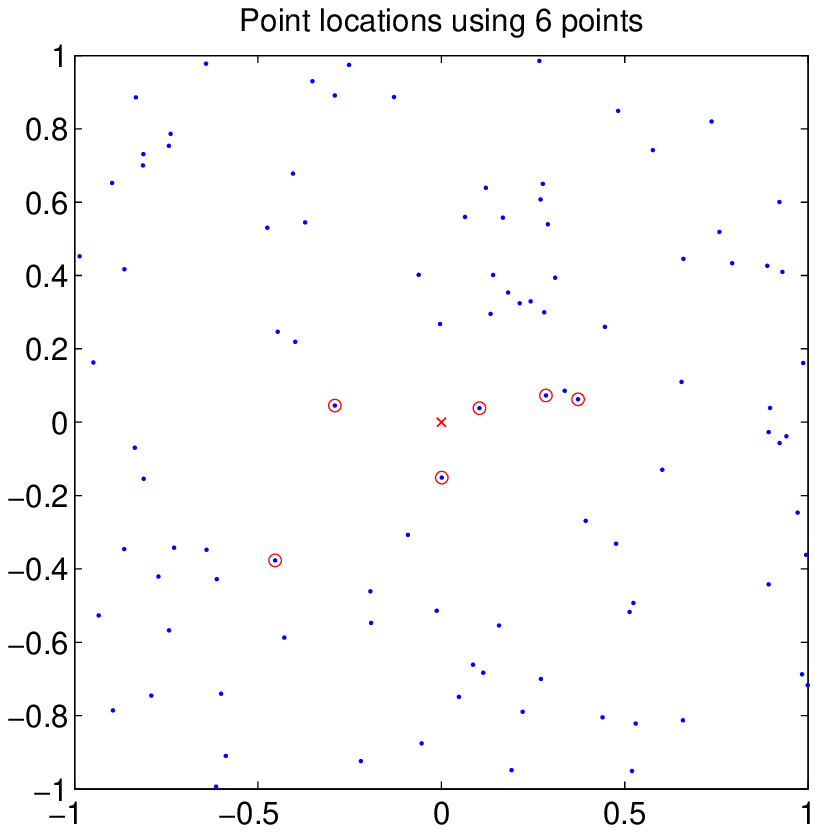}
\includegraphics[width=\RSw,height=\RSh]{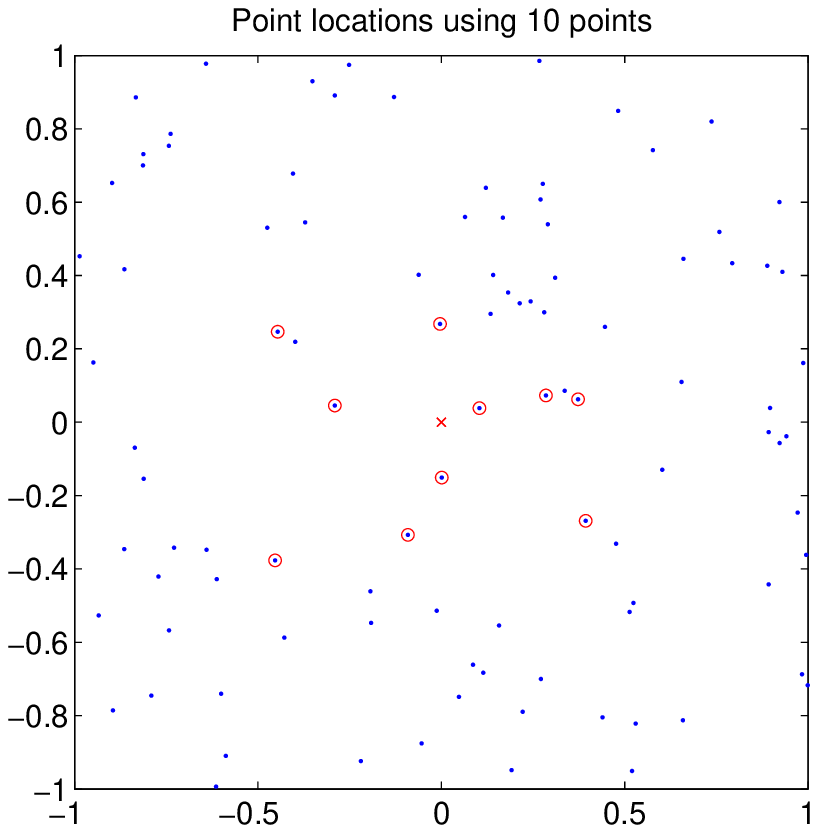}
\includegraphics[width=\RSw,height=\RSh]{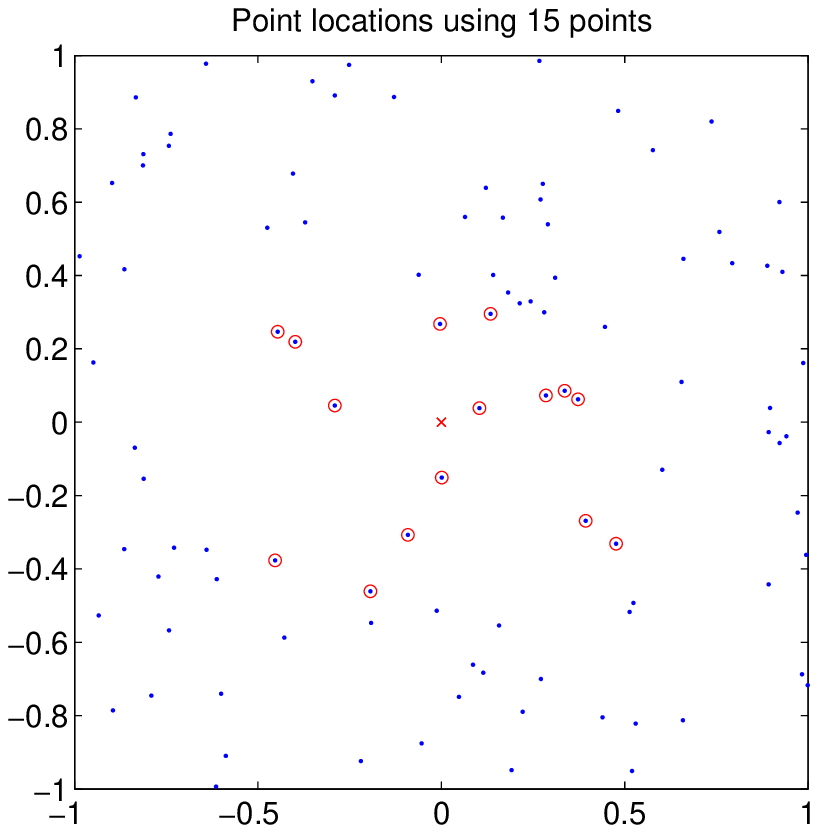}
\end{center} 
    \caption{Local recovery in $W_2^3(\R^2)$ at the origin using the greedy point selection
      strategy.\RSlabel{FigLocalm3}}
\end{figure} 
  \begin{figure}[hbtp]
    \begin{center}
 \includegraphics[width=\RSw,height=\RSh]{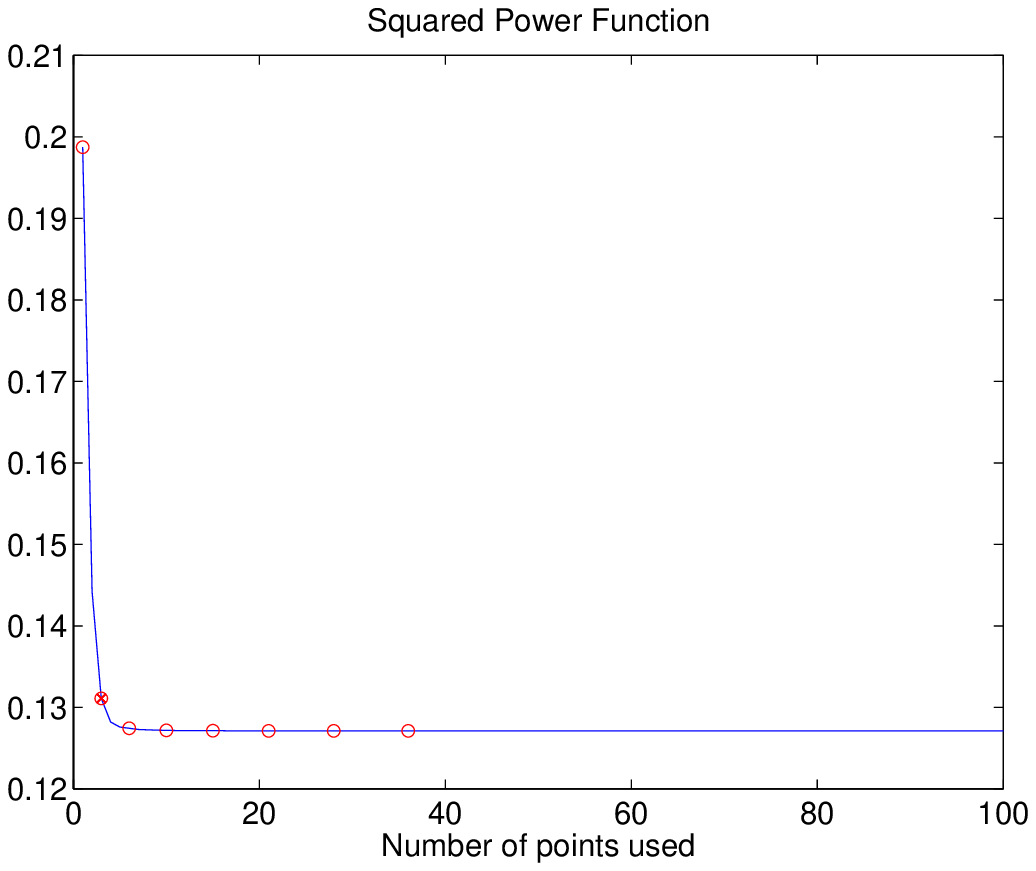}
 \includegraphics[width=\RSw,height=\RSh]{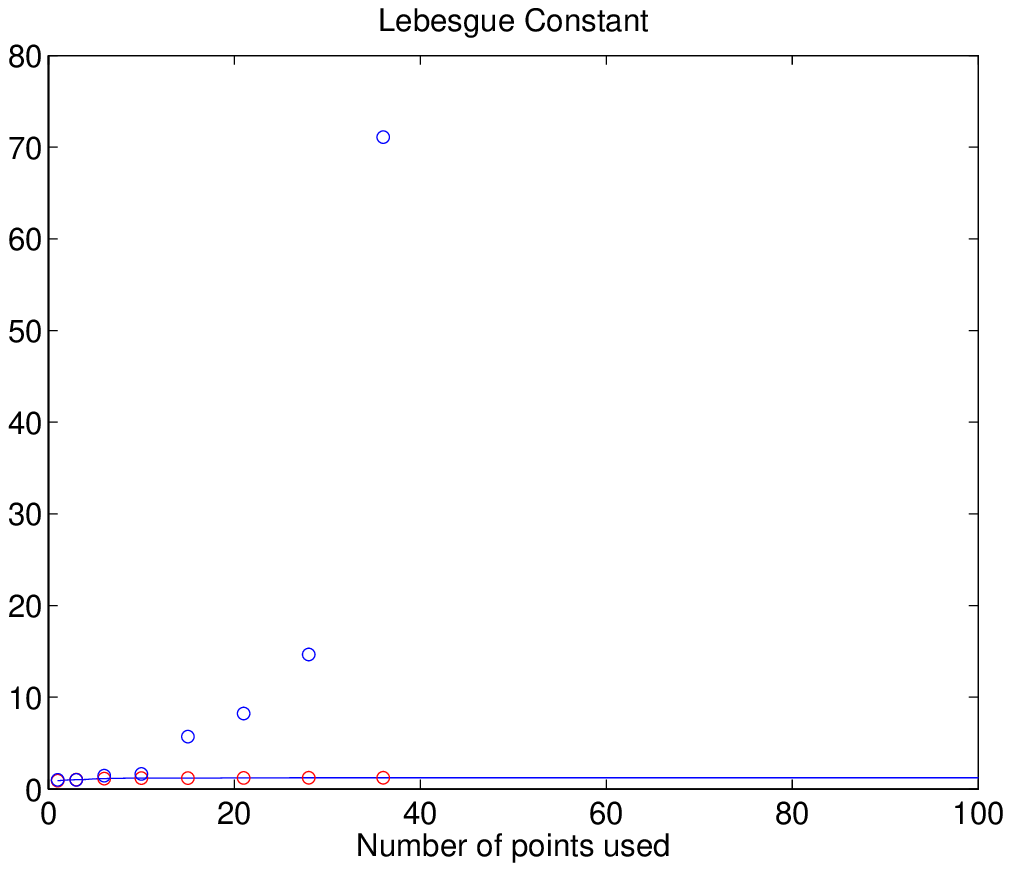}\\
\includegraphics[width=\RSw,height=\RSh]{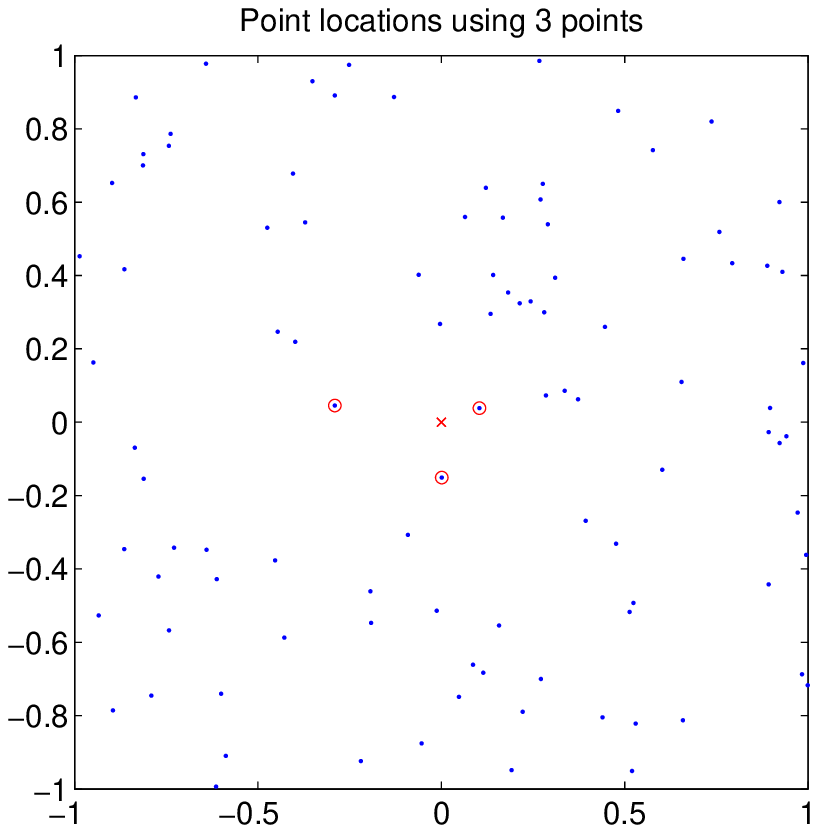}
\includegraphics[width=\RSw,height=\RSh]{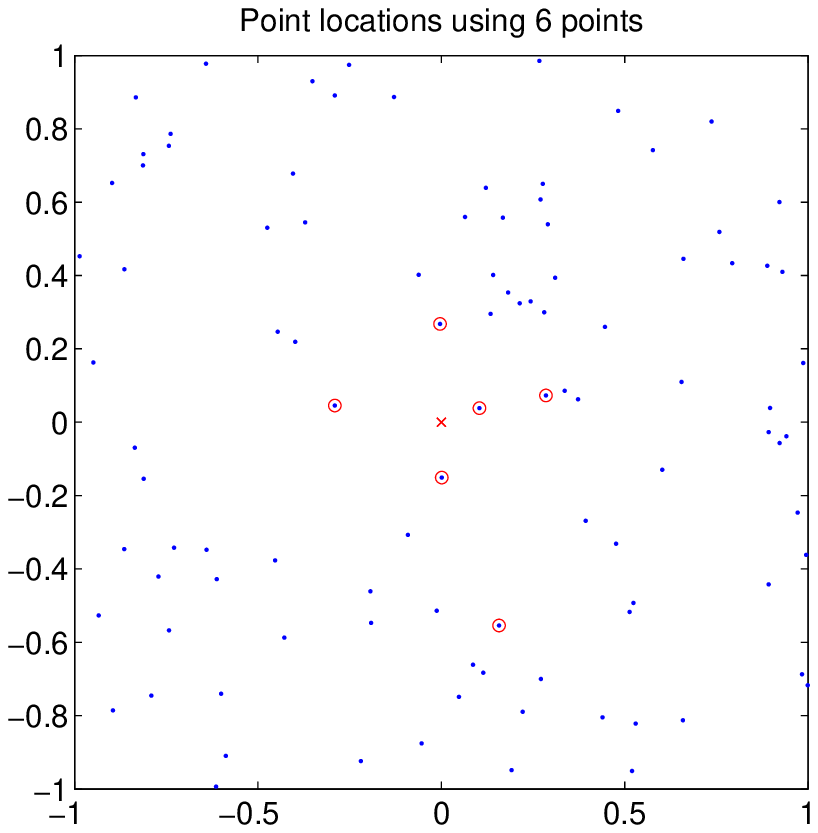}
\includegraphics[width=\RSw,height=\RSh]{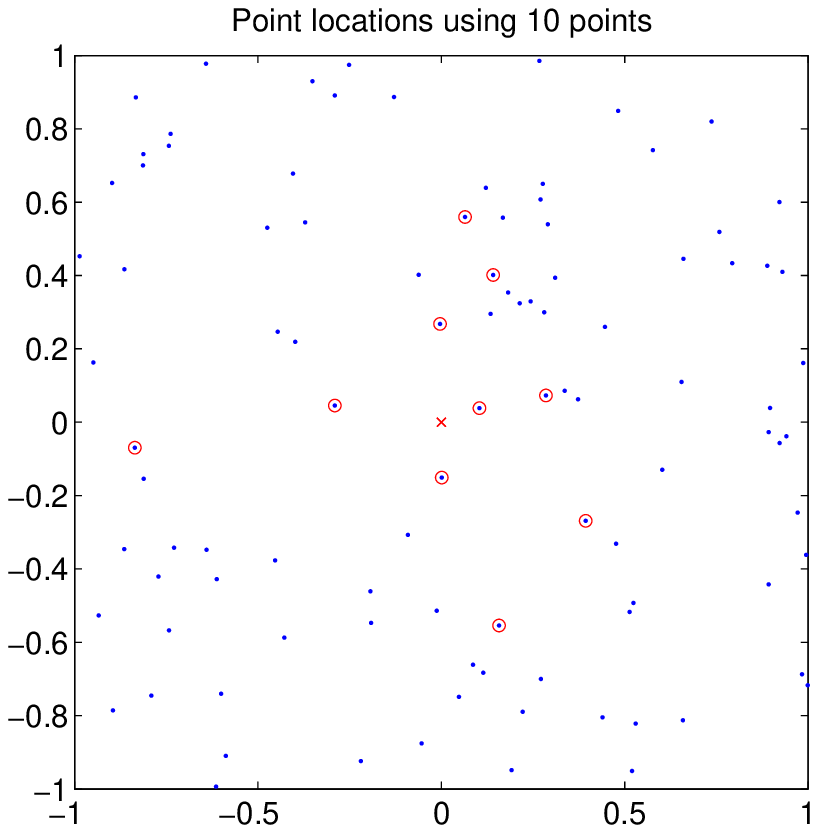}
\end{center} 
    \caption{Local recovery in $W_2^{3/2}(\R^2)$ at the origin using the greedy point selection
      strategy. Three points are enough.\RSlabel{FigLocalm1p5}}
\end{figure} 
  \begin{figure}[hbtp]
    \begin{center}
 \includegraphics[width=\RSw,height=\RSh]{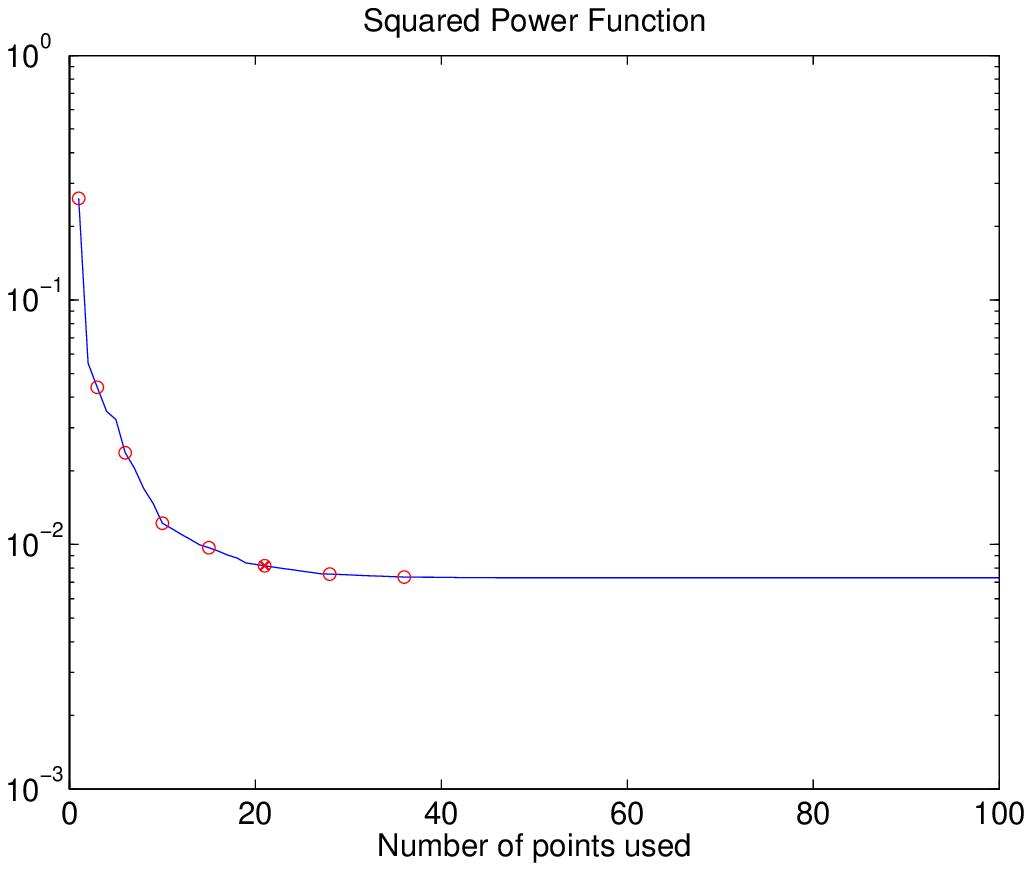}
 \includegraphics[width=\RSw,height=\RSh]{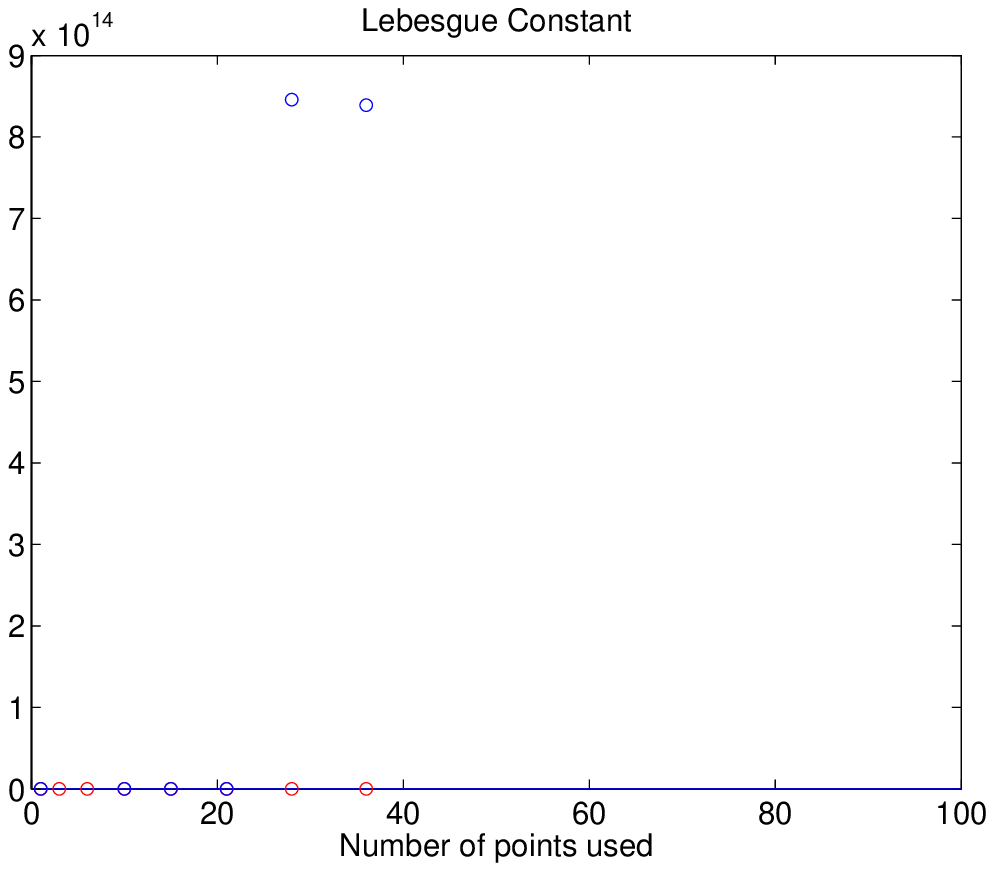}\\
\includegraphics[width=\RSw,height=\RSh]{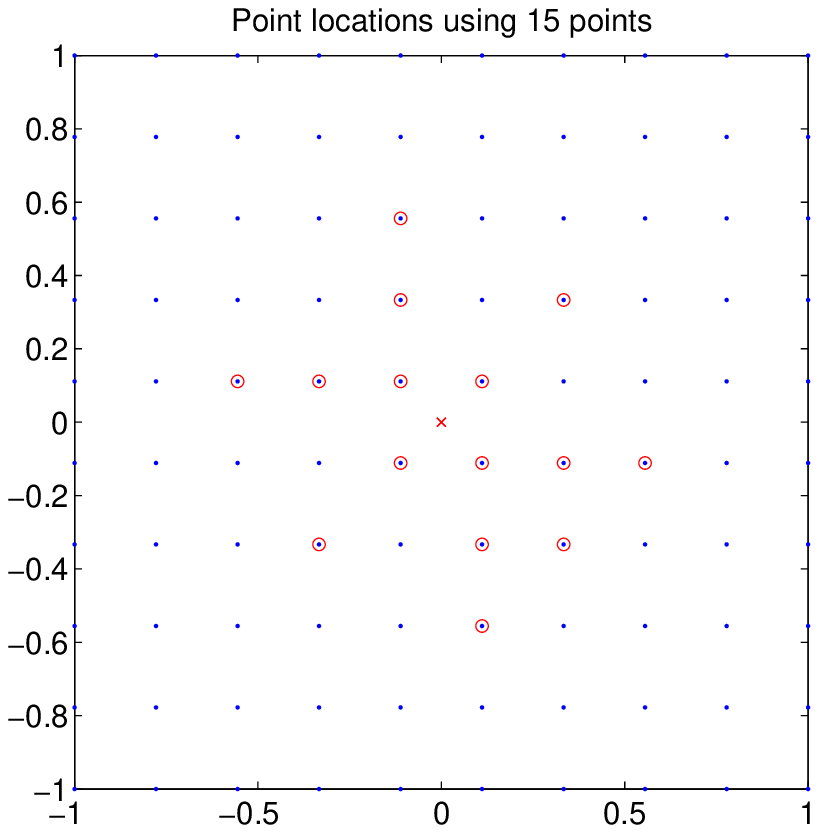}
\includegraphics[width=\RSw,height=\RSh]{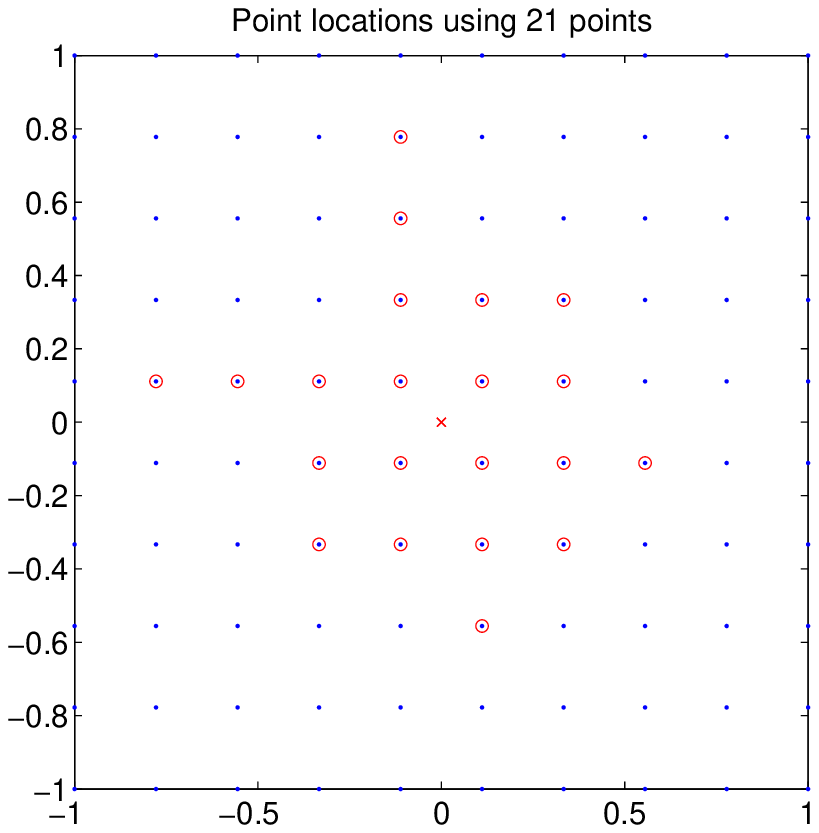}
\includegraphics[width=\RSw,height=\RSh]{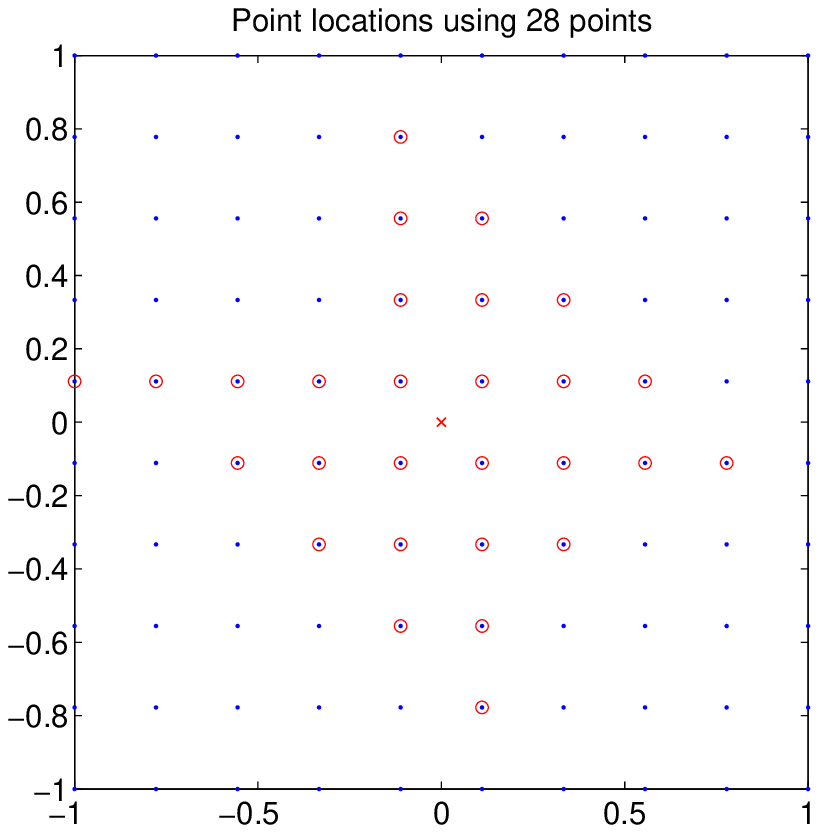}
\end{center} 
    \caption{Local recovery in $W_2^{6}(\R^2)$ using the greedy point selection
      strategy at scale $c=0.1$ on 100 regular points. One should use 21 points.
      \RSlabel{FigLocalm6}}
\end{figure} 
  \begin{figure}[hbtp]
    \begin{center}
 \includegraphics[width=\RSw,height=\RSh]{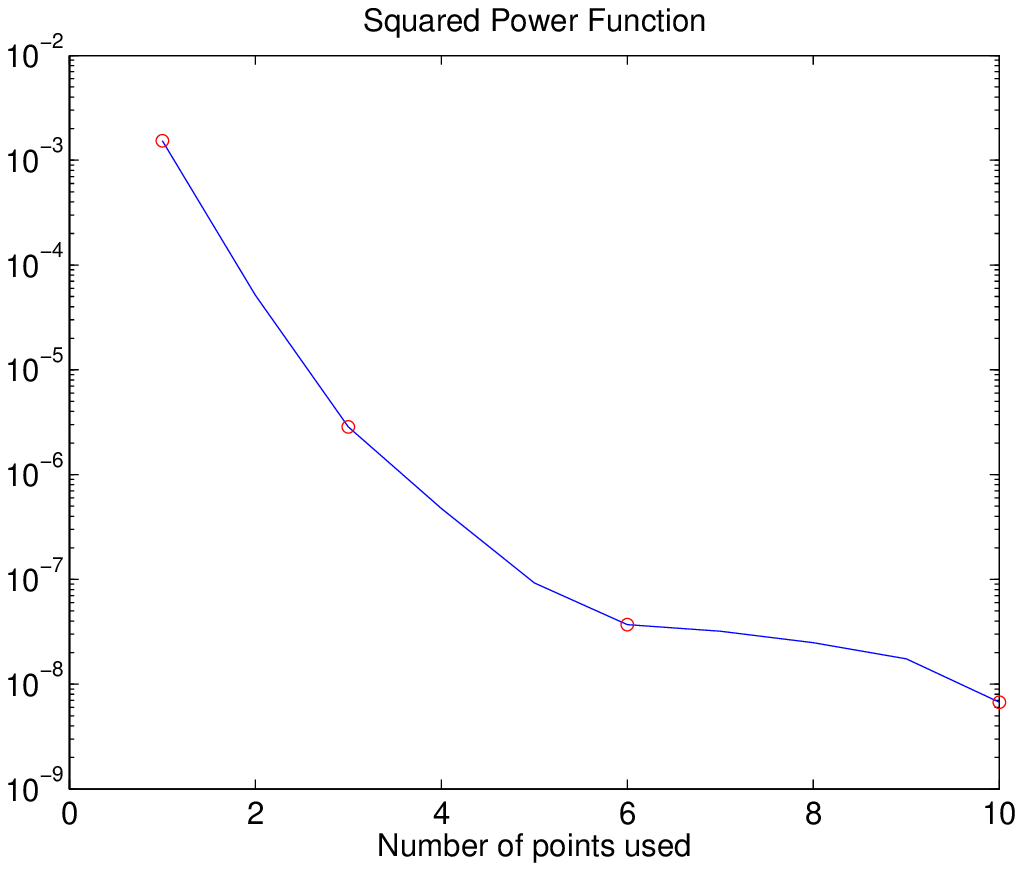}
 \includegraphics[width=\RSw,height=\RSh]{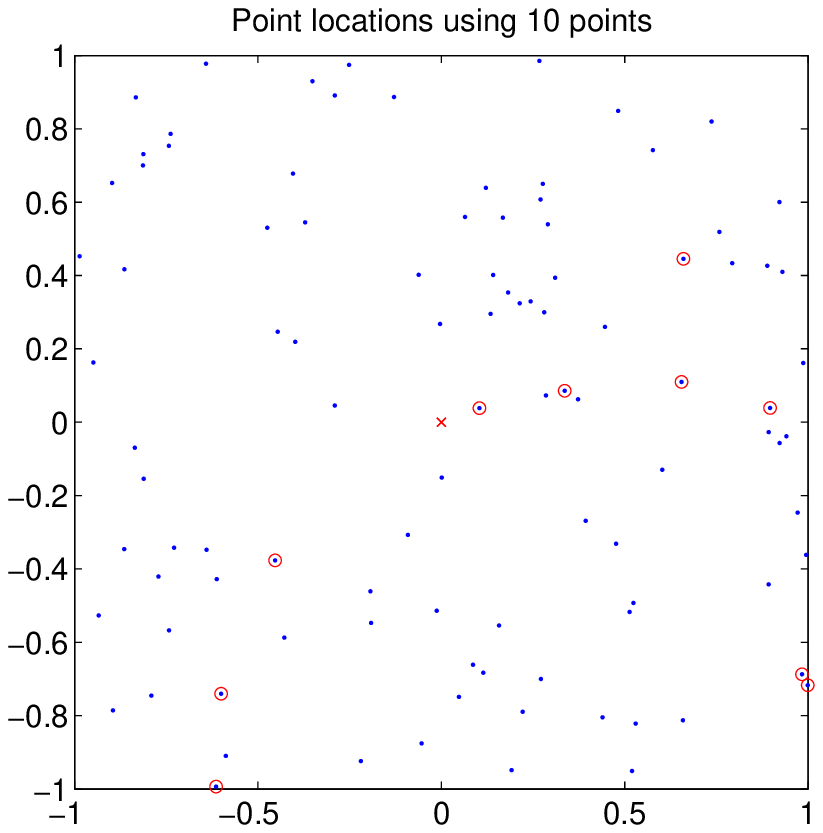}
\end{center} 
    \caption{Local recovery in $W_2^{6}(\R^2)$ at the origin using the greedy point selection
      strategy at scale $c=1$ up to $P^2<1.e-8$ on irregular points.
      One should use 10 points only.
      \RSlabel{FigLocalm6scal1}}
\end{figure} 
  \begin{figure}[hbtp]
    \begin{center}
 
 \includegraphics[width=\RSw,height=\RSh]{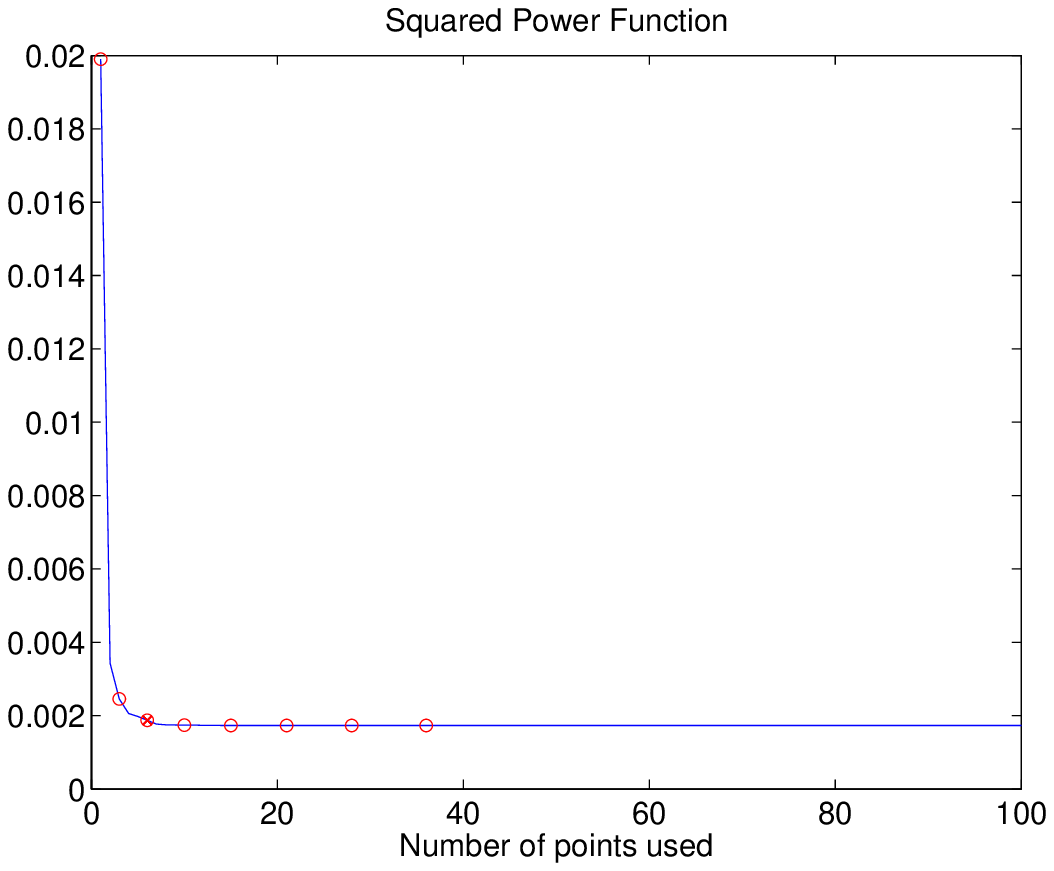}\\

\includegraphics[width=\RSw,height=\RSh]{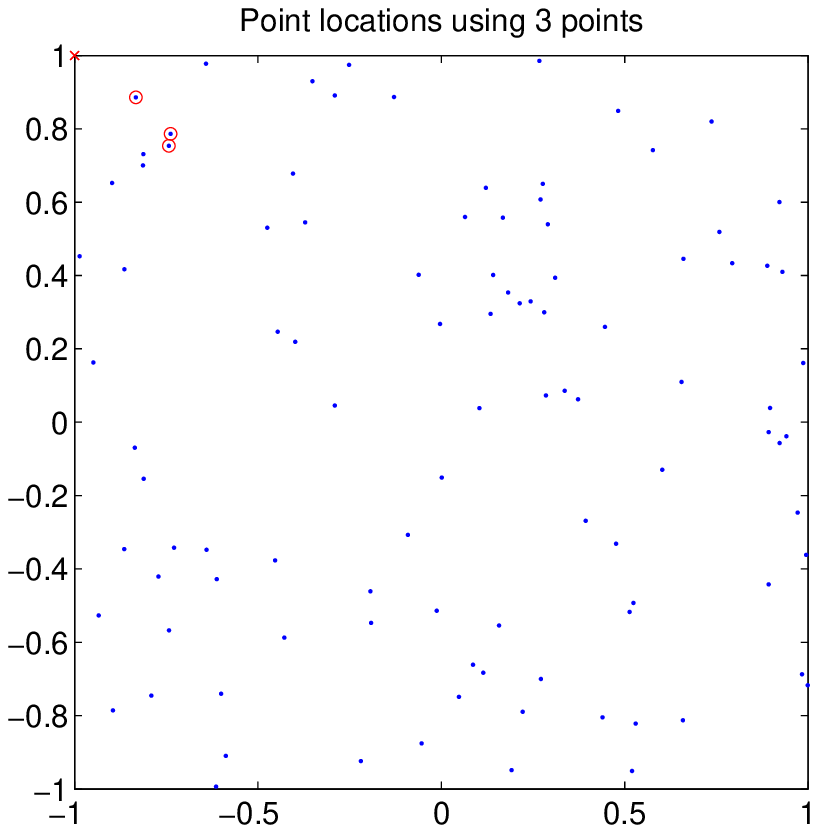}
\includegraphics[width=\RSw,height=\RSh]{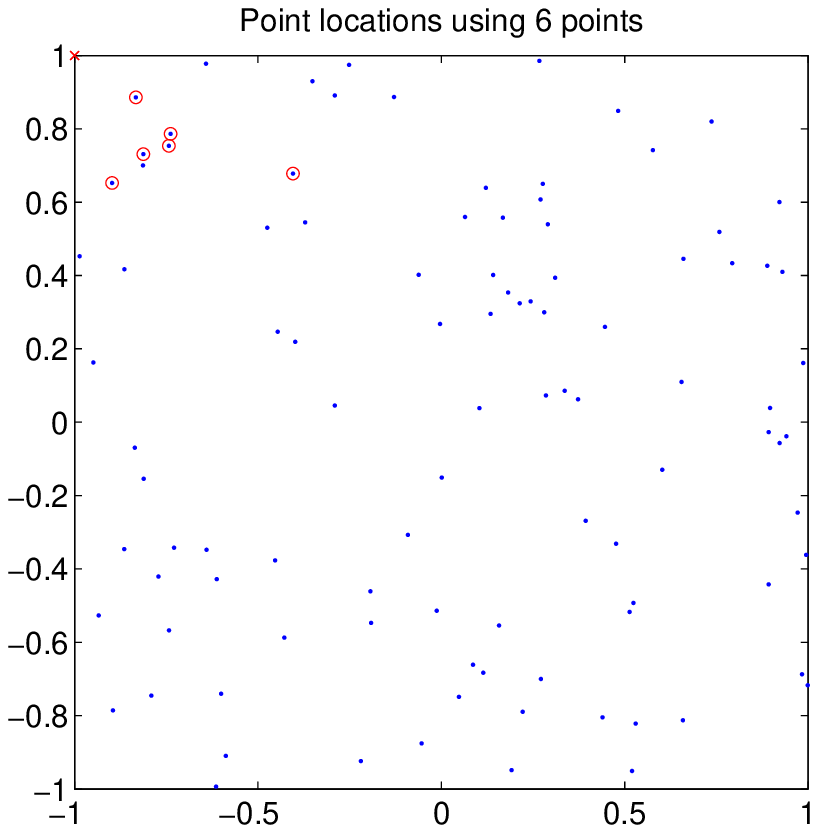}
\includegraphics[width=\RSw,height=\RSh]{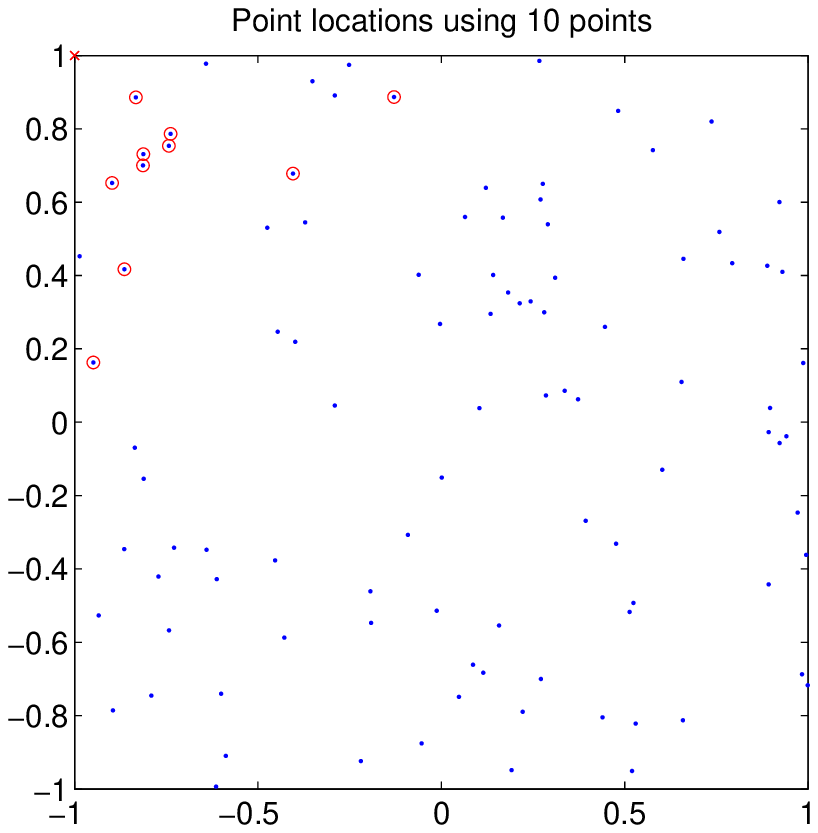}\end{center} 
    \caption{Local recovery in $W_2^3(\R^2)$ at the top left corner using the greedy point selection
      strategy. \RSlabel{FigLocalm3corner}}
\end{figure} 
\biglf
If the Sobolev smoothness order $m$ is fixed, users should aim at the optimal
rate $h^{m-d/2}$ in terms of the fill distance $h$, see the end of Section
\RSref{SecRec}.
In the polynomial situation, one then needs reproduction of order
$q=\lceil m-d/2 \rceil$ and at least $Q={q-1+d \choose d}$ points
\RScite{davydov-schaback:2019-1}.
But in order
to allow greedy point selection it is
recommended to offer more points. Since for
corner points in $\R^2$ three quadrants are missing, one should take at least
$2^d\cdot Q$ points. For $d=2$, we shall use $n=5\cdot Q$ nearest
neighbours in what follows,
and all 2D examples work fine with this choice.
\section{Global Examples}\RSlabel{SecGlobExa}
We use the results of the previous section and the same scattered point
set of $n=100$ points offered for selection.
But now we run the local greedy algorithm on each point of a $51\times 51=2601$
evaluation
grid to simulate an upsampling from rather irregular input data locations to
a regular grid. Readers should see these examples as zoom-ins for much larger
cases. We stick to a $51\times 51$ grid, because we want reasonable plots and a
comparison with the global interpolation on $n$ points.
The sequence of examples matches
the sequence of the previous section. In all cases, we take $5\cdot Q$ points
for $q=\lceil m-d/2\rceil$ and $Q={q-1+d \choose d}$ when we work with the
Mat\'ern kernel generating $W_2^m(\R^d)$. This aims at optimal
convergence rates
$h^{m-d/2}$ in Sobolev spaces.
\biglf
The case $m=3$ requires $q=2$ and $Q=6$ for getting
the expected rate of $h^2$.  We offer the $5\cdot Q=30$ nearest points for each
evaluation point $z$, and let the algorithm select the 6 best ones
to form the set $X(z)$ of Section \RSref{SecRec}, following Figure
\RSref{FigLocalm3}.
The results are in Figure \RSref{FigCholGlob2m3},
plotted over all 2601 evaluation
points.  The top left plot shows the squared Power Function when  
the full interpolation problem is solved on all $n=100$ points. The local
case is top right, working on 6 selected points, and the values must be larger than the
first plot everywhere. The lower left plot shows the difference. It is only
large near the boundary, where the global case has 0.007 and the local case has
0.009. Also, the Lebesgue constants only blow up near the corners. Note
that we allowed 30 points for a
selection of 6. Offering more points does not help. 
\biglf
   \def\RSh{5.5cm}
  \def\RSw{5.5cm}
 \begin{figure}[hbtp]
    \begin{center}
 \includegraphics[width=\RSw,height=\RSh]{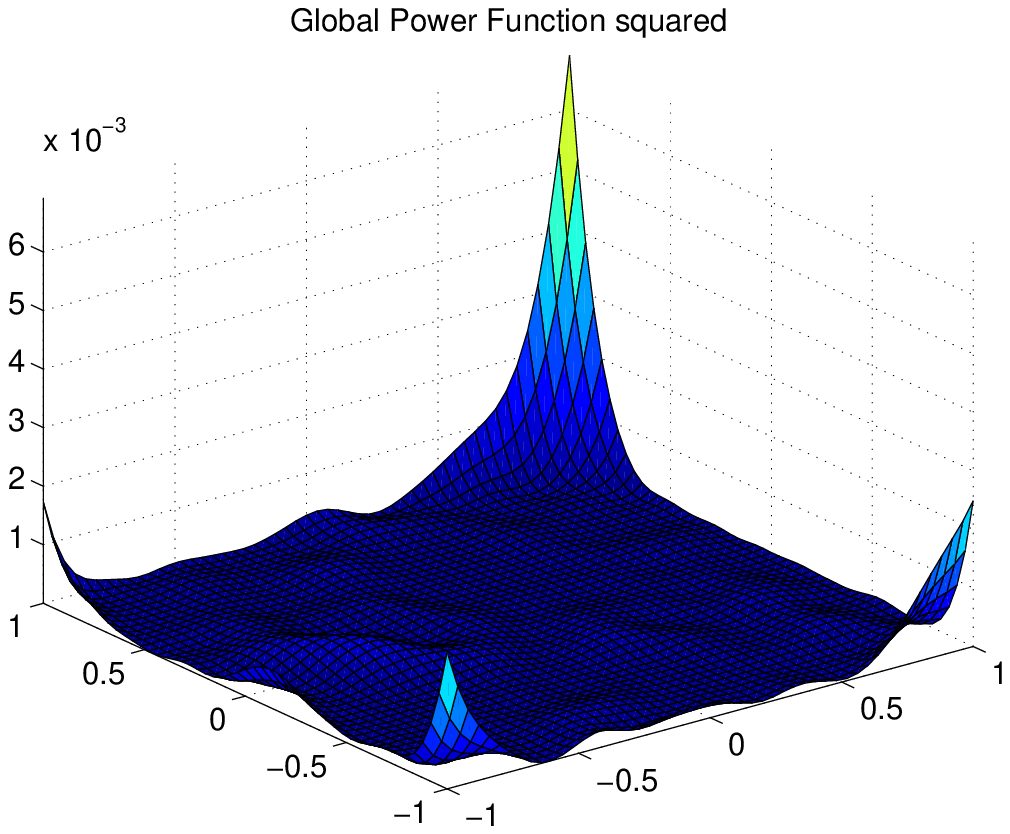}
 \includegraphics[width=\RSw,height=\RSh]{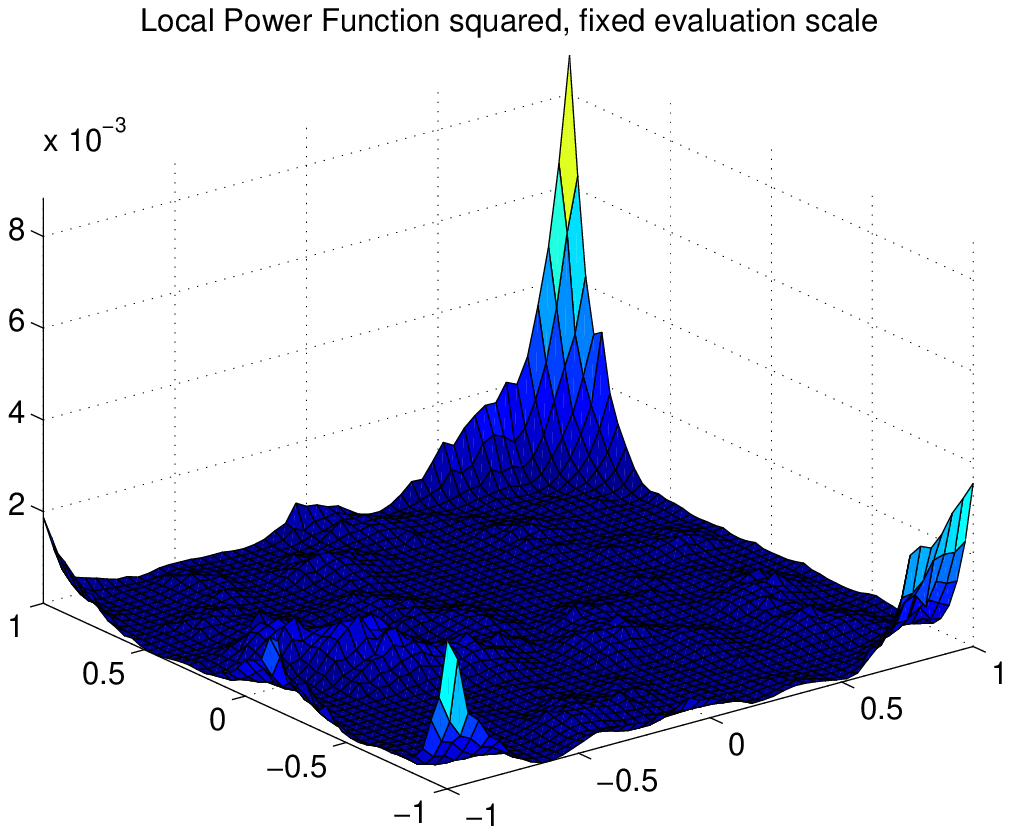}\\
 \includegraphics[width=\RSw,height=\RSh]{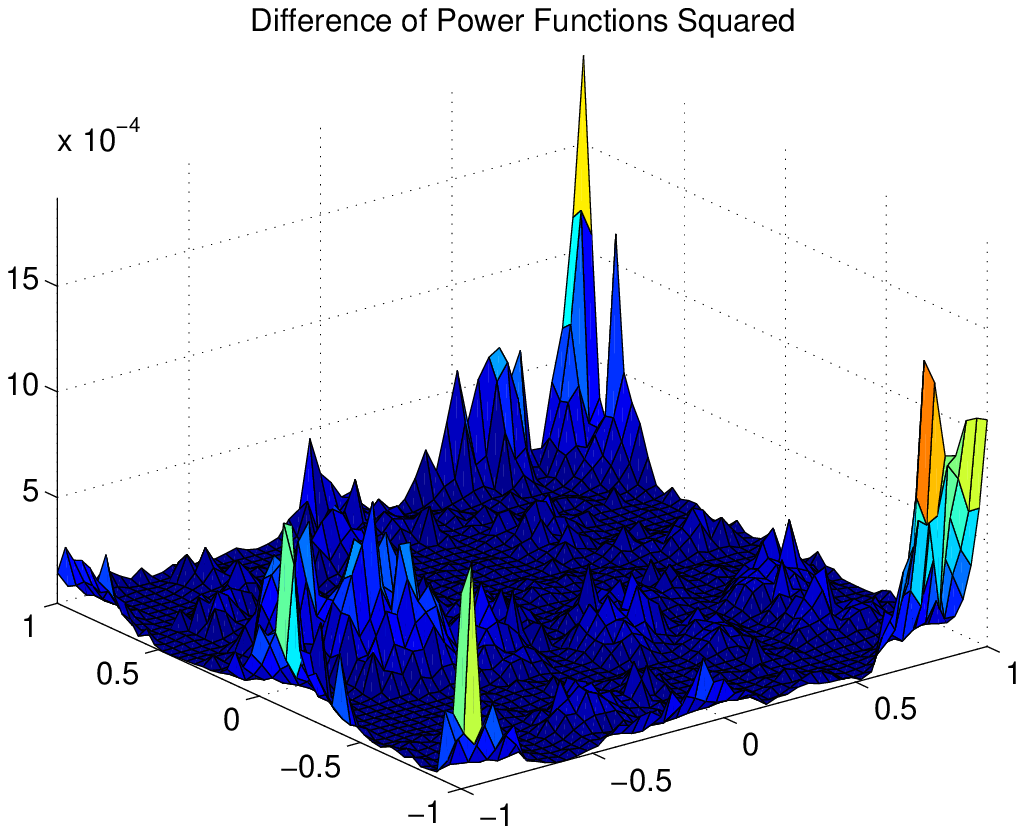}
 \includegraphics[width=\RSw,height=\RSh]{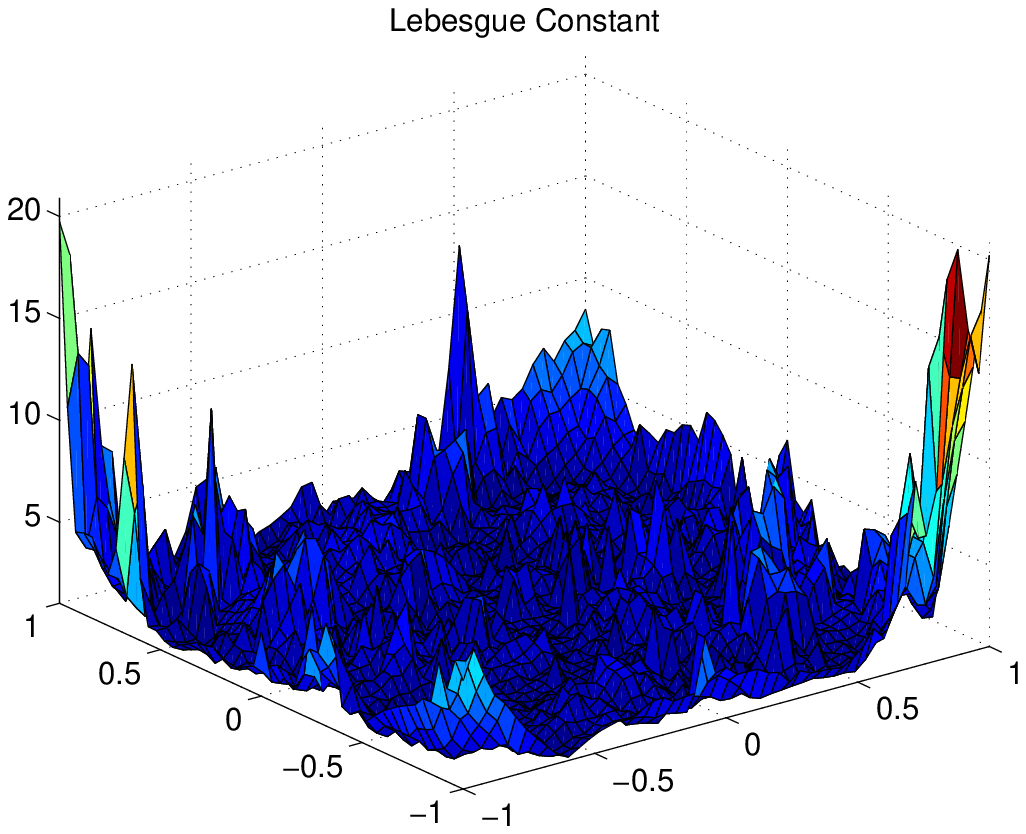}
\end{center} 
    \caption{Local recovery in $W_2^3(\R^2)$ on 2601 points in $[-1,+1]^2$
      using the greedy point selection
      strategy, with selection of 6 points out of 30.\RSlabel{FigCholGlob2m3}}
\end{figure} 
 \biglf
 The case $m=1.5$ is selecting $Q=3$ points out of $5Q=15$
 and produces Figure
 \RSref{FigCholGlob2m1p5}.
 Working locally does not lose more than about 10\%. The case $m=6$ at scale
 $0.1$ on 100 regular data points is in Figure \RSref{FigGlobalm6}. This runs as 
 expected.
 \def\RSh{5.5cm}
  \def\RSw{5.5cm}
 \begin{figure}[hbtp]
    \begin{center}
 \includegraphics[width=\RSw,height=\RSh]{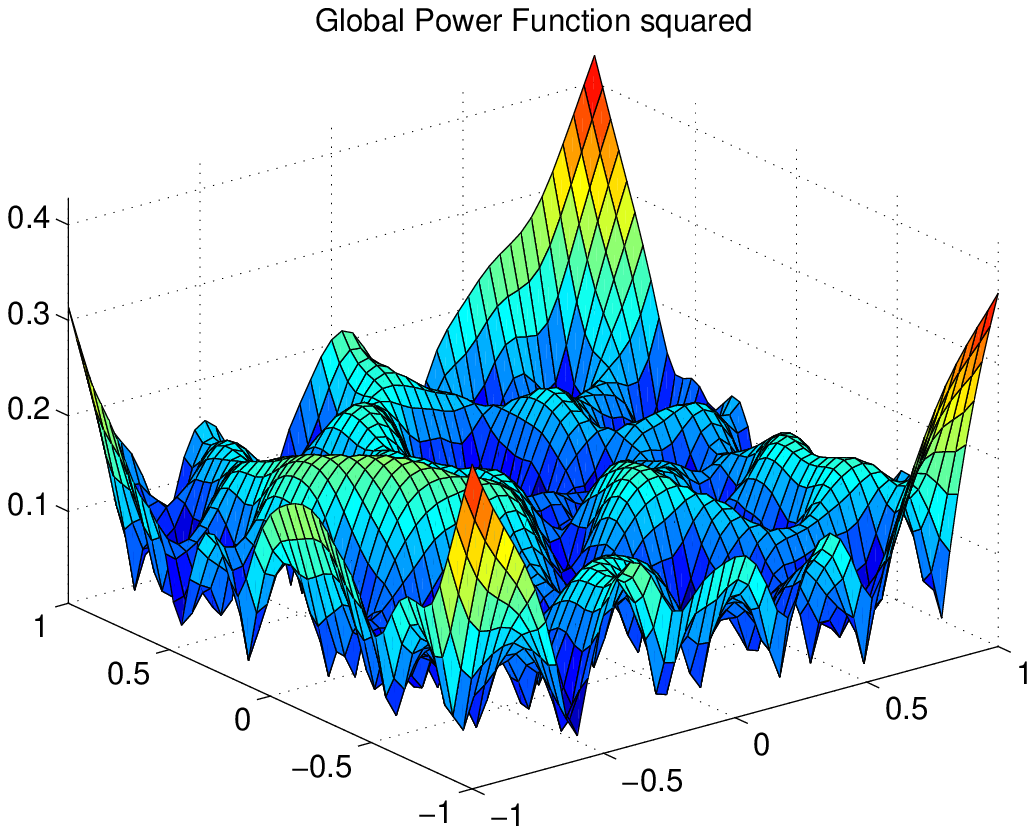}
 \includegraphics[width=\RSw,height=\RSh]{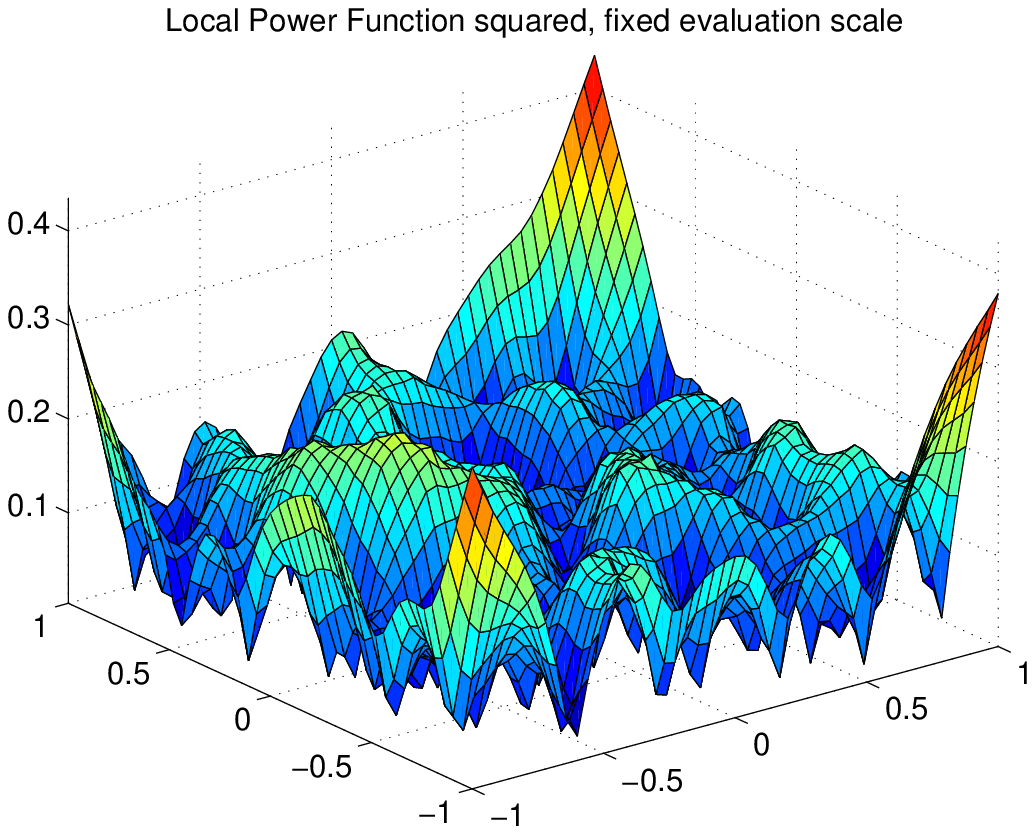}
 \includegraphics[width=\RSw,height=\RSh]{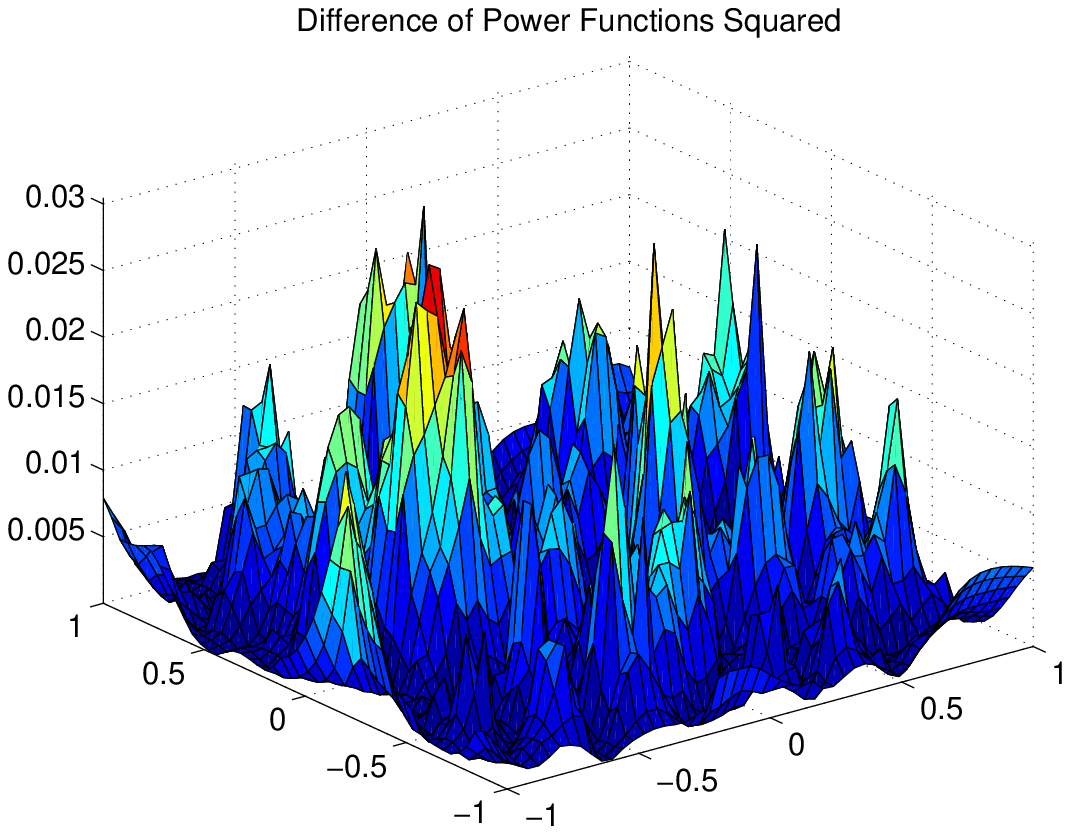}
 \includegraphics[width=\RSw,height=\RSh]{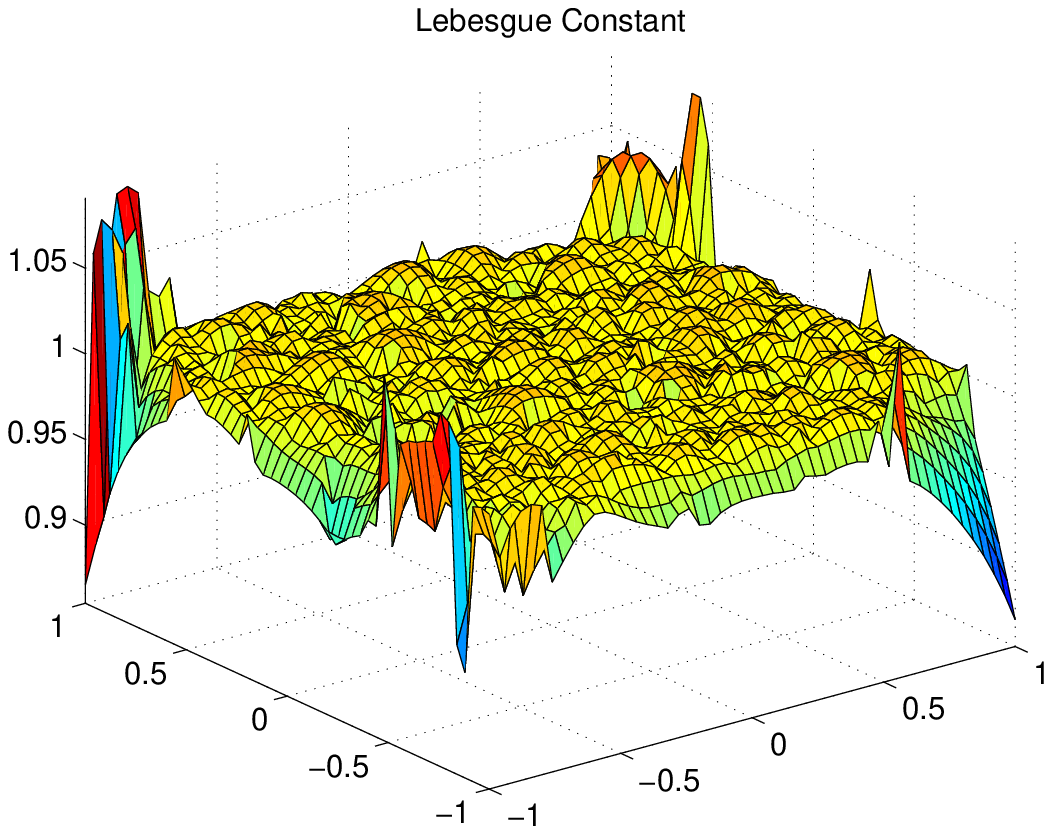}
\end{center} 
    \caption{Local recovery in $W_2^{3/2}(\R^2)$ on 2601 points in $[-1,+1]^2$
      using the greedy point selection
      strategy, with selection of 3 points out of 15.\RSlabel{FigCholGlob2m1p5}}
\end{figure}
 \begin{figure}[hbtp]
    \begin{center}
 \includegraphics[width=\RSw,height=\RSh]{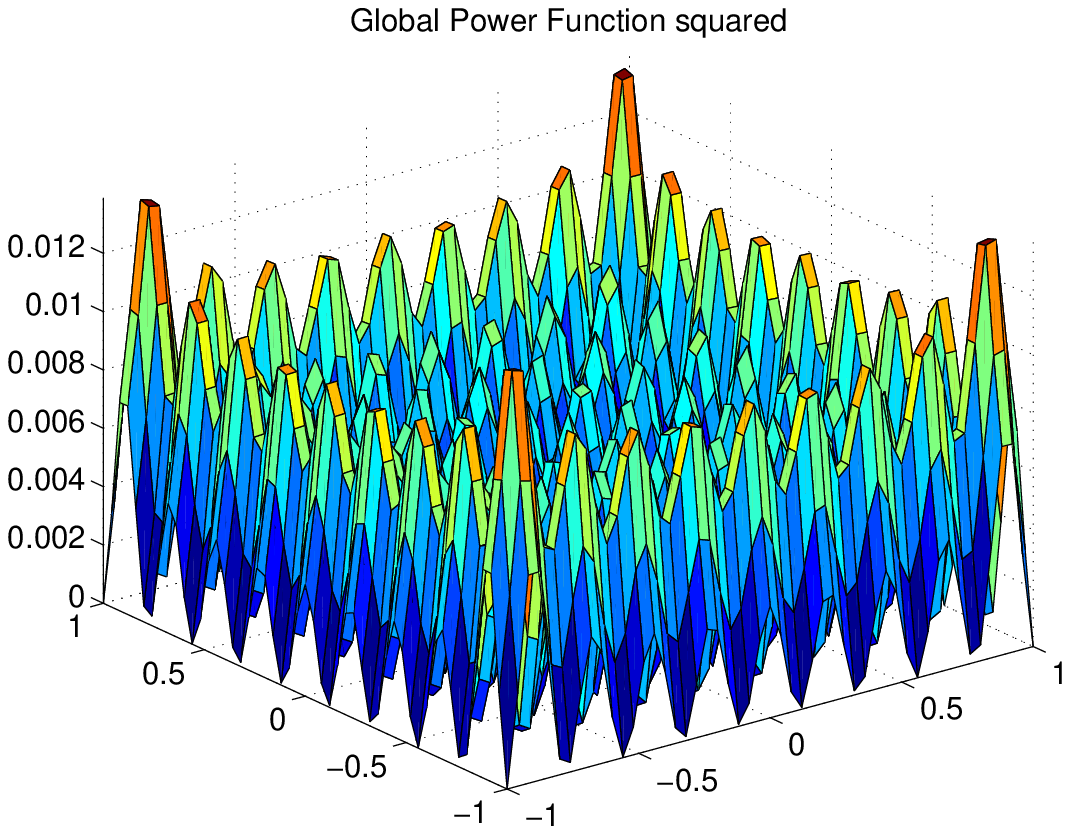}
 \includegraphics[width=\RSw,height=\RSh]{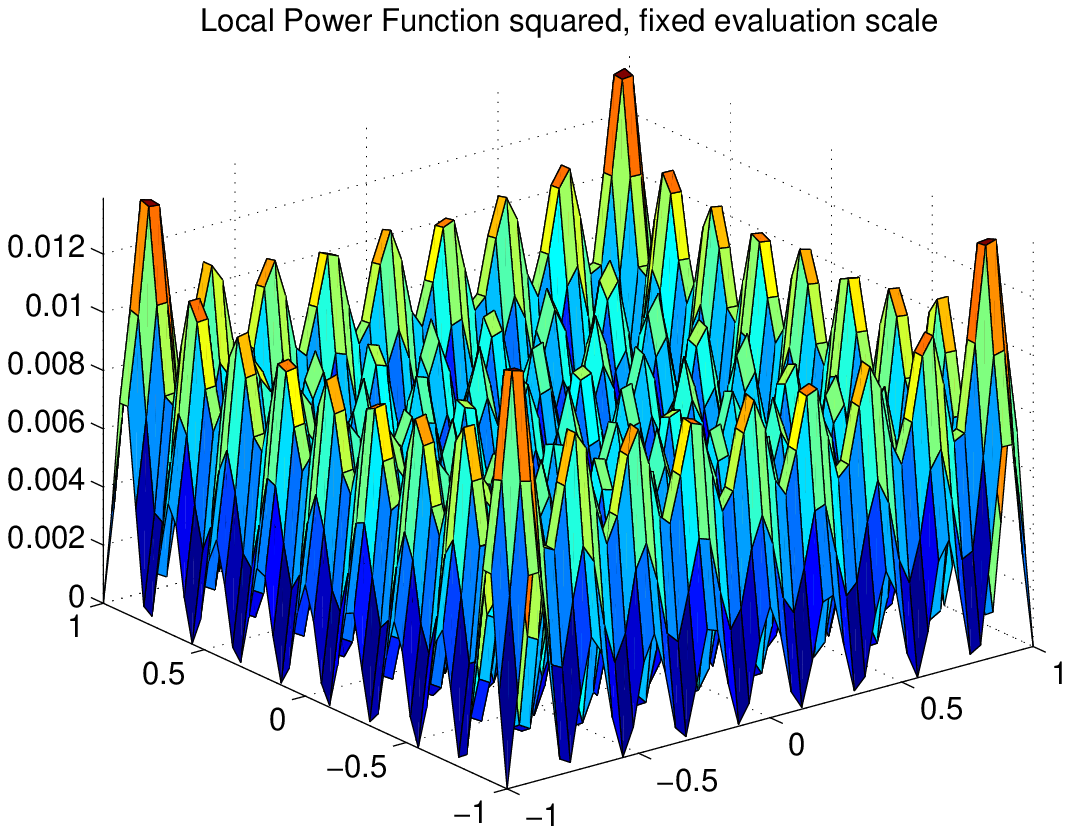}\\
 \includegraphics[width=\RSw,height=\RSh]{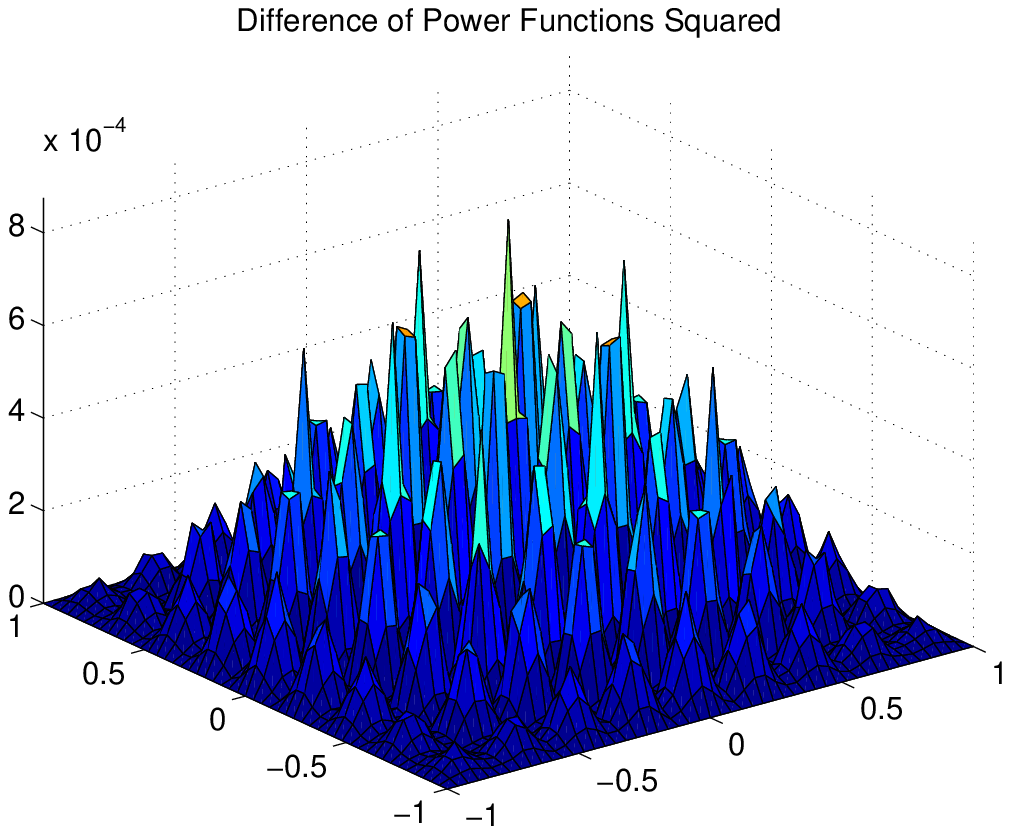}
 \includegraphics[width=\RSw,height=\RSh]{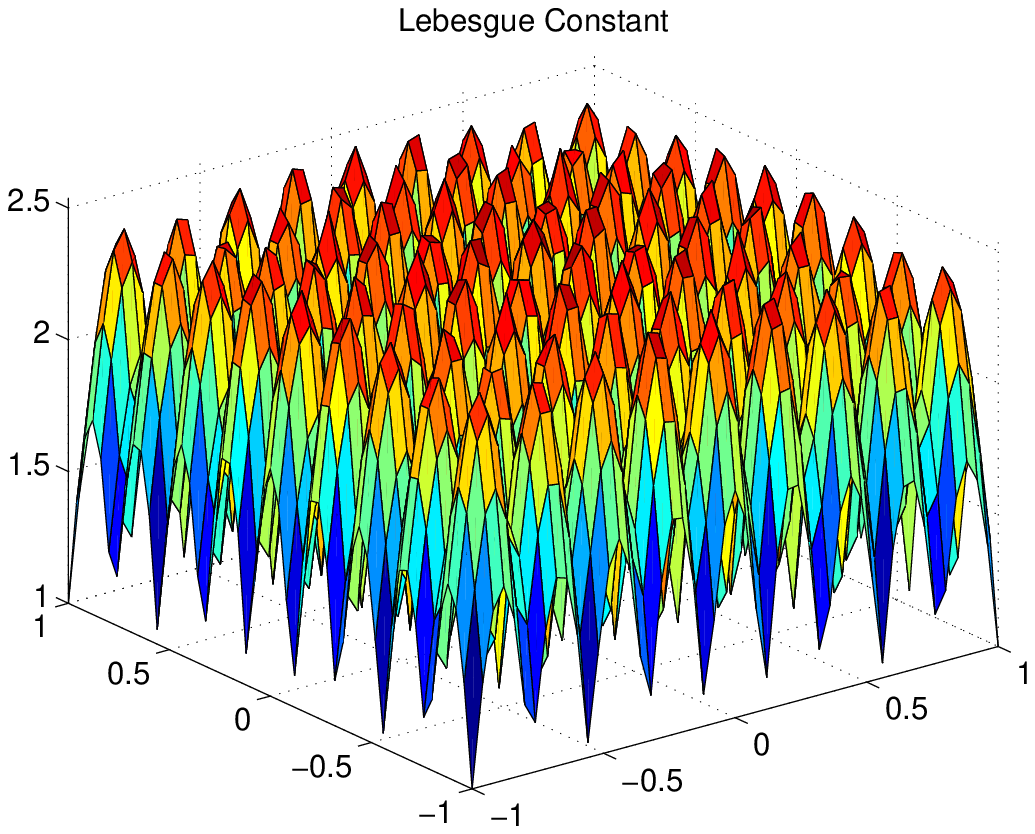}
\end{center} 
    \caption{Local recovery in $W_2^6(\R^2)$ on 2601 regular points in
      $[-1,+1]^2$
      at scale $0.1$ 
      using the greedy point selection
      strategy on 100 regular data points,
      with selection of 21 points out of 100.\RSlabel{FigGlobalm6}}
\end{figure} 
\biglf
 Then we present results for scale
 1.0 that leads to a condition estimate
 of $7.9\cdot 10^{18}$ for the full problem. But
 we stop the greedy algorithm at a tolerance
 of $10^{-8}$ for the squared Power
 Function.
 Here, the global Power Function values may be
 polluted due to the bad condition of the full kernel matrix. The corresponding
 plots are omitted. But the right-hand plot shows that the algorithm often gets
 away with less than 21 points to reach the tolerance. We suggested 10 points
 for Figure \RSref{FigLocalm6scal1},
 but there we worked at a single point in the
 interior.  
 
  \def\RSh{5.5cm}
  \def\RSw{5.5cm}
  \begin{figure}[hbtp]
    \begin{center}
 \includegraphics[width=\RSw,height=\RSh]{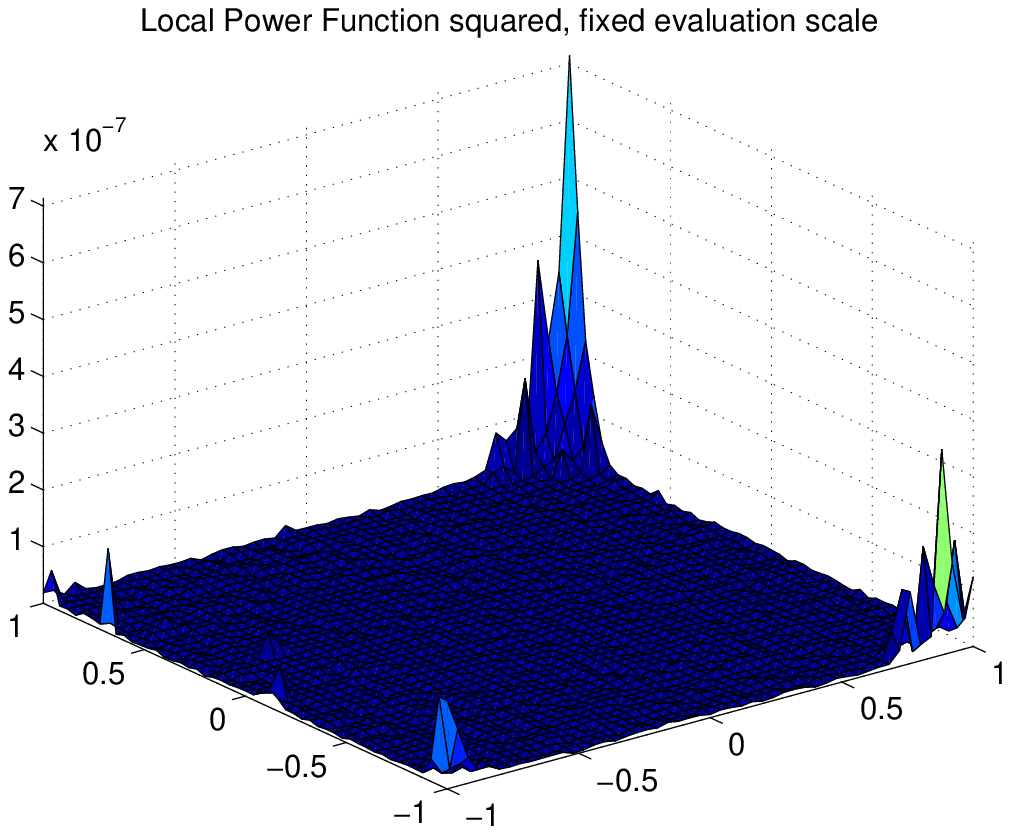}
 \includegraphics[width=\RSw,height=\RSh]{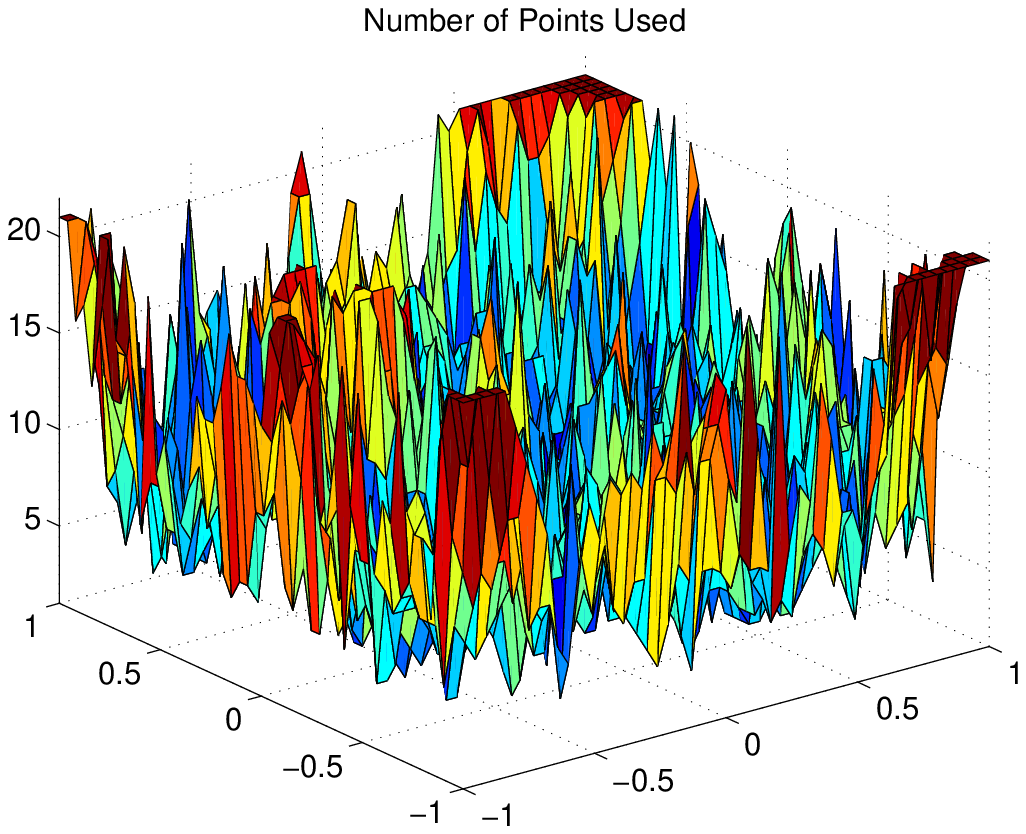}
\end{center} 
    \caption{Local recovery in $W_2^{6}(\R^2)$ on 2601 points in $[-1,+1]^2$
      using the greedy point selection
      strategy, with selection of 21 points out of 100, and with a lower bound
      of $10^{-8}$ for the squared Power Function. .\RSlabel{FigCholGlob2m6thr}}
  \end{figure}
  \section{Function Reproduction}\RSlabel{SecFunRep}
  We now repeat the above examples, but focus on recovery of a smooth function,
  here the {\tt peaks} function of MATLAB. The goal is to show the
  discontinuities induced by working locally.
  \biglf
  Again, we start with $m=3$, and let Figure \RSref{FigCholGlob2m3peaks} show
  results for the situation of Figure \RSref{FigCholGlob2m3}. Note that the
  local version selects only 6 points out of the 30 nearest neighbours.
  Similarly, Figure \RSref{FigCholGlob2m3peaks} presents the case of $m=3/2$,
  selecting 3 out of 15 points.
   \def\RSh{5.5cm}
  \def\RSw{5.5cm}
 \begin{figure}[hbtp]
    \begin{center}
 \includegraphics[width=\RSw,height=\RSh]{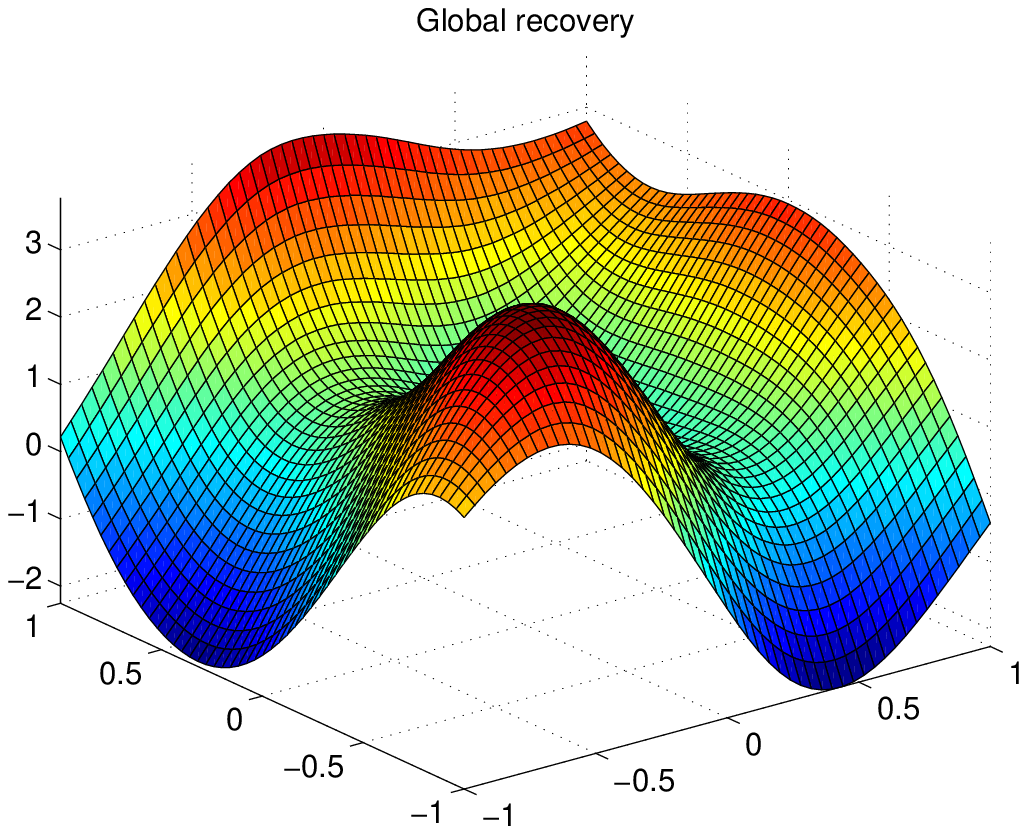}
 \includegraphics[width=\RSw,height=\RSh]{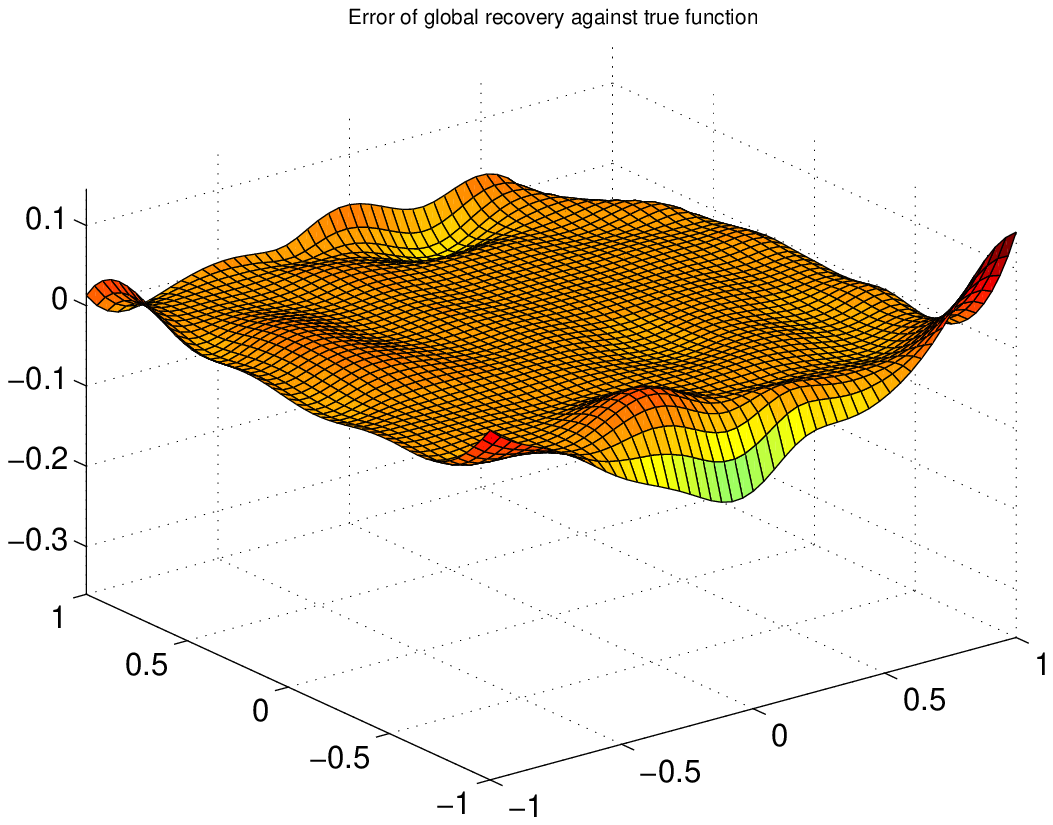}\\
 \includegraphics[width=\RSw,height=\RSh]{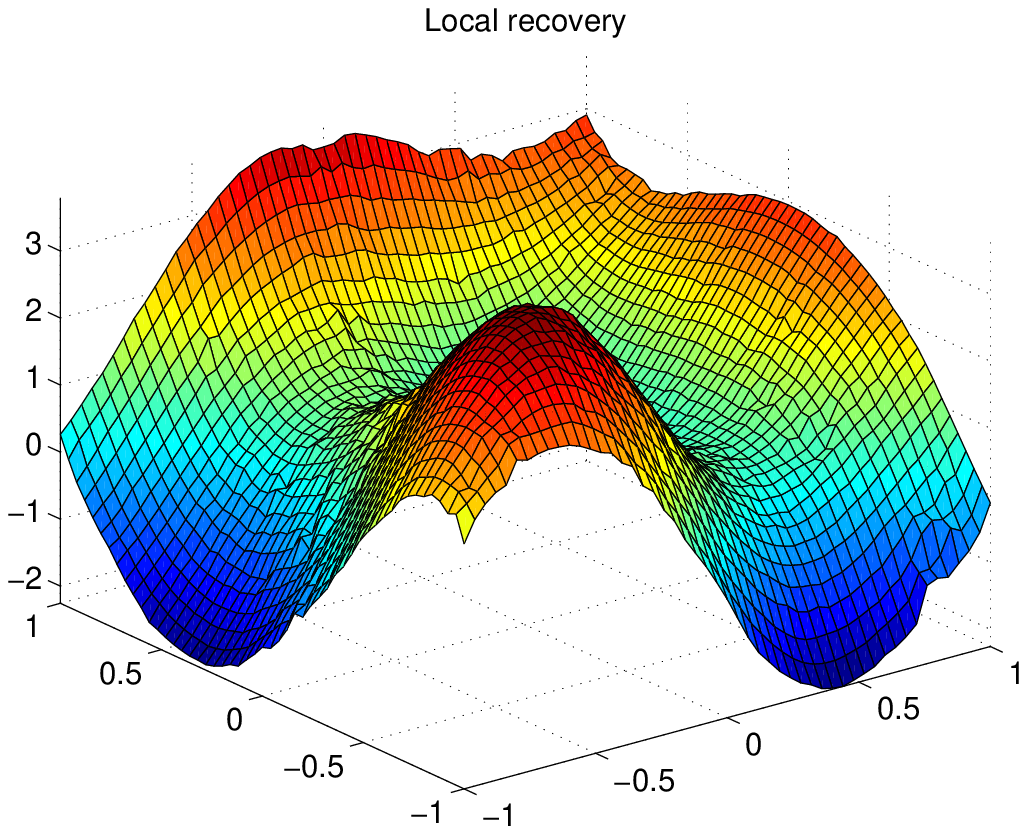}
 \includegraphics[width=\RSw,height=\RSh]{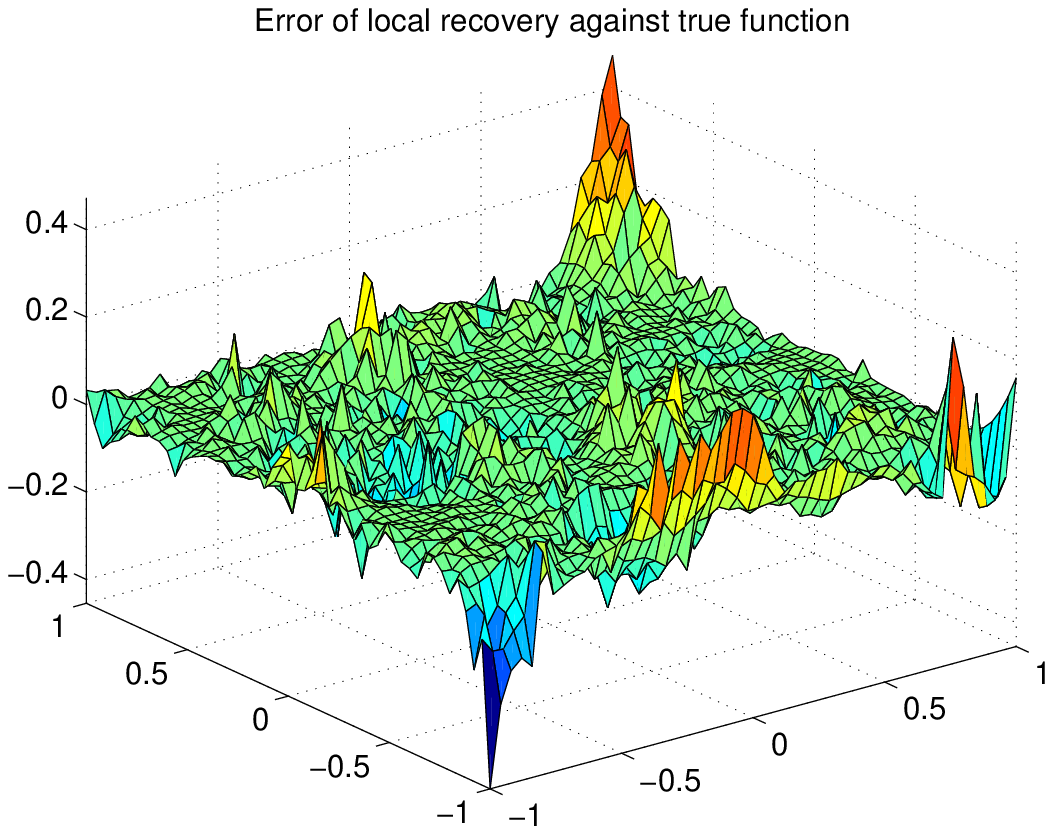}
\end{center} 
    \caption{Local recovery of the {\tt peaks} function  in $W_2^3(\R^2)$, on 2601 points in $[-1,+1]^2$
      using the greedy point selection
      strategy, with selection of 6 points out of 30.\RSlabel{FigCholGlob2m3peaks}}
\end{figure} 
\begin{figure}[hbtp]
    \begin{center}
 \includegraphics[width=\RSw,height=\RSh]{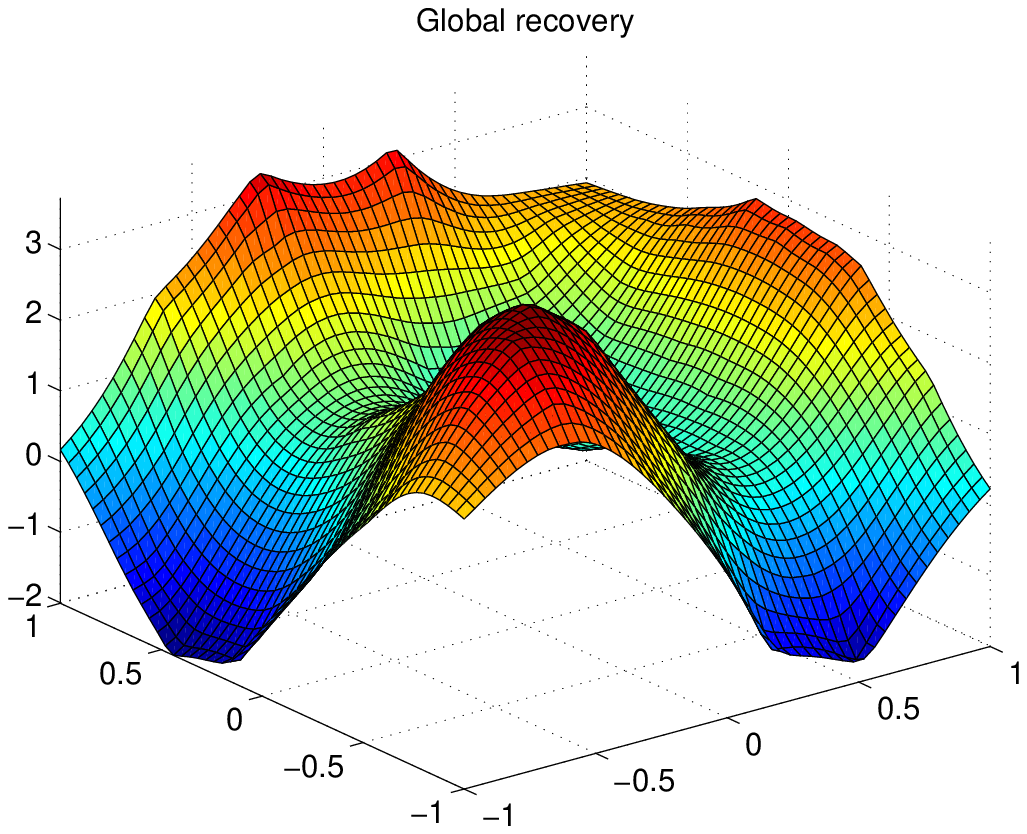}
 \includegraphics[width=\RSw,height=\RSh]{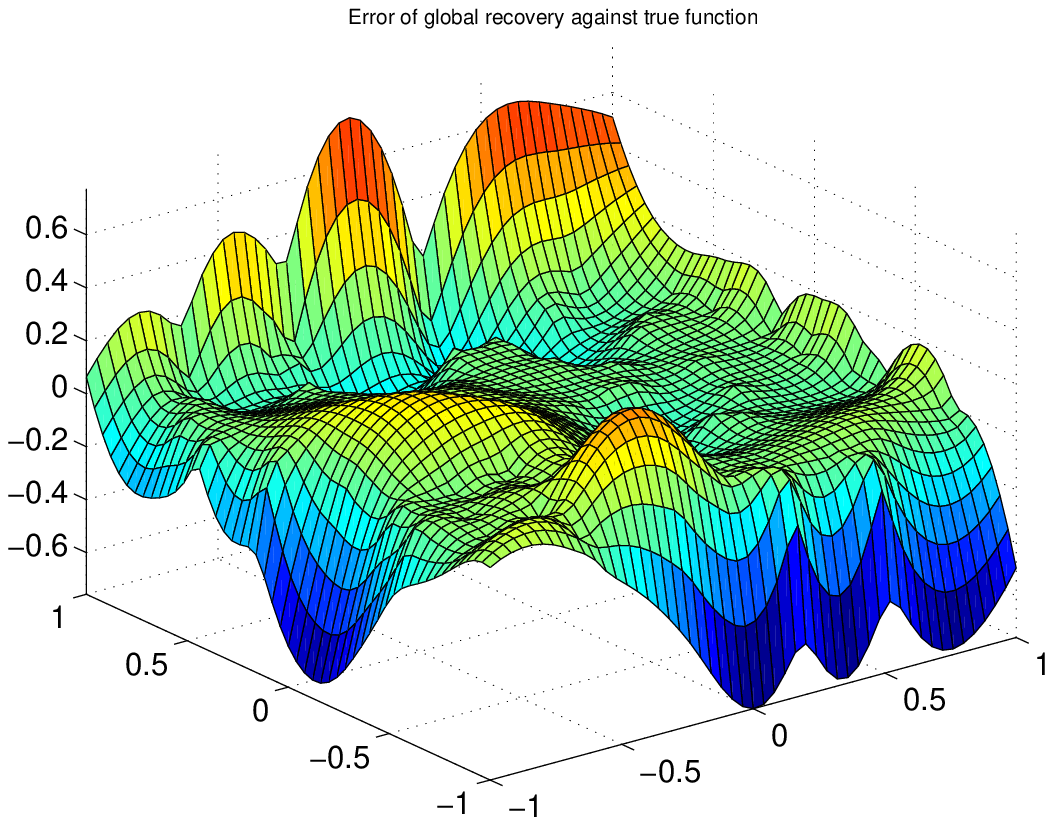}\\
 \includegraphics[width=\RSw,height=\RSh]{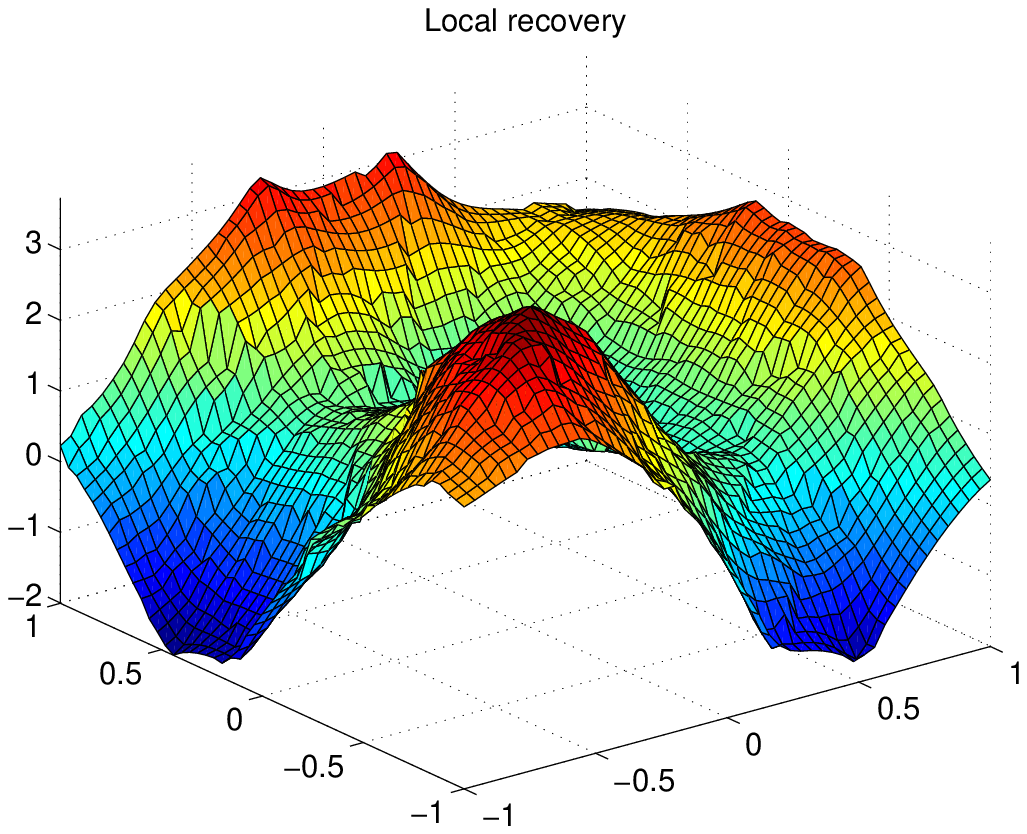}
 \includegraphics[width=\RSw,height=\RSh]{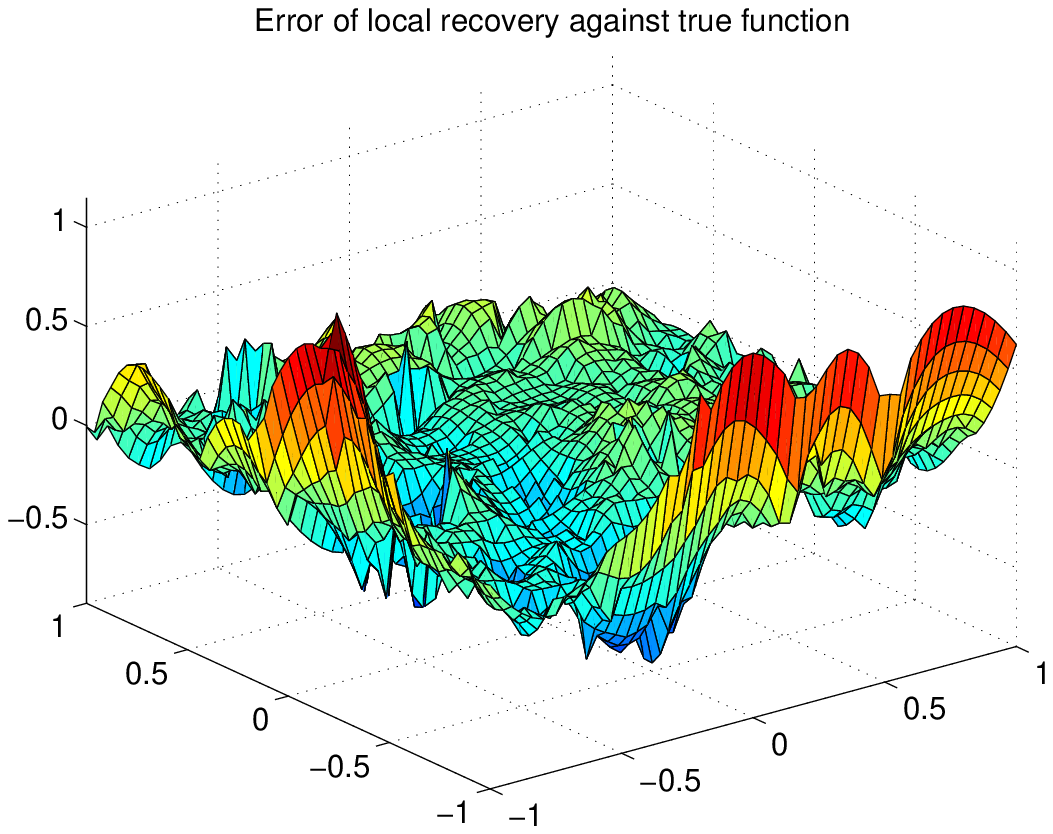}
\end{center} 
    \caption{Local recovery of the {\tt peaks} function  in $W_2^{3/2}(\R^2)$, on 2601 points in $[-1,+1]^2$
      using the greedy point selection
      strategy, with selection of 3 points out of 15.\RSlabel{FigCholGlob2m1p5peaks}}
\end{figure} 
\begin{figure}[hbtp]
    \begin{center}
 \includegraphics[width=\RSw,height=\RSh]{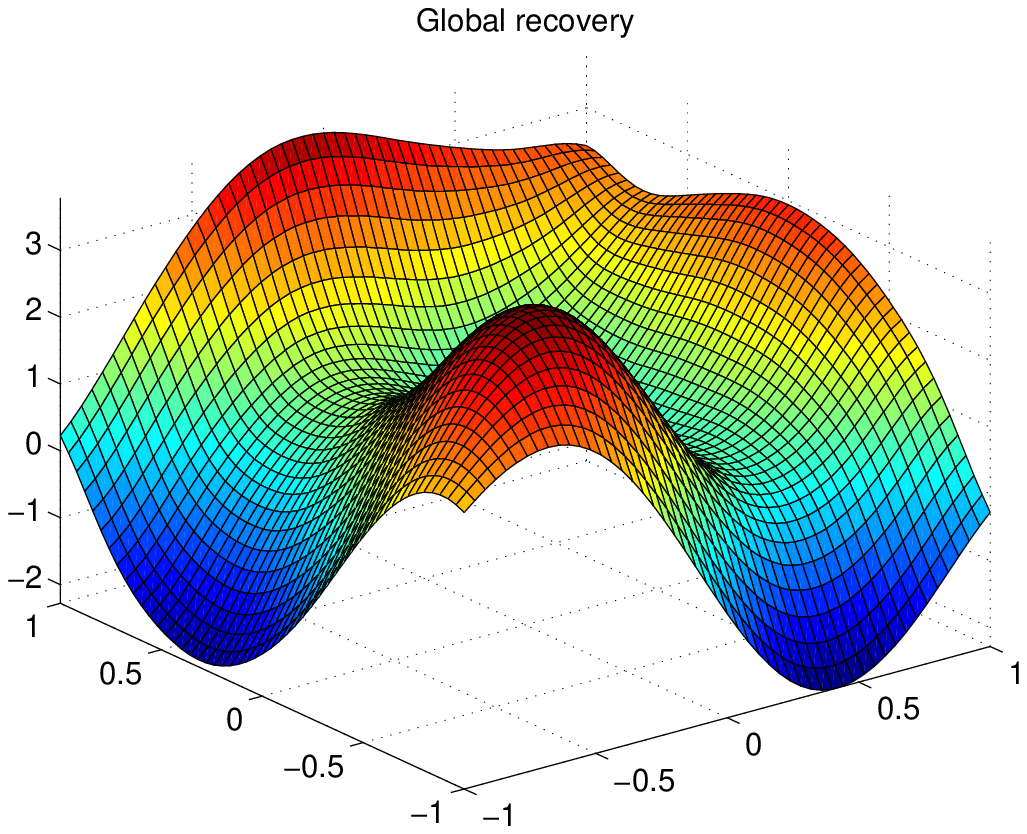}
 \includegraphics[width=\RSw,height=\RSh]{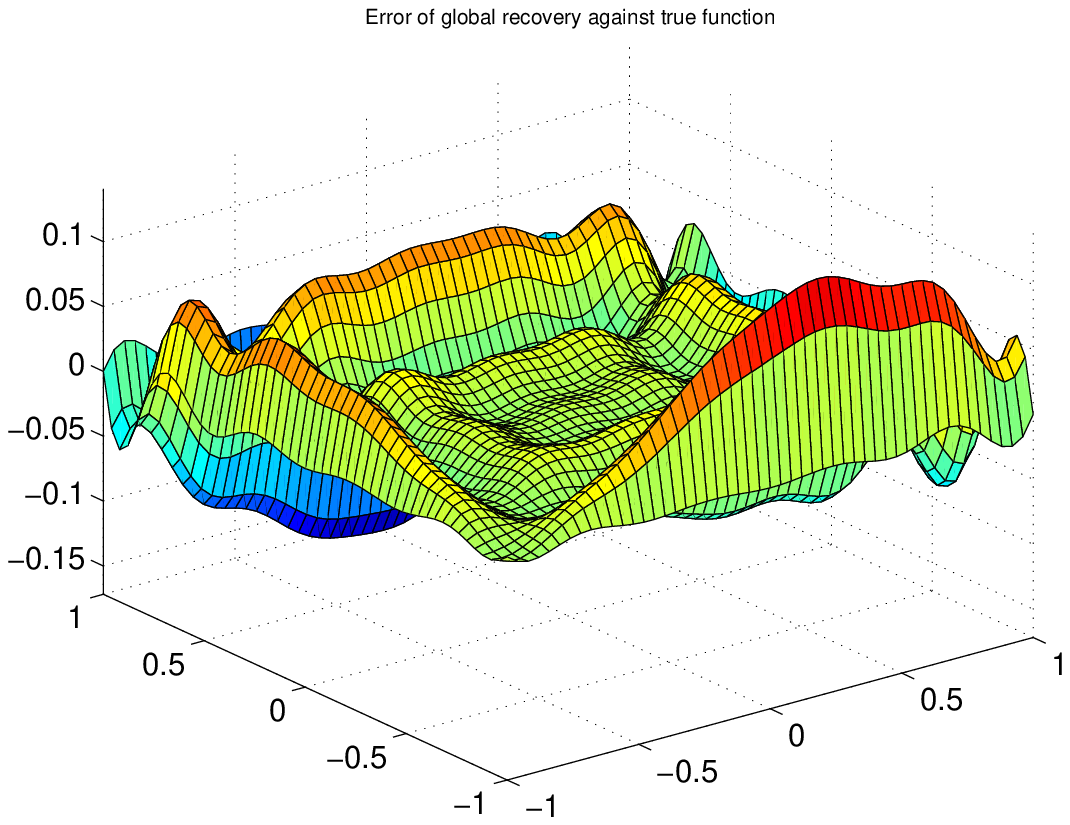}\\
 \includegraphics[width=\RSw,height=\RSh]{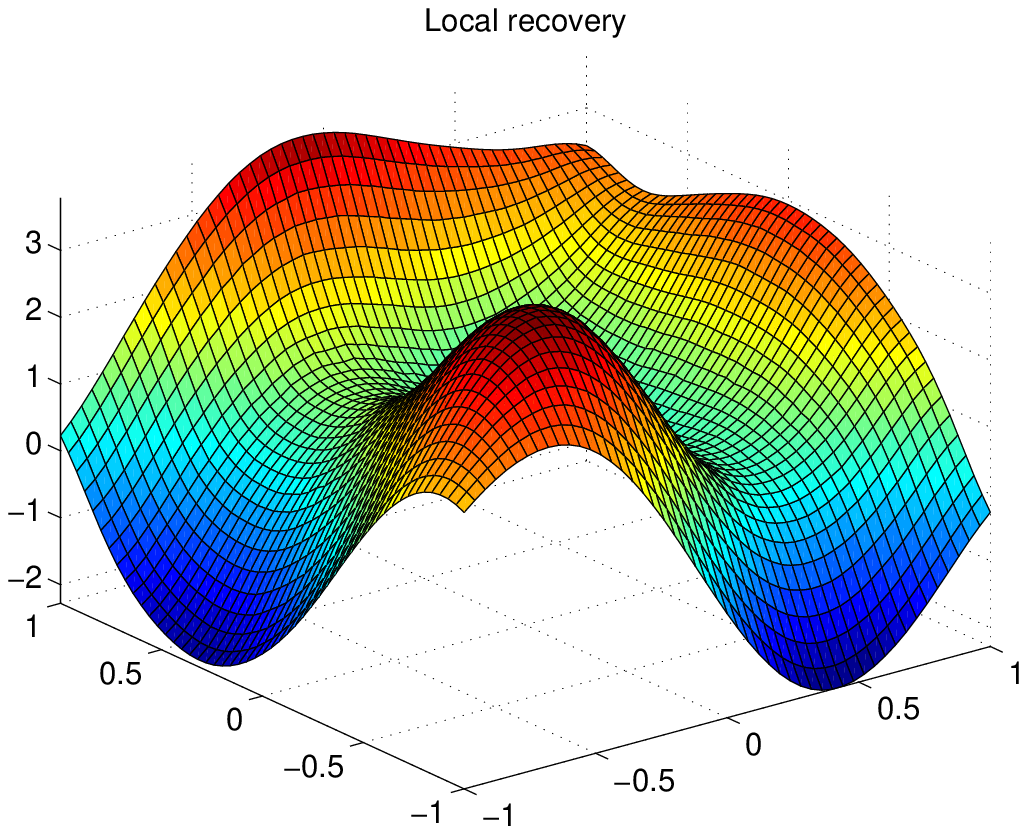}
 \includegraphics[width=\RSw,height=\RSh]{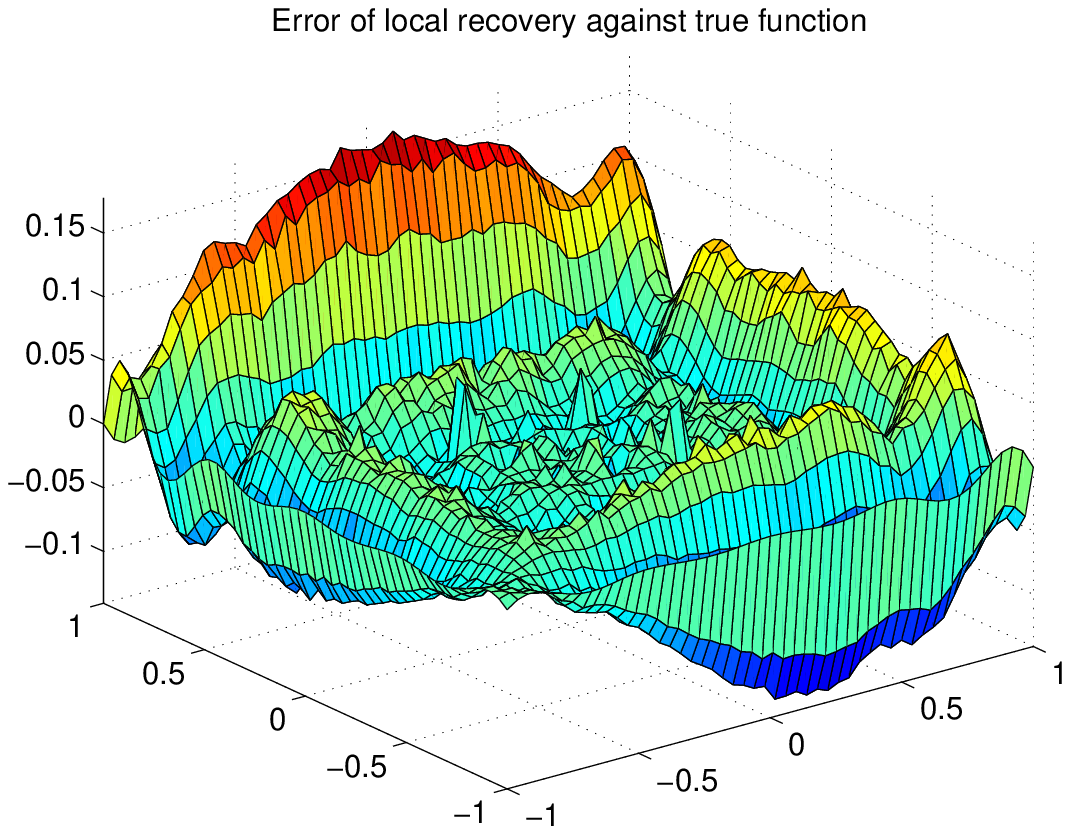}
\end{center} 
    \caption{Local recovery of the {\tt peaks} function  in $W_2^6(\R^2)$ on 2601 regular points in
      $[-1,+1]^2$
      at scale $0.1$ 
      using the greedy point selection
      strategy on 100 regular data points,
      with selection of 21 points out of 100.\RSlabel{FigGlobalm6peaks}}
\end{figure} 
\begin{figure}[hbtp]
    \begin{center}
 \includegraphics[width=\RSw,height=\RSh]{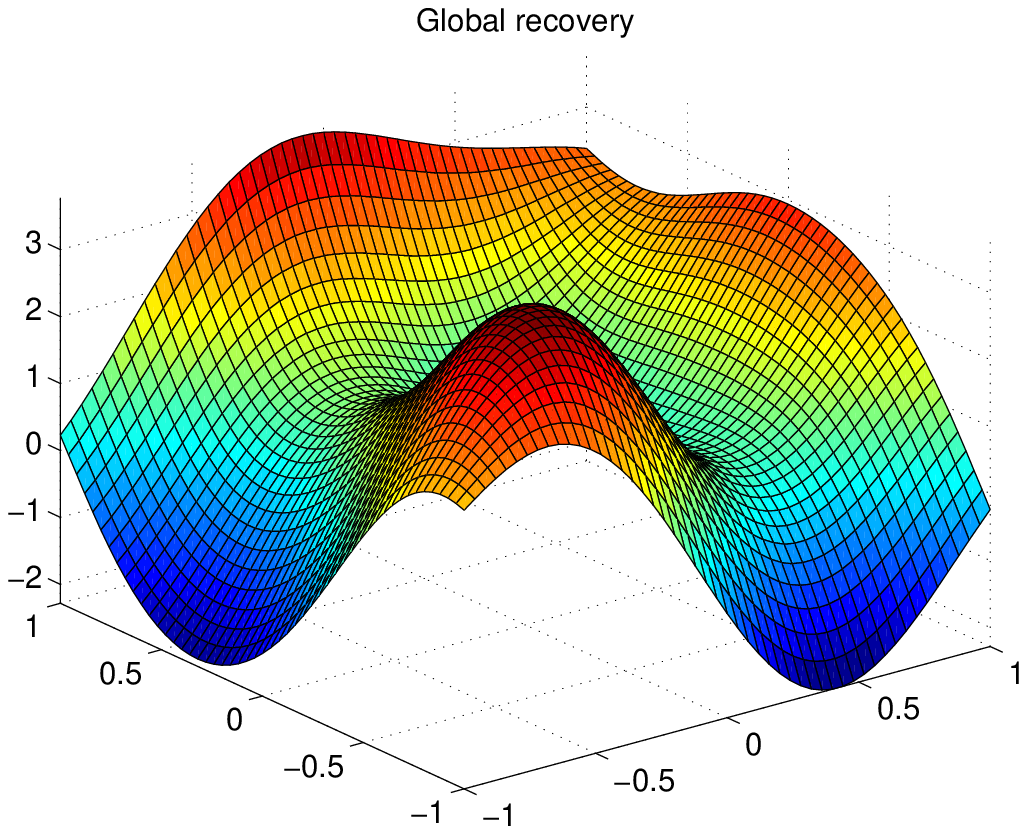}
 \includegraphics[width=\RSw,height=\RSh]{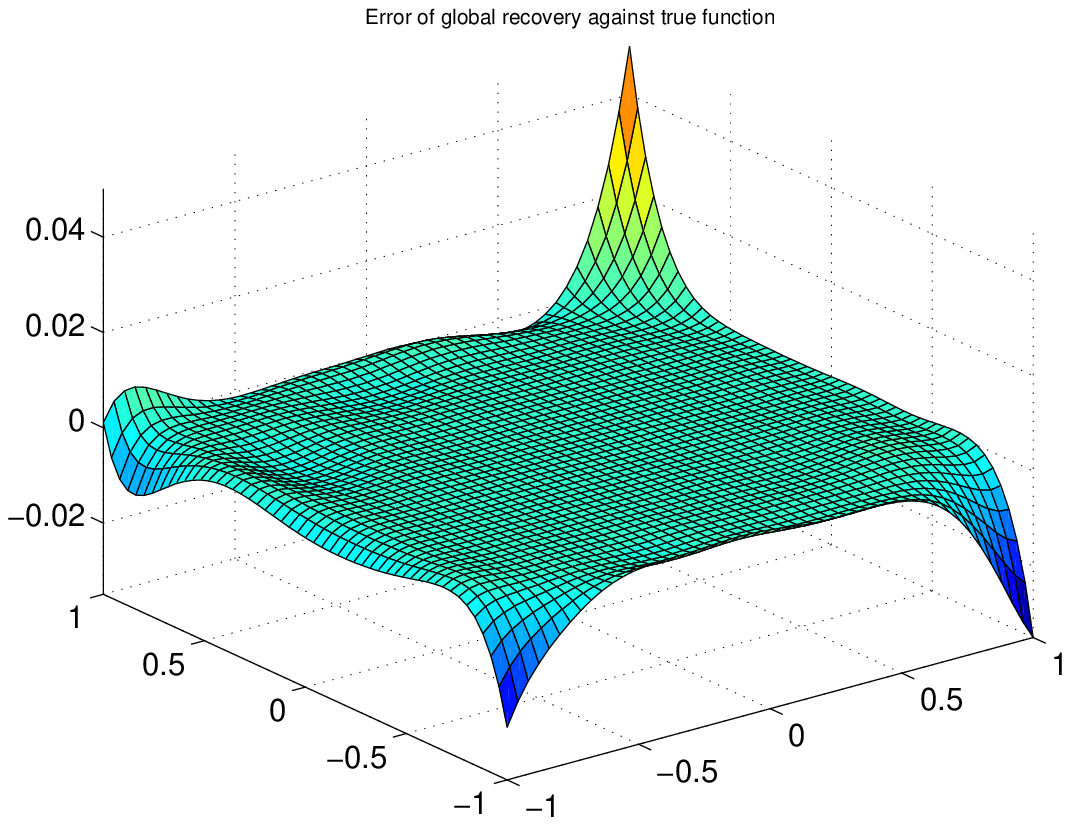}\\
 \includegraphics[width=\RSw,height=\RSh]{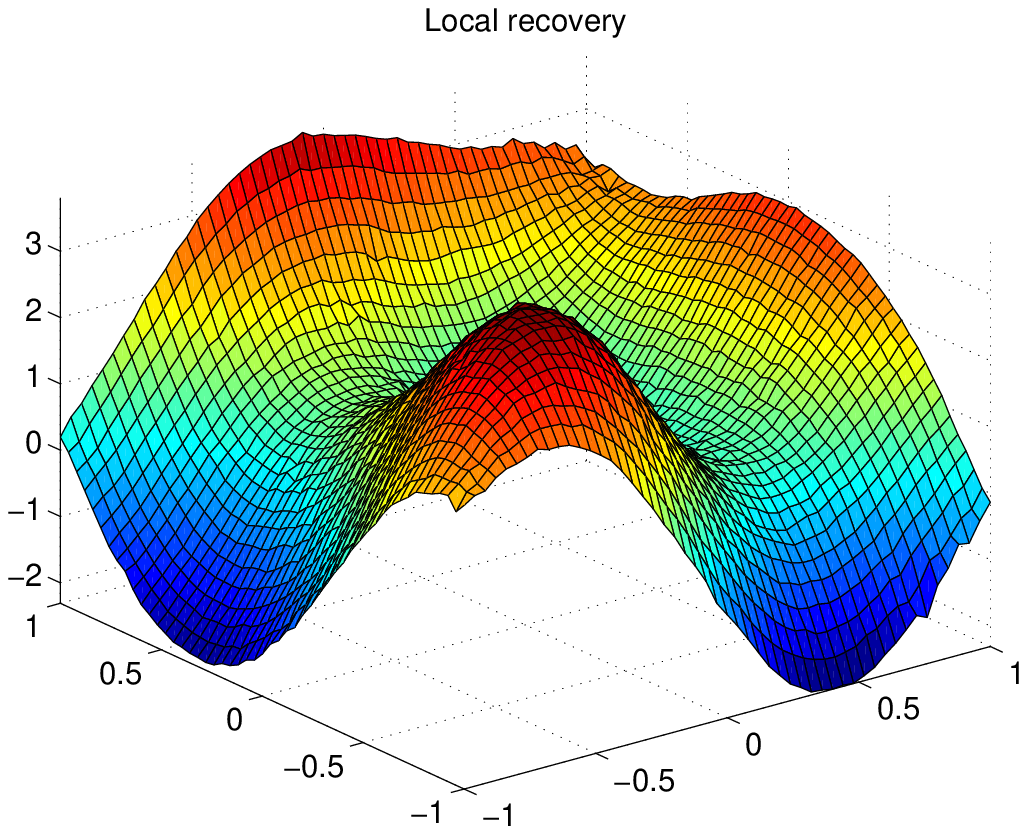}
 \includegraphics[width=\RSw,height=\RSh]{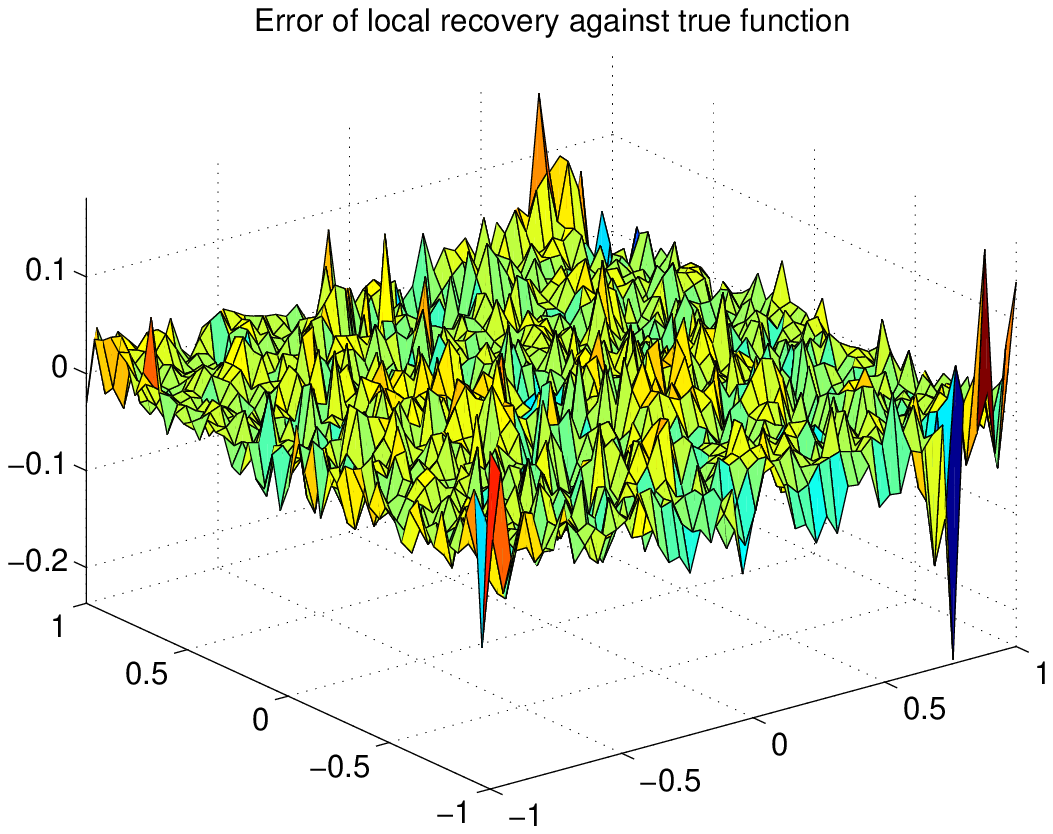}\\
\includegraphics[width=\RSw,height=\RSh]{GlobalPtsUn100-jk0-10m60-thr8-sca1.0-none.eps}
\end{center} 
    \caption{Local recovery of the {\tt peaks} function  in $W_2^6(\R^2)$ on 2601 regular points in
      $[-1,+1]^2$
      at scale $1.0$ 
      using the greedy point selection
      strategy on 100 irregular data points, with a threshold of $10^{-8}$ on the
      squared Power Function. The
      final plot shows the number of points used locally. 
      \RSlabel{FigGlobalm6thrpeaks}}
\end{figure} 
 \biglf
 To see the discontinuities in a close-up, we keep the 100 irregular data points
 in $[-1,+1]^2$, but now we place 2601 regular evaluation points into the
 interval 
 $[0,0.4]^2$
 containing only 13 irregular data points. Note that we can re-evaluate the
 global interpolant on the 2601 local points to have a zoom effect, but the local interpolant has to
 be recalculated, and it gets different and more exact. Figure
 \RSref{FigZoomInm3} shows the case $m=3$ again, to be compared to Figures \RSref{FigCholGlob2m3}
 and 
 \RSref{FigCholGlob2m3peaks}. 
 The top left plot shows the data points (red circles) and the evaluation points
 (blue dots), filling $[0,0.4]$ completely. The squared Power Functions differ
 by a factor of about 2, but are now around $10^{-5}$ instead of $10^{-2}$ in
 Figure \RSref{FigCholGlob2m3}. The second row and the last plot show that the error of the function recovery is about 0.1,
 while it is around 0.5 in Figure \RSref{FigCholGlob2m3peaks}. The Lebesgue
 constants get better, see the lower left plot. 
    \def\RSh{3.5cm}
  \def\RSw{3.5cm}
 
\begin{figure}[hbtp]
    \begin{center}
 \includegraphics[width=\RSw,height=\RSh]{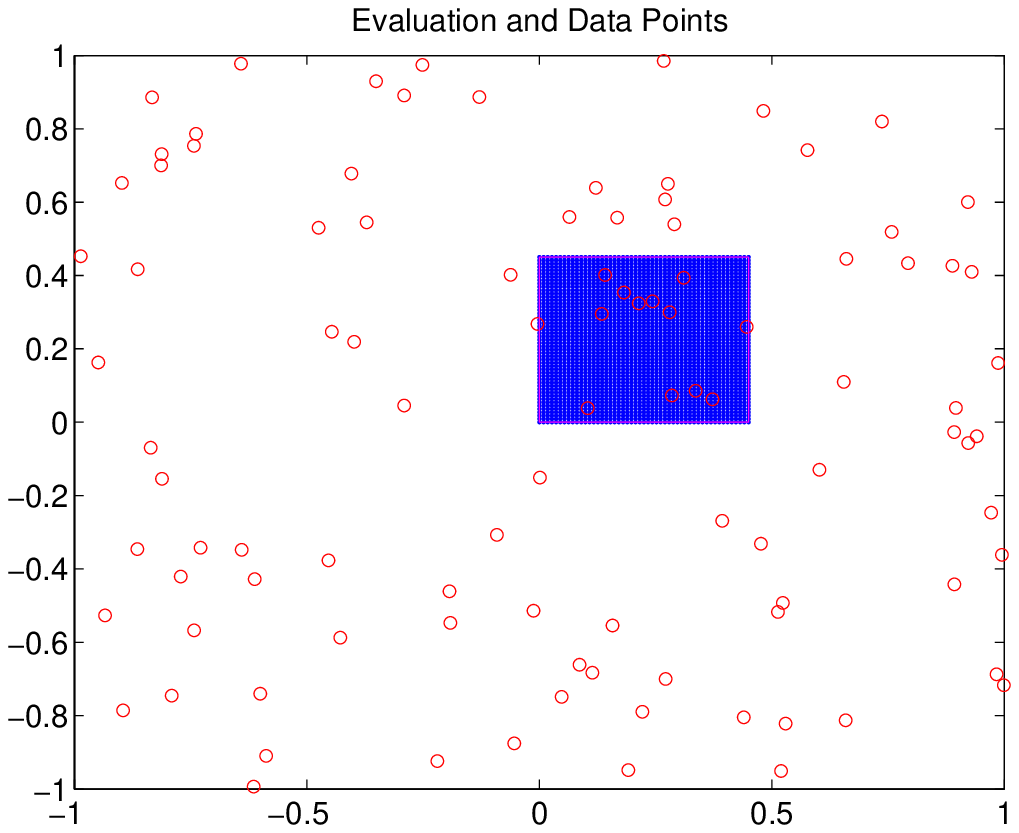}
 \includegraphics[width=\RSw,height=\RSh]{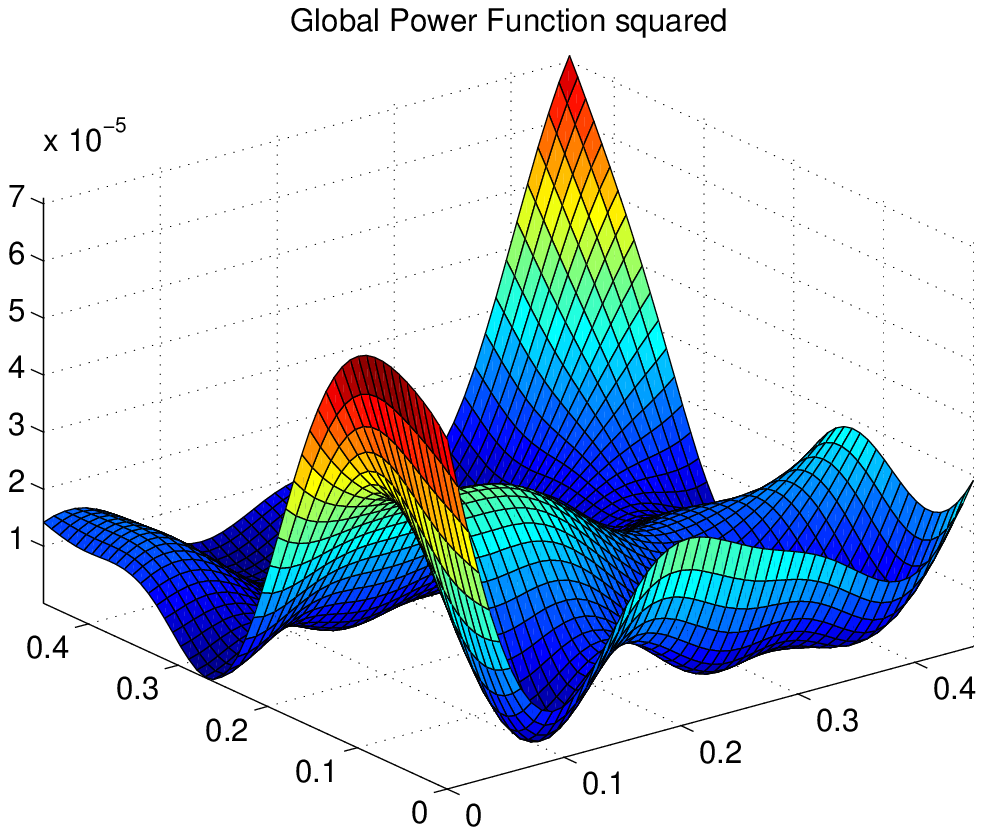}
 \includegraphics[width=\RSw,height=\RSh]{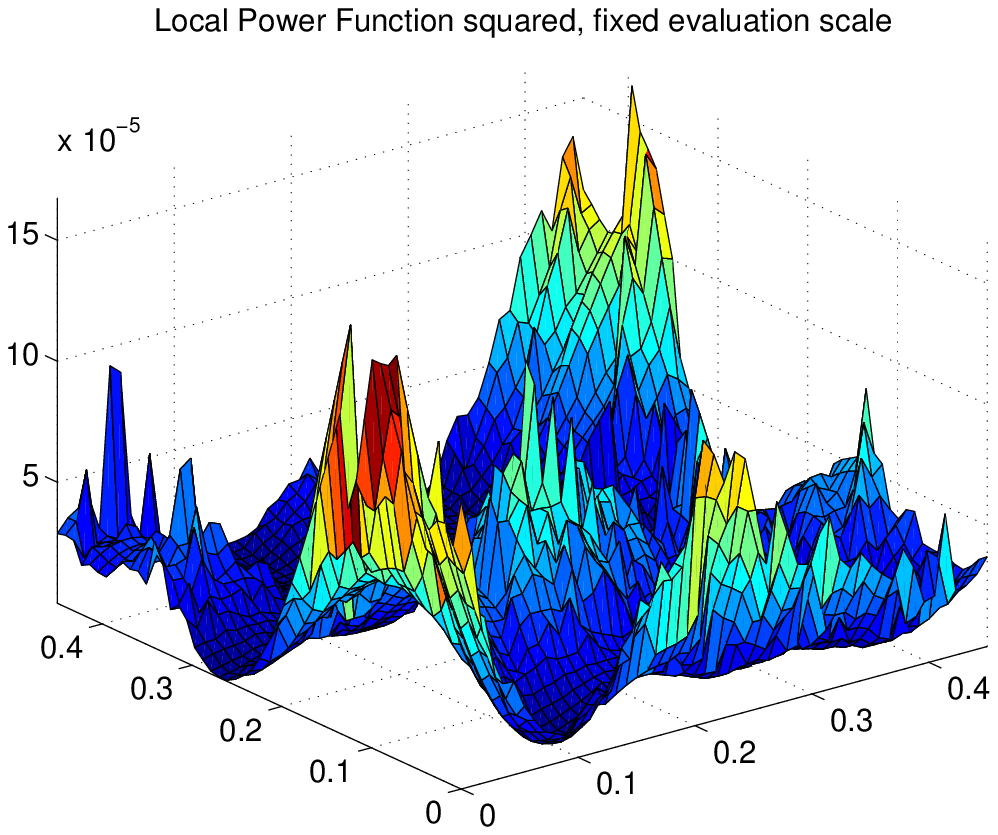}\\
 \includegraphics[width=\RSw,height=\RSh]{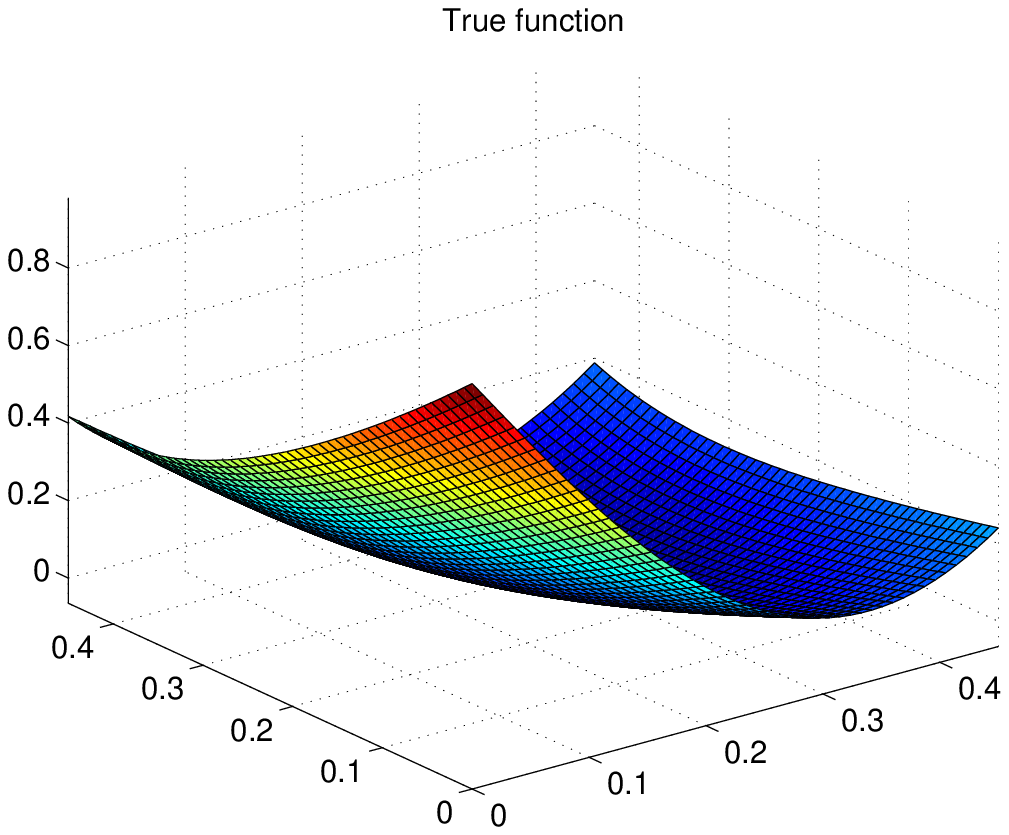}
\includegraphics[width=\RSw,height=\RSh]{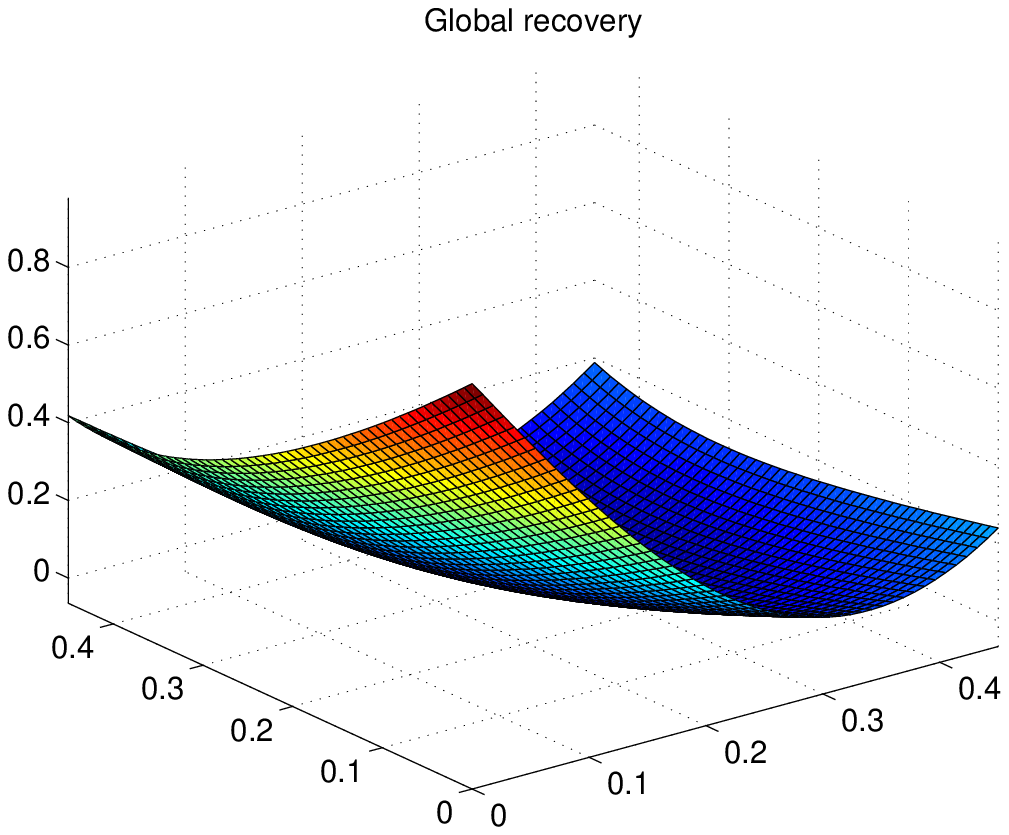}
\includegraphics[width=\RSw,height=\RSh]{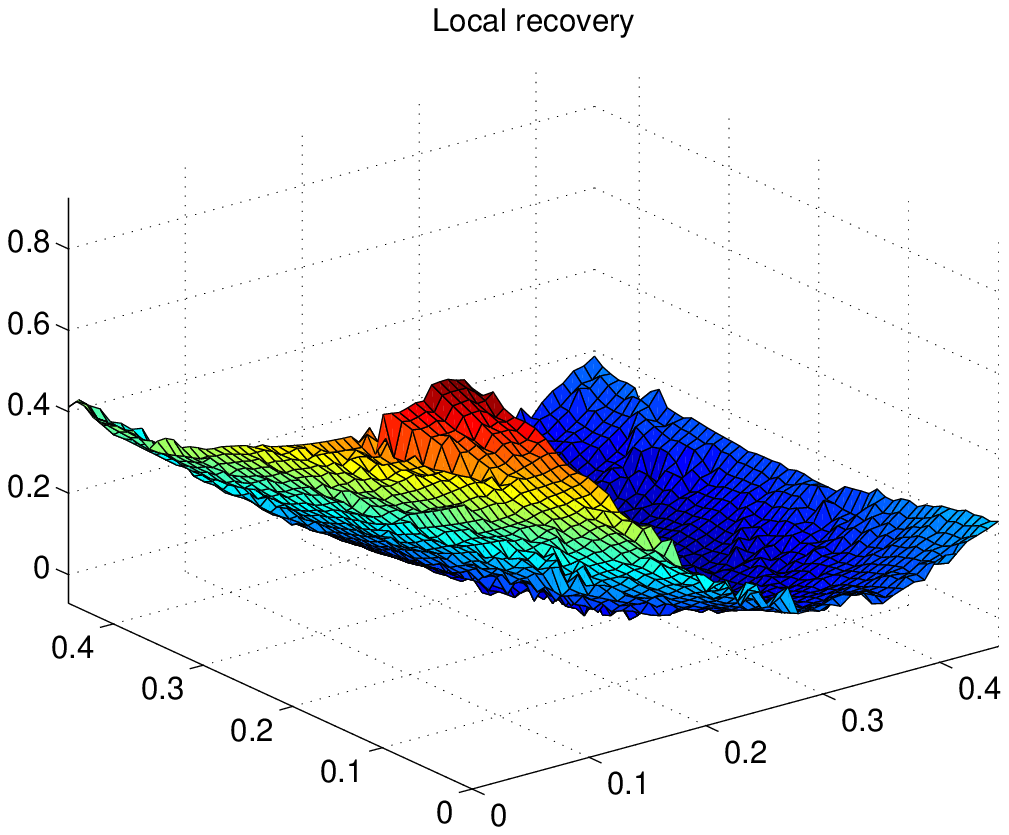}\\
\includegraphics[width=\RSw,height=\RSh]{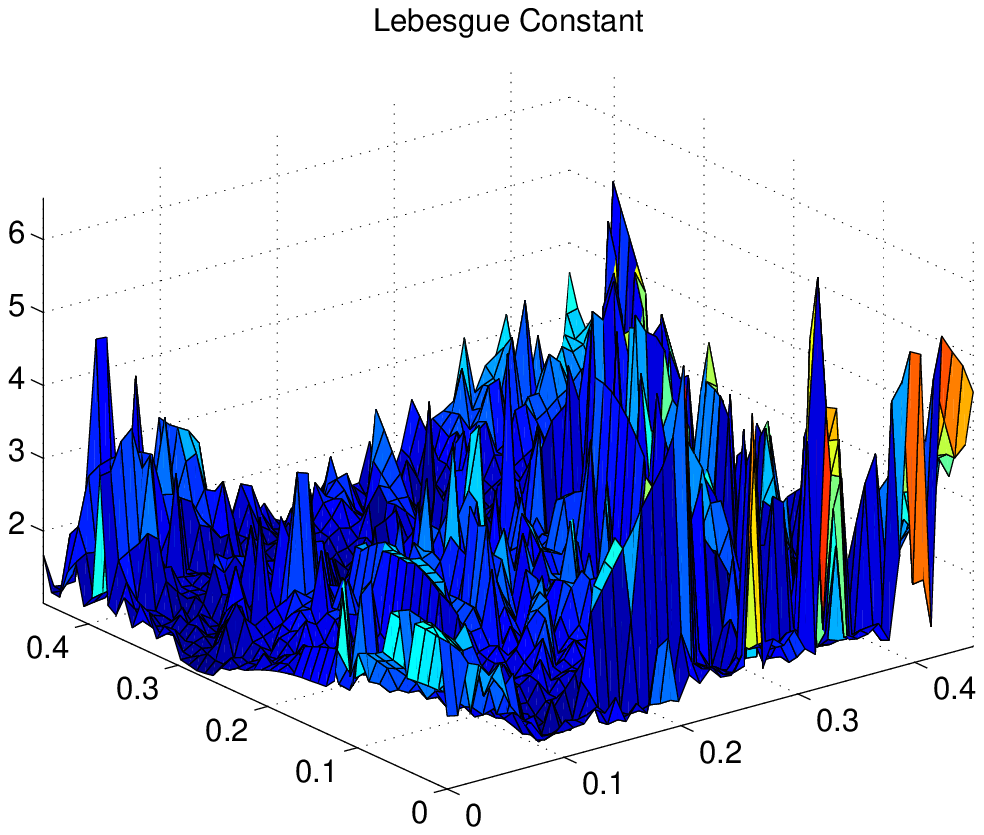}
\includegraphics[width=\RSw,height=\RSh]{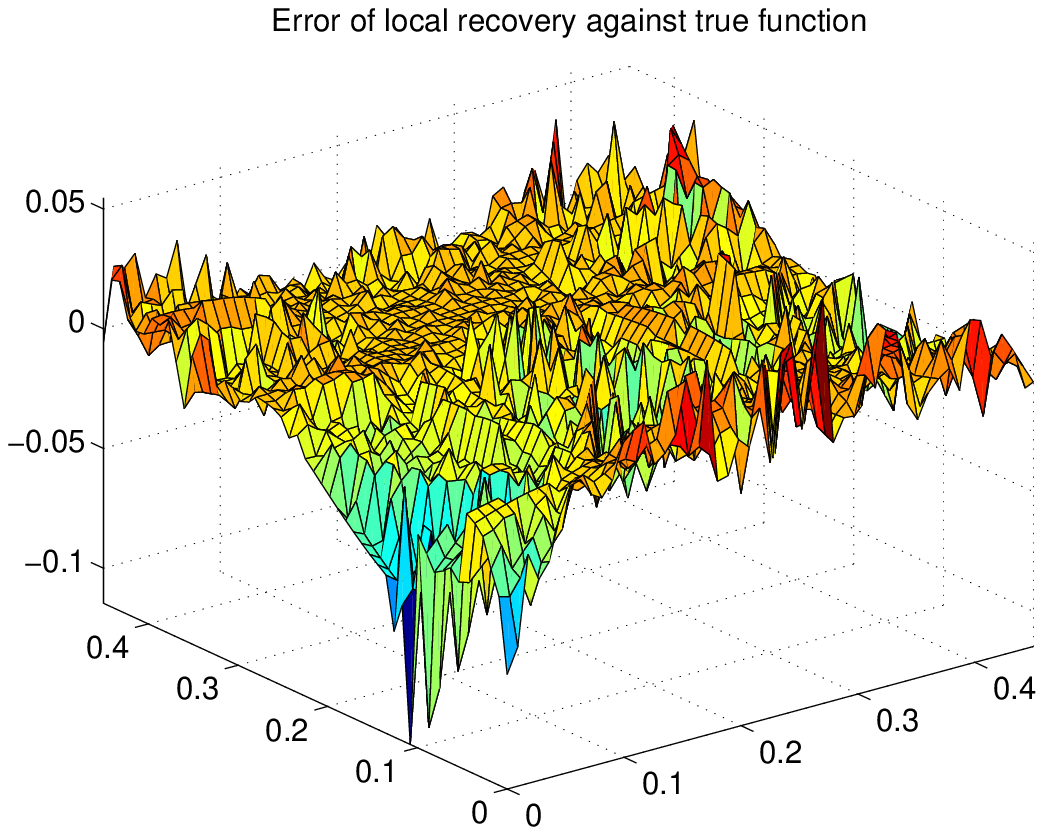}
\end{center} 
    \caption{Local recovery of the {\tt peaks} function  in $W_2^3(\R^2)$ on 2601 regular points in
      $[0,0.4]^2$
      at scale $1.0$ 
      using the greedy point selection
      strategy on 100 irregular data points in $[-1,+1]^2$.
      \RSlabel{FigZoomInm3}}
\end{figure}
\biglf
Figure \RSref{FigZoomInm6} is for comparison to Figures \RSref{FigCholGlob2m6thr}
and \RSref{FigGlobalm6thrpeaks}.
Since the iteration stops at the $10^{-8}$ threshold for the squared Power
Function, the top right plot is chaotic, while the centre plot shows that  using
more points would go down to $10^{-11}$. The actual ploints used are in the
final plot, staying well below 21 proposed for $m=6$. Due to the threshold, the
noise in the right centre plot is roughly constant everywhere, but still better
than in Figure \RSref{FigGlobalm6thrpeaks}.   
\begin{figure}[hbtp]
    \begin{center}
 \includegraphics[width=\RSw,height=\RSh]{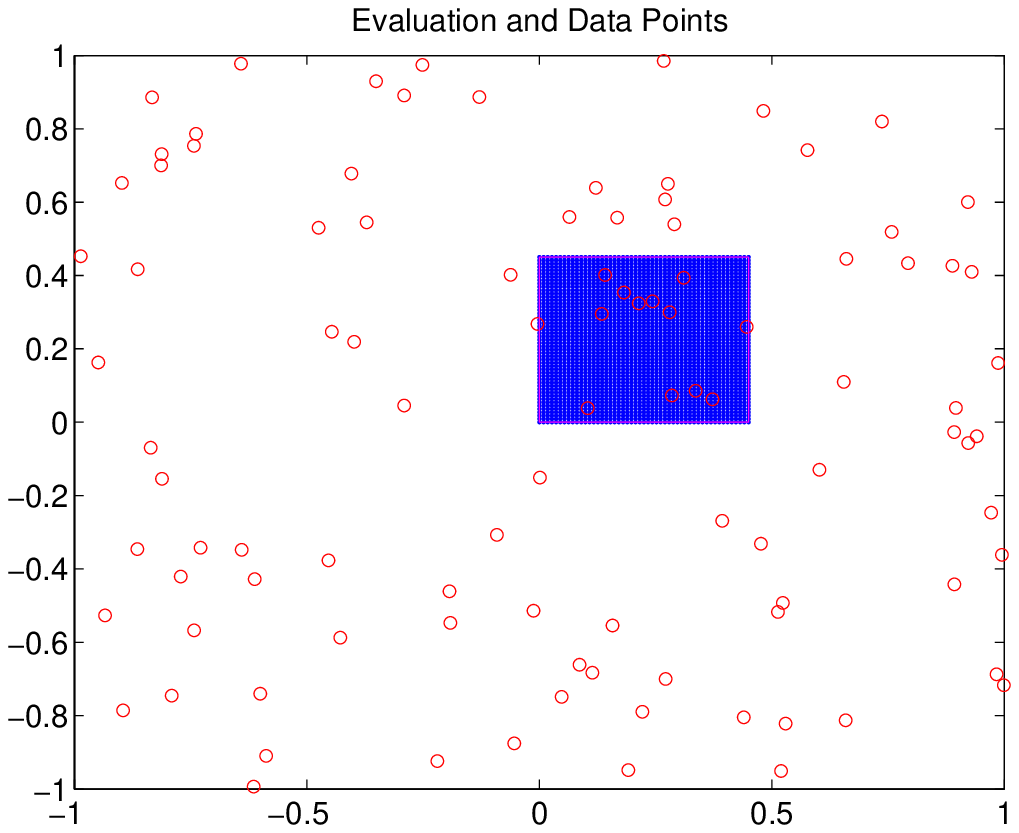}
 \includegraphics[width=\RSw,height=\RSh]{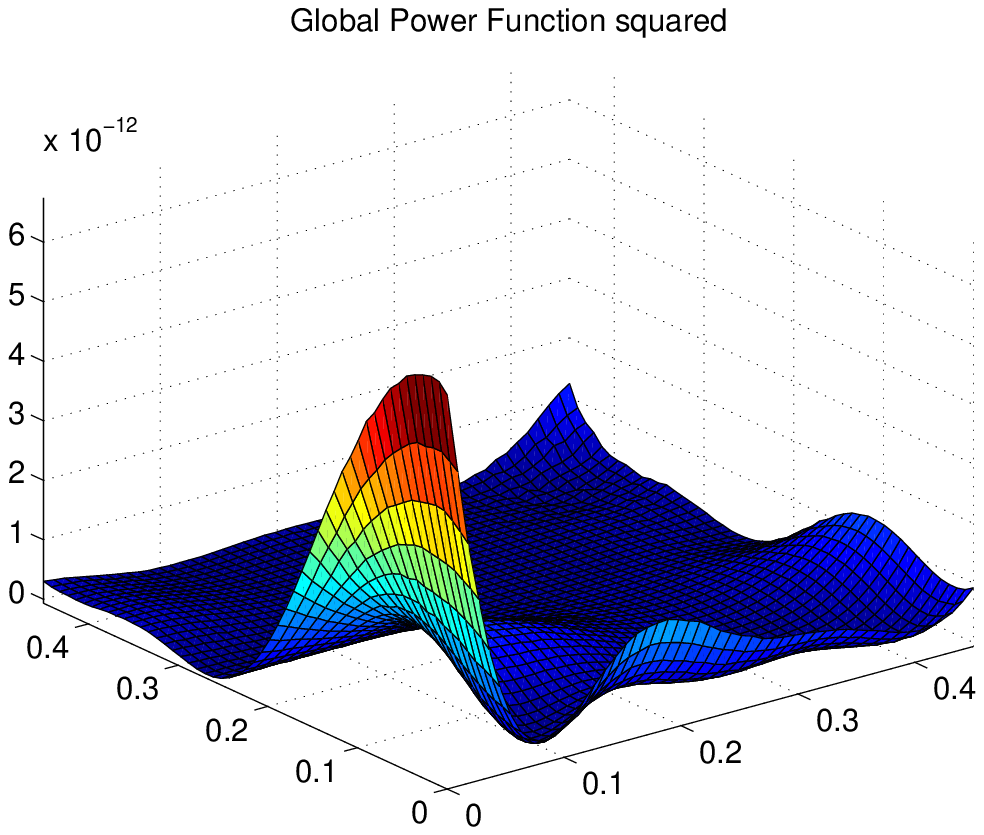}
 \includegraphics[width=\RSw,height=\RSh]{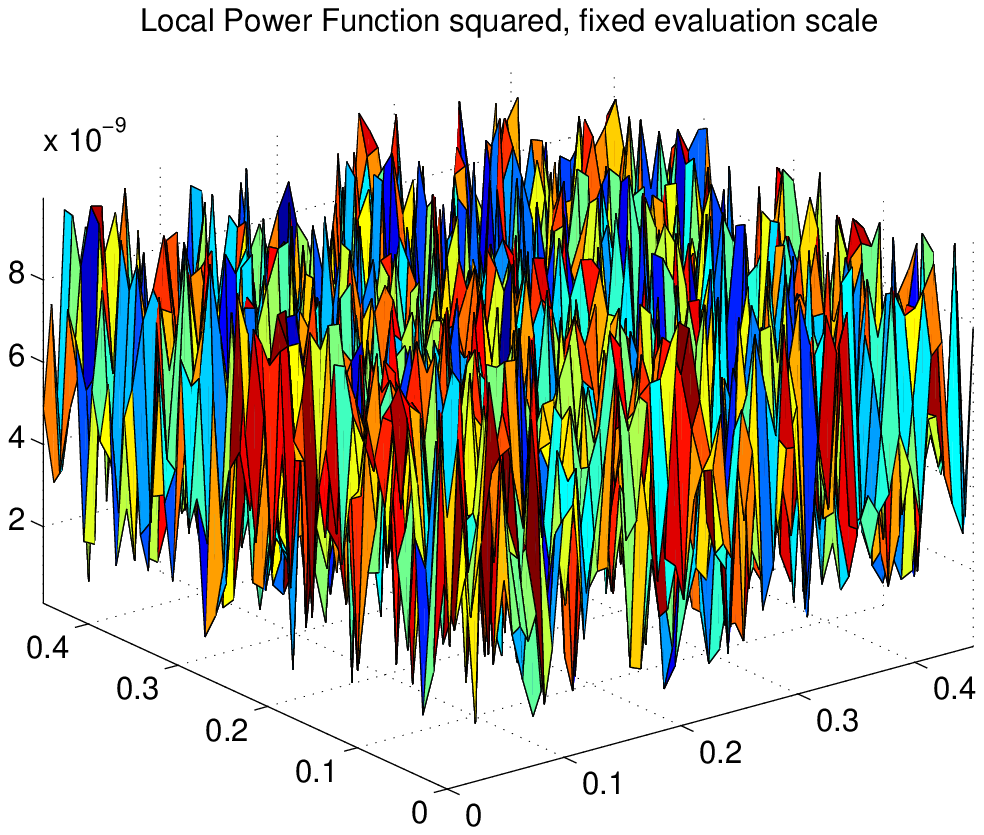}\\
 \includegraphics[width=\RSw,height=\RSh]{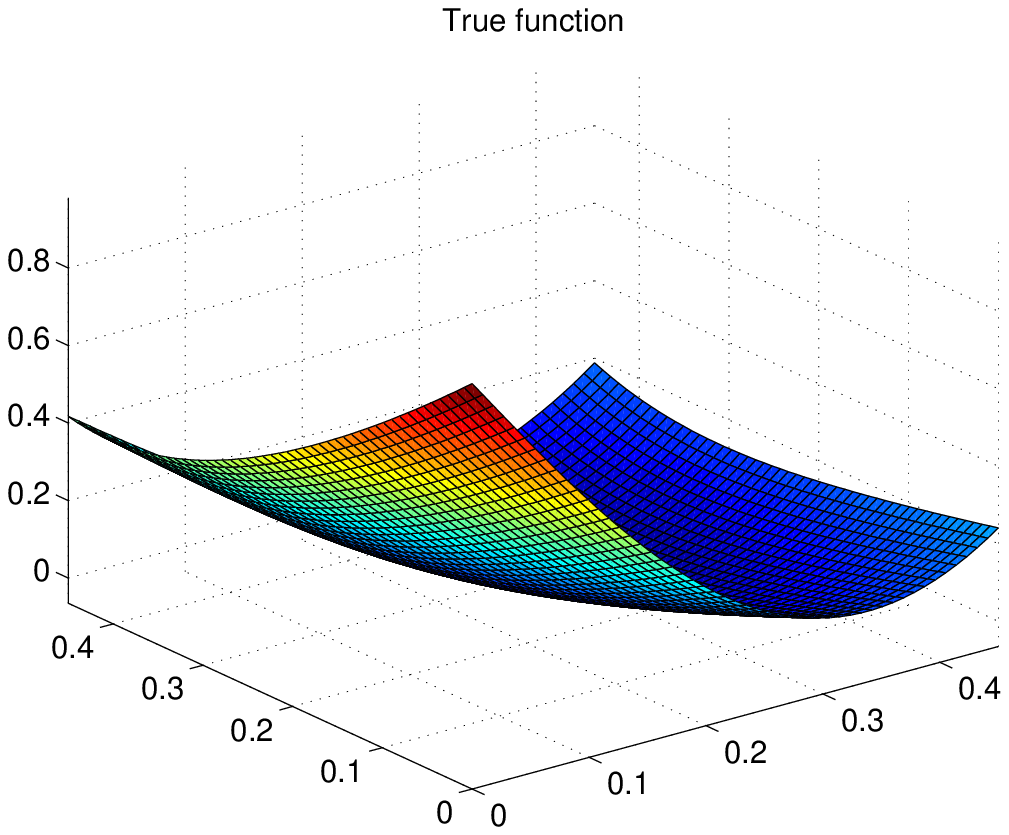}
\includegraphics[width=\RSw,height=\RSh]{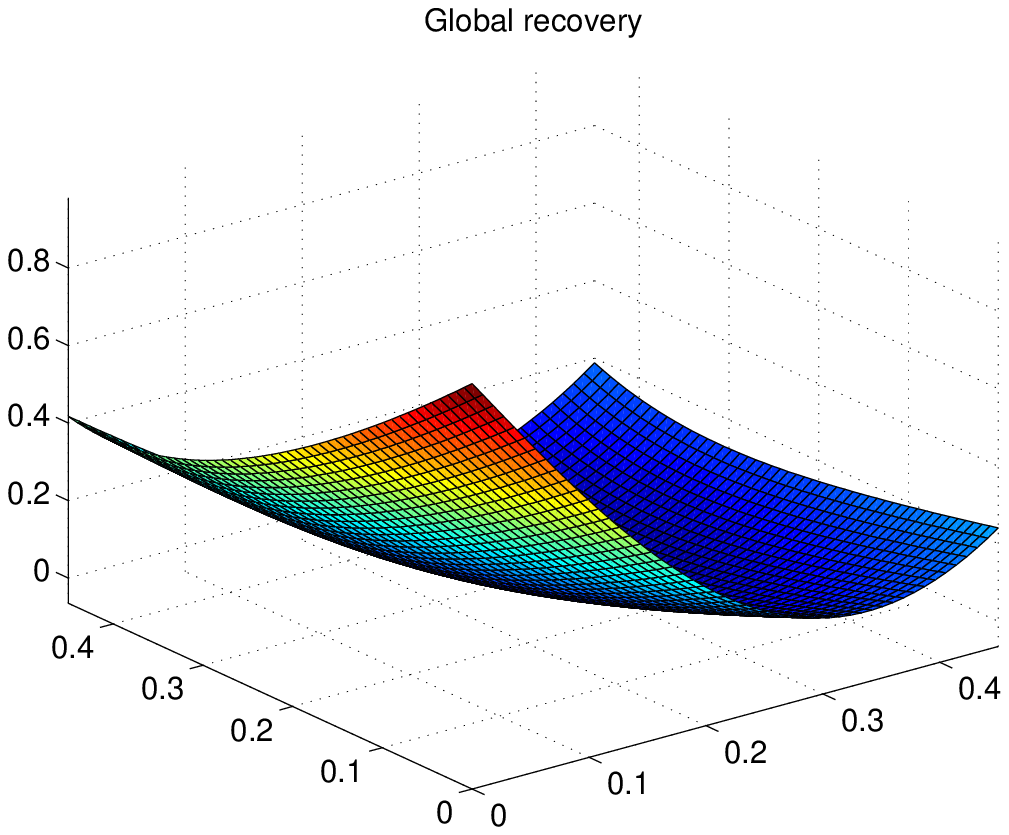}
\includegraphics[width=\RSw,height=\RSh]{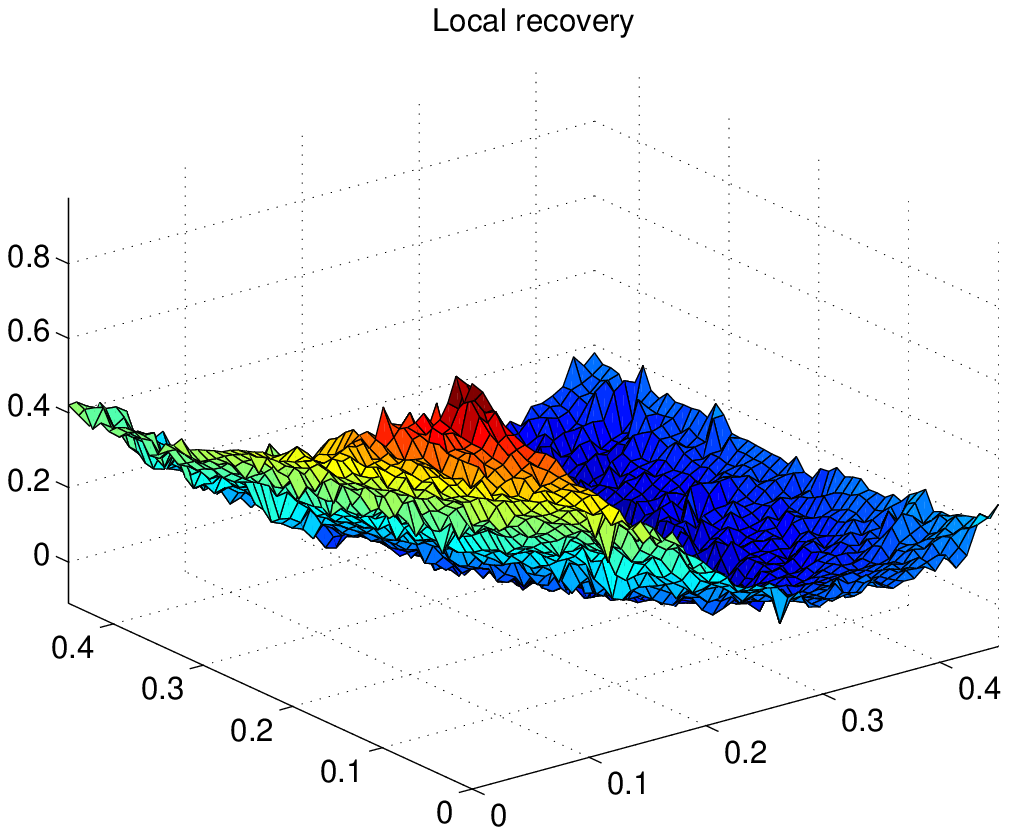}\\
\includegraphics[width=\RSw,height=\RSh]{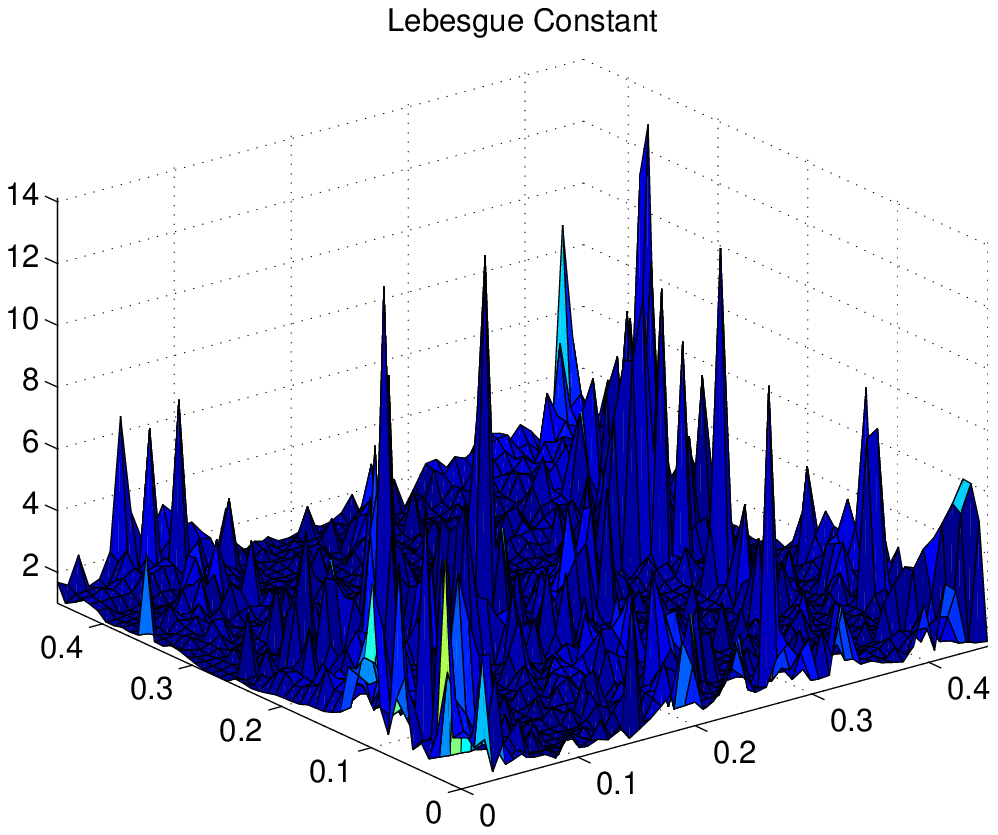}
\includegraphics[width=\RSw,height=\RSh]{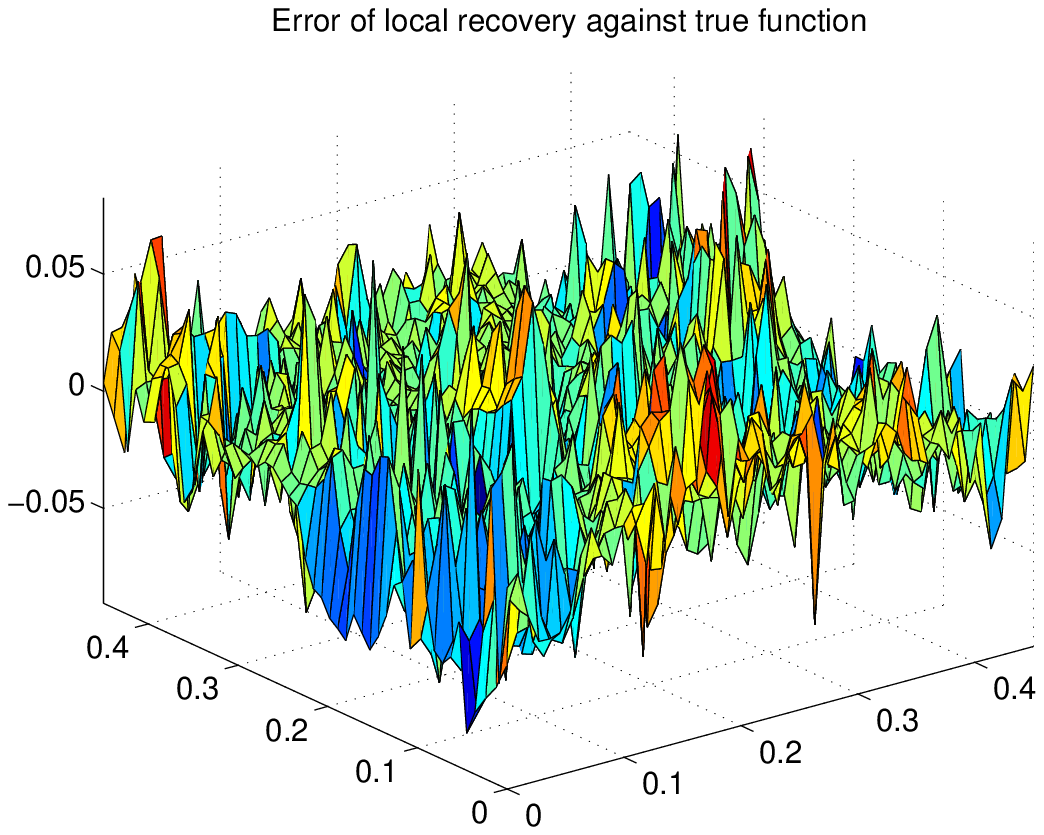}
\includegraphics[width=\RSw,height=\RSh]{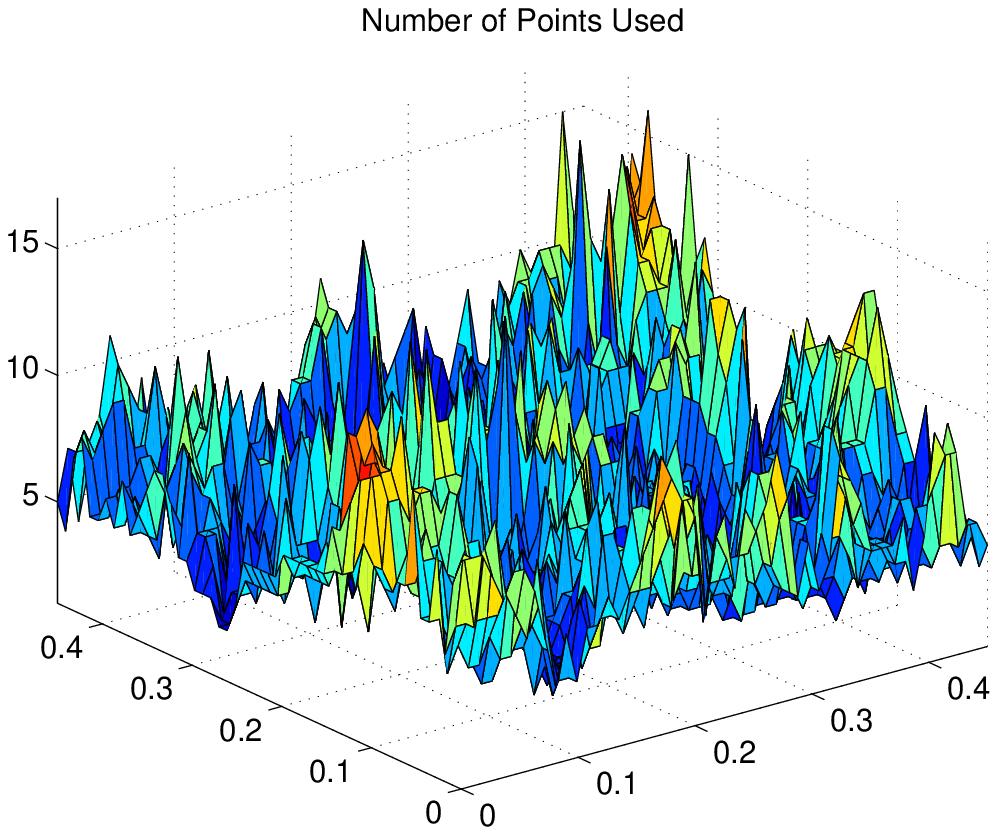}
\end{center} 
    \caption{Local recovery of the {\tt peaks} function  in $W_2^6(\R^2)$ on 2601 regular points in
      $[0,0.4]^2$
      at scale $1.0$ 
      using the greedy point selection
      strategy on 100 irregular data points in $[-1,+1]^2$. The squared Power
      function tolerance is $10^{-8}$.
      \RSlabel{FigZoomInm6}}
\end{figure}
\section{Convergence Rates}\RSlabel{SecConv}
The previous sections tried to achieve $\calO(h^{m-d/2})$ convergence
in $L_\infty$ in Sobolev spaces $W_2^m(\R^d)$ using well-selected
sets of $Q={ q-1+d \choose d}$ points for $q=\lceil m-d/2\rceil$.
This number of points is necessary to let the set be
in general position for polynomials in $\P_q^d$, and this is necessary
for convergence like $h^q$. These rates hold pointwise by construction,
and Figure \RSref{FigRatesHall} shows that they hold in the large.
\biglf
A fixed evaluation grid of $441=21 \times 21$ was used in all cases.
Since the worst situation should occur in corners or near boundaries,
there are no serious changes when using a finer evaluation grid.
The greedy point selection  method was run on each $z$
of the evaluation grid, offering $5\cdot Q$ data points of random sets $X_N$
with $N$ up to 10.000. To avoid additional randomness, the sets $X_N$ were
nested for increasing $N$. The plots show the maximum of the Power Function on
the evaluation set. The expected rates are attained well, see the dotted lines
marking the expected convergence rate. In all cases, the Power Function gets
small enough for letting any function approximation  be exact within
plot precision.
\section{Stability Issues}\RSlabel{SecStab}
The case $m=6$ is unexpectedly unstable already for reasonable $h$,
needing a closer look.
Inspecting the test runs, instabilities come up  when the squared Power Function
reaches machine precision, making the decision \eref{eqdecision} unsafe.
This occurs even when the local triangular
$Q\times Q$ matrix of \eref{eqNewtonrep} is still
within standard condition limits. 
The squared Power Function behaves like $h^{2m-d}$ in $L_\infty$
by standard convergence theory, and then fill distances of order
$$
h< 10^{-\frac{15}{2m-d}}=:\overline{h}_{2m-d}
$$
can be expected to fail in double precision. Of course, scaling
and multipliers plays a major part here,
but the ratios
$$
\overline{h}_{2}\;:\;
\overline{h}_{5}\;:\;
\overline{h}_{10} \;\;\;\hbox{    like     }\;\;\;
0.0000000316\;:\;
0.001 \;:\;
0.0316 \;;
$$
matching $m=1.5,\,3$, and $6$ in $\R^2$
indicate that large $m$ will fail unexpectedly early. Figure
\RSref{FigRatesHall} goes down to
10.000 points with $h=0.02$, illustrating the above argument. It can be expected
that the case $m=3/2$ runs up to using approximately 60 million points. 
\biglf
Note that this applies also to global interpolation, making small $h$ for
large $m$ hazardous, if not impossible. But the condition of the global
kernel matrix does not enter here, in contrast to the standard stability
theory of kernel-based interpolation.
The well-known phenomena like rank loss or bad condition
are local, not global. Like in \RScite{schaback:1994-2}, they are
connected to how well a kernel can be approximated by polynomials.
\biglf
To circumvent the stability problems,
the additional green line
shows the results when the threshold for the Power Function is set to
$10^{-6}$.

\begin{figure}[hbtp]
    \begin{center}
      \includegraphics[width=6.0cm,height=6.0cm]{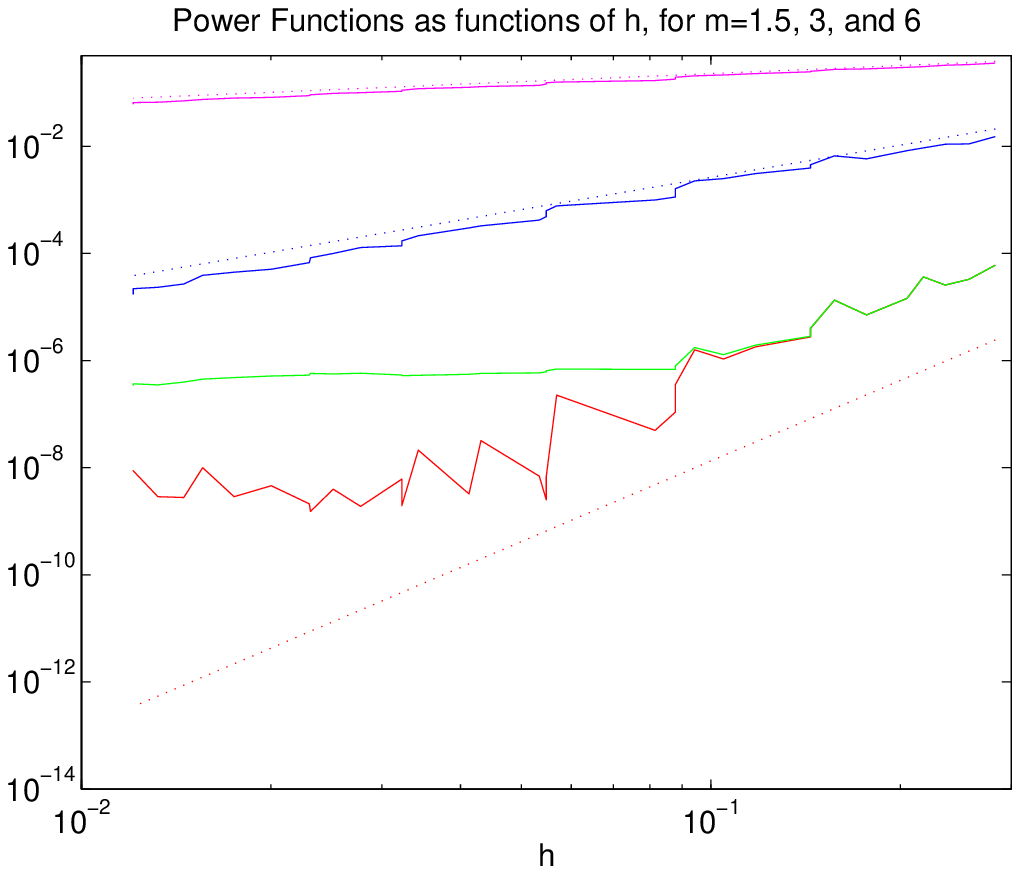}
\end{center} 
    \caption{Power Functions as functions of $h$ for $m=1.5$ (magenta), $3$
      (blue), and $6$ (red)
      based on minimal point sets. The dotted lines
      indicate the expected theoretical rates $h^{m-d/2}$.
      \RSlabel{FigRatesHall}}
\end{figure}

\section{Conclusions and Open Problems}\RSlabel{SecCOP}
The suggested selection algorithm for local point sets works as
expected, leading to optimal convergence rates in Sobolev spaces
using the minimal possible number of points for that purpose.
Resulting functions are discontinuous, but since they converge in $L_\infty$ to
smooth functions, they can be called {\em asymptotically smooth}.
\biglf
However, the method can run into instabilities for large data sets
roughly at the limits of stability of the global problem.
Since the method does not produce {\em scalable stencils}
in the sense of \RScite{davydov-schaback:2019-1}, the instability problems
may be overcome by going to scale-invariant techniques. This needs further work,
along with comparisons to Moving Least Squares or Shepard-type techniques.
\biglf
Furthermore, it is open whether the technique always produces
sets that are in general position with respect to polynomials,
and how the Lebesgue constants behave. In contrast to global techniques,
the method can
vary the kernel scale locally without sacrificing local convergence rates.
This needs further work as well. 
\biglf
The method can easily be extended to finding
good point sets for local approximations of derivatives
\cite{wu:1992-1,davydov-schaback:2018-2,
  davydov-schaback:2019-1}, and this
will have some influence on meshless RBF-FD methods 
\cite{wright-fornberg:2006-1,
fornberg-et-al:2013-1,
larsson-et-al:2013-1,barnett:2015-1,
davydov-schaback:2018-2,flyer-et-al:2016-1}. 
\biglf
Also, the implications for the ``flat limit'' situation,
see
e.g. \cite{driscoll-fornberg:2002-1,fornberg-et-al:2004-1,%
  larsson-fornberg:2005-1,schaback:2005-2,schaback:2008-2,song-et-al:2012,
fasshauer-mccourt:2015-1} are worth investigating.
\biglf
There are no conflicts of interest.
\bibliographystyle{plain}

\end{document}